\documentclass[leqno,11pt]{article}

\usepackage{amssymb}
\usepackage{amsmath}
\usepackage{amsthm}
\usepackage[dvipdfmx]{graphicx}
\usepackage{bm}
\usepackage{cases}

\makeatletter
    
    \@addtoreset{equation}{section}
\makeatother

 \newtheorem{thm}{Theorem}[section]
 
 \newtheorem{lem}[thm]{Lemma}
 \newtheorem{prop}[thm]{Proposition}

\begin{document}

\begin{center}
\bf\Large
Uniform Estimates for the Flow of a Viscous Incompressible Fluid
 down an Inclined Plane in the Thin Film Regime
\end{center}

\begin{center}
{\em
Dedicated to Professor Shuichi Kawashima on the occasion of his 60th birthday
}
\end{center}

\begin{center}
Hiroki Ueno, Akinori Shiraishi, and Tatsuo Iguchi
\end{center}

\begin{abstract}
We consider a two-dimensional motion of a thin film flowing down an inclined plane 
under the influence of the gravity and the surface tension. 
In order to investigate the stability of such flow, it is hard to treat the Navier--Stokes equations directly, 
so that a thin film approximation is often used. 
It is an approximation obtained by the perturbation expansion with respect to the aspect ratio $\delta$ 
of the film under the thin film regime $\delta\ll1$. 
Our purpose is to give a mathematically rigorous justification of the thin film approximation 
by establishing an error estimate between the solution of the Navier--Stokes equations and 
those of approximate equations. 
To this end, in this paper we derive a uniform estimate for the solution of the Navier--Stokes equations 
with respect to $\delta$ under appropriate assumptions.
\end{abstract}

\section{Introduction}
In this paper, we consider a two-dimensional motion of a liquid film of a viscous and incompressible fluid 
flowing down an inclined plane under the influence of the gravity and the surface tension on the interface. 
The motion can be mathematically formulated as a free boundary problem for the incompressible Navier--Stokes equations. 
We assume that the domain $\Omega(t)$ occupied by the liquid at time $t\ge0$, the liquid surface $\Gamma(t)$, 
and the rigid plane $\Sigma$ are of the forms 
$$
\left\{
\begin{array}{l} 
\Omega (t)=\{(x,y)\in\mathbb{R}^2\ |\ 0<y<h_0+\eta(x,t)\}, \\
\Gamma(t)=\{(x,y)\in\mathbb{R}^2\ |\ y=h_0+\eta(x,t)\}, \\
\Sigma=\{(x,y)\in\mathbb{R}^2|\ y=0\},
\end{array}
\right.
$$
where $h_0$ is the mean thickness of the liquid film and $\eta(x,t)$ is the amplitude of the liquid surface. 
Here we choose a coordinate system $(x, y)$ so that $x$ axis is down and $y$ axis is normal to the plane. 
\bigskip
\begin{figure}[h]
\begin{center}
\includegraphics{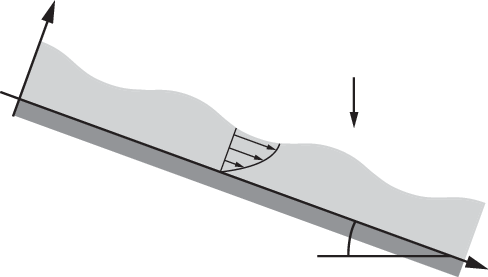}
\end{center}
\caption{Sketch of a thin liquid film flowing down an inclined plane}
\end{figure}
\setlength{\unitlength}{1pt}
\begin{picture}(0,0)
\put(300,-165){$x$}
\put(95,-30){$y$}
\put(110,-60){$\Gamma(t)$}
\put(135,-90){$\Omega(t)$}
\put(100,-100){$\Sigma$}
\put(240,-80){$\bm{g}$}
\put(220,-145){$\alpha$}
\end{picture}
The motion of the liquid is described by the velocity $\bm{u}=(u,v)^{\rm T}$ and 
the pressure $p$ satisfying the Navier--Stokes equations 
\begin{equation}\label{ns0}
\left\{
 \begin{array}{lcl}
  \rho\bigl(\bm{u}_t+(\bm{u}\cdot\nabla)\bm{u}\bigr)=\nabla\cdot\mathbf{P}+\rho g(\sin\alpha, -\cos\alpha)^{\rm T}
   & \mbox{in} & \Omega(t), \; t>0, \\
  \nabla\cdot\bm{u}=0 & \mbox{in} & \Omega(t), \; t>0,
 \end{array}
\right.
\end{equation}
where $\mathbf{P}=-p\mathbf{I}+2\mu\mathbf{D}$ is the stress tensor, 
$\mathbf{D}=\frac{1}{2}\bigl(D\bm{u}+(D\bm{u})^{\rm T}\bigr)$ is the deformation tensor, 
$\mathbf{I}$ is the unit matrix, 
$\rho$ is a constant density of the liquid, $g$ is the acceleration of the gravity, 
$\alpha$ is the angle of inclination, and $\mu$ is the shear viscosity coefficient. 
The dynamical and kinematic conditions on the liquid surface are 
\begin{equation}\label{1bc0}
\left\{
 \begin{array}{lcl}
  \mathbf{P}\bm{n}=-p_0\bm{n}+\sigma H\bm{n} & \mbox{on} & \Gamma(t), \; t>0, \\
  \eta_t+u\eta_x-v=0 & \mbox{on} & \Gamma(t), \; t>0,
 \end{array}
 \right.
\end{equation}
where $\bm{n}$ is the unit outward normal vector to the liquid surface, that is, 
$\bm{n}=\frac{1}{\sqrt{1+\eta_{x}}}(-\eta_x,1)^{\rm T}$, 
$p_0$ is a constant atmospheric pressure, $\sigma$ is the surface tension coefficient, 
and $H$ is the twice mean curvature of the liquid surface, that is, 
$H=\Bigl(\frac{\eta_x}{\sqrt{1+\eta_{x}^{2}}}\Bigr)_x$. 
The boundary condition on the rigid plane is the non-slip condition 
\begin{equation}\label{0bc0}
\bm{u}=\bm{0} \quad\mbox{on}\quad \Sigma, \; t>0.
\end{equation}

\medskip
These equations have a laminar steady solution of the form 
\begin{equation}\label{laminar flow}
\eta=0, \quad u=(\rho g\sin\alpha/2\mu)(2h_0y-y^2), \quad v=0, \quad p=p_0-\rho g\cos\alpha(y-h_0),
\end{equation}
which is called the Nusselt flat film solution. Throughout this paper, we assume that the flow is downward $l_0$-periodic or approaches asymptotically this flat film solution at spacially infinity.

Concerning the instability of this laminar flow, there are vast research literatures from the physical 
and  engineering point of view. 
The first investigation of the wave motion of thin film including the effect of the surface tension 
was provided by Kapitza \cite{Kapitza}. 
In Particular, he considered the case where the liquid film flows down a vertical wall, 
that is, the case $\alpha=\frac{\pi}{2}$. 
Yih \cite{Yih1} first formulated the linear stability problem of the laminar flow of the liquid film 
flowing down an inclined plane as an eigenvalue problem for the complex phase velocity, 
more specifically, the Orr-Sommerfeld problem although he neglected the effect of the surface tension. 
Benjamin \cite{Benjamin} took into account the effect of the surface tension and showed that 
the critical Reynolds number is given by ${\rm R}_c=\frac{5}{4}\frac1{\tan\alpha}$  by expanding the normal mode solution in powers of $y$. 
(In his original paper, the critical Reynolds number was given by ${\rm R}_c=\frac{5}{6}\frac1{\tan\alpha}$. 
This difference comes from the definition of the Reynolds number, that is, Benjamin used the average 
speed of the Nusselt flat film solution, whereas we use the speed of the solution on the 
liquid surface as in Benney \cite{Benney}.) 
Later, Yih \cite{Yih2} showed the same condition by expanding the normal mode solution in powers of 
the aspect ratio of the film which will be denoted by $\delta$ in this article. 
An approach taking into account the nonlinearity was first given by Mei \cite{Mei} and Benney \cite{Benney}. 
While Mei considered the gravity waves, Benney considered the capillary-gravity waves and he recovered 
Benjamin's and Yih's linear stability theories. 

Using the mean thickness of the liquid $h_0$, the characteristic scale of the streamwise direction $l_0$, 
and the typical amplitude of the liquid surface $a_0$, Benney introduced two non-dimensional parameters 
$\delta$ and $\varepsilon$ defined by 
\[
\delta=\frac{h_0}{l_0},\quad\varepsilon=\frac{a_0}{h_0},
\]
respectively.
We note that we do not determine  a characteristic scale $l_0$ in $x$ a priori because $l_0$ is a typical wavelength of a nontrivial wave pattern which arises as a consequence of a destabilization and  $l_0$ itself is an object of scientific interest.
While theoretically the destabilization appears as a long wave instability in the case $\delta\to0$, which corresponds to  the case $l_0\to\infty$,  experimentally $l_0$ is often determined by observing waves generated by an external vibrator.
Benney derived the following single nonlinear evolution equation 
\begin{align}\label{Benney}
\eta_t=
& A(1+\varepsilon\eta)\eta_x + \delta\big(B(1+\varepsilon\eta)\eta_{xx}
  + \varepsilon C(1+\varepsilon\eta)\eta_x^2\big) \\
& + \delta^2\big(D(1+\varepsilon\eta)\eta_{xxx} + \varepsilon E(1+\varepsilon\eta)\eta_x\eta_{xx}
  + \varepsilon^2 F(1+\varepsilon\eta)\eta_x^3\big) \notag \\
& + \delta^3\big(G(1+\varepsilon\eta)\eta_{xxxx} + \varepsilon H(1+\varepsilon\eta)\eta_x\eta_{xxx}
  + \varepsilon I(1+\varepsilon\eta)\eta_{xx}^2 \notag \\
& \phantom{+\delta^2\big(}+\varepsilon^2J(1+\varepsilon\eta)\eta_x^2\eta_{xx}
  + \varepsilon^3K(1+\varepsilon\eta)\eta_x^4\big)+O(\delta^4) \notag
\end{align}
with polynomials $A,B,\ldots,K$ in $1+\varepsilon\eta$ by the method of a perturbation expansion of the solution $(u, v, p)$ with respect to $\delta$ 
under the thin film regime $\delta\ll1$.

Thereafter, several authors have followed the Benney's approach. 
We note that if the Weber number ${\rm W}$ satisfies the condition ${\rm W}=O(1)$, 
the effect of the surface tension does not appear until the term of $O(\delta^3)$ in \eqref{Benney}. 
Since Benney considered the case ${\rm W}=O(1)$ and calculated the terms up to $O(\delta^2)$, 
the effect of the surface tension was omitted in his stability analysis. 
Consequently, his results showed that linearly unstable waves grow more rapidly in the nonlinear range. 
Nakaya \cite{Nakaya} computed the terms up to $O(\delta^3)$ and showed that the surface tension has 
a stabilization effect in the development of the monochromatic waves. 
On the other hand, Gjevik \cite{Gjevik} incorporated the effect of the surface tension into the equation 
by assuming the condition ${\rm W}=O(\delta^{-2})$ and investigated the growth of an initially unstable 
periodic surface perturbation and its nonlinear interaction with the higher harmonics. 
Their results imply that the surface tension plays an important role in investigating 
the stability of surface waves, which have already been pointed out by Kapitza \cite{Kapitza}. 
We remark that the condition ${\rm W}=O(\delta^{-2})$ holds for many kinds of fluid such as water and alcohol at normal temperature. 
Moreover, several authors extended the Benney's results to the three-dimensional case. 
Roskes \cite{Roskes} calculated the terms up to $O(\delta^2)$ and investigated the interactions 
between two-dimensional and three-dimensional weakly nonlinear waves on the liquid film under the condition 
${\rm W}=O(1)$, which implies that he did not consider the effect of the surface tension. 
Atherton and Homsy \cite{Atherton} and Lin and Krishna \cite{Lin} calculated the terms up to $O(\delta)$ 
and $O(\delta^2)$, respectively, under the condition  ${\rm W}=O(\delta^{-2})$, namely, 
they took the effect of the surface tension in the equation in three-dimensional case. 
Furthermore, while they considered the case where ${\rm R}=O(1)$, 
Topper and Kawahara \cite{Kawahara} derived approximate equations under the conditions 
 ${\rm W}=O(\delta^{-2})$ and ${\rm R}=O(\delta)$. 
More details or a list of useful references about the thin film approximation can be found in 
\cite{Chang, Craster, Kalli, Oron}.

\medskip
Many approximate equations are obtained from \eqref{Benney}. 
For example, by neglecting the terms of $O(\delta^2+\varepsilon^2)$, we obtain the Burgers equation 
\begin{equation}\label{Burgers}
\eta_t = -2\eta_x-4\varepsilon\eta\eta_x+\delta B(1)\eta_{xx} 
\end{equation}
with $B(1)=\frac{8}{15}\bigl(\frac54\frac1{\tan\alpha}-{\rm R}\bigr)$, from which we can recover the 
Benjamin's critical Reynolds number ${\rm R}_c=\frac54\frac1{\tan\alpha}$. 
By neglecting the terms of $O(\delta^3+\varepsilon\delta+\varepsilon^2)$,
we obtain the KdV--Burgers equation 
\begin{equation}\label{KdVB}
\eta_t = -2\eta_x-4\varepsilon\eta\eta_x+\delta B(1)\eta_{xx}+\delta^2 D(1)\eta_{xxx}, 
\end{equation}
which was named by Johnson \cite{Johnson}. 
Here, $D(1)=-2-\frac{22}{63}R^2+\frac{40}{63}\frac{\rm R}{\tan\alpha}$. 
By neglecting the terms of $O(\delta^4+\varepsilon\delta+\varepsilon^2)$ and assuming ${\rm R}>{\rm R}_c$, 
we obtain the so-called KdV--Kuramoto--Sivashinsky equation or Kawahara equation 
(see \cite{Kawahara}) 
\begin{equation}\label{KdVKS}
\eta_t = -2\eta_x-4\varepsilon\eta\eta_x+\delta B(1)\eta_{xx}+\delta^2 D(1)\eta_{xxx}+\delta^3 G(1)\eta_{xxxx}
\end{equation}
with $G(1)=-\frac23\frac{\rm W}{\sin\alpha}-\frac{157}{56}{\rm R}
   -\frac{8}{45}\frac{\rm R}{\tan^2\alpha}
   +\frac{138904}{155925}\frac{{\rm R}^2}{\tan\alpha}
   -\frac{1213952}{2027025}{\rm R}^3$. 
Note that the effect of the surface tension, namely, the Weber number ${\rm W}$, first appears in the coefficient of the fourth order derivative term in the case ${\rm W}=O(1)$. 
Moreover, by assuming $\varepsilon=1$, that is, the strongly nonlinear case and ${\rm W}=\delta^2\widetilde{\rm W}$ and neglecting the terms of $O(\delta^2)$, we obtain the so-called Benney equation (see \cite{Gjevik})
\begin{equation}\label{Benney2}
\eta_t=\bigg[-\frac23(1+\eta)^3+\delta\bigg\{\frac{2}{3\tan\alpha}(1+\eta)^3\eta_x-\frac{8\rm R}{15}(1+\eta)^6\eta_x-\frac{2\widetilde{\rm W}}{3\sin\alpha}(1+\eta)^3\eta_{xxx}\bigg\}\bigg]_x.
\end{equation}

\medskip
Now, our purpose is to give a mathematically rigorous justification of these thin film approximations 
by establishing the error estimate between the solution of the Navier--Stokes equations \eqref{ns0}--\eqref{0bc0} and 
those of approximate equations \eqref{Burgers}--\eqref{Benney2}.
More specifically, we will estimate a norm of a difference between the solution $\eta^\delta$ of Navier--Stokes  equations and the solution $\eta^{\rm app}$ of approximate equations \eqref{Burgers}--\eqref{Benney2} and show that the norm goes to 0 as $\delta\to0$.
In order to carry out the justification, the most difficult task is to derive a uniform estimate for the solution 
of the Navier--Stokes equations with respect to $\delta$ in the thin film regime $\delta\ll1$. 
In this paper, we will focus on deriving a uniform estimate of the solution with respect to $\delta$ when the Reynolds number, the angle of inclination, and the initial date are sufficiently small under the conditions $O(1)\le{\rm W}\le O(\delta^{-2})$, ${\rm R}=O(1)$, $\alpha=O(1)$, and $x\in \mathbb{T}$ or $\mathbb{R}$. 
In  \cite{Ueno}, we will give the mathematically rigorous justification of the thin film approximations. 
We remark that Bresch and Noble \cite{Bresch} have already derived a uniform estimate of the solution with respect to $\delta$ by assuming ${\rm W}=O(\delta^{-2})$, ${\rm R}=O(\delta)$, $\alpha=O(\sqrt{\delta})$, $x\in \mathbb{T}$,  and that  initial data are sufficiently small. 
Their assumptions on $\rm R$ and $\alpha$ are too restrictive when we consider the asymptotic behavior of the solution as $\delta\to0$. 
Moreover, they assumed $\varepsilon=\delta$ and excluded the case of $\varepsilon=1$, so that their uniform estimate cannot be applied to the justification for the Benney equation. Therefore, our results are not included in their works.

\medskip
Concerning a mathematical analysis of the problem, Teramoto \cite{Teramoto} showed that the initial value 
problem to the Navier--Stokes equations \eqref{ns0}--\eqref{0bc0} has a unique solution globally in time 
under the assumptions that the Reynolds number and the initial data are sufficiently small. 
Nishida, Teramoto, and Win \cite{Nishida} showed the exponential stability of the Nusselt flat film solution under the assumptions that the angle of inclination is sufficiently small 
and the flow is downward periodic in addition to the assumptions in \cite{Teramoto}. 
Furthermore, Uecker \cite{Uecker} studied the asymptotic behavior for $t\to\infty$ of the solution in the case of  $x\in\mathbb{R}$ and showed that the perturbations of the Nusselt flat film solution decay like the self-similar solution of the Burgers equation under the assumptions that the initial data are sufficiently small and ${\rm R}<{\rm R}_c$. 
However, they did not consider the $\delta$ scaling because they non-dimensionalized $x$ and $y$ components by using the same unit length $h_0$.

\medskip
The plan of this paper is as follows.
In Section 2, we rewrite the problem in a non-dimensional form and transform the problem 
in a time dependent domain to a problem in a time independent domain by using an appropriate diffeomorphism. 
Then, we give our main theorem in this paper and we remark that an outline of the proof is same as \cite{Nishida}.
In Section 3, we derive energy estimates to the transformed equations.
Only by following \cite{Nishida}, we cannot obtain a uniform estimate in $\delta$  because we cannot control lower order terms just by using energies derived in [15].
Hence,  we introduce an essentially new energy function  in order to control lower order terms which is one of difficulties to obtain a uniform boundedness of the solution in $\delta$.
Therefore, Section 3 is a key section in this paper. 
In Section 4, we give estimates for the pressure. 
In order to obtain a uniform estimate in $\delta$, we need to carefully estimate the pressure, while in \cite{Nishida} there was no need to use such a estimate.
In Section 5, we estimate carefully nonlinear terms appeared in the right-hand side of the energy inequality so that we can get a uniform estimate in $\delta$. 
Finally, combining the estimates obtained in the last three sections, 
we derive a uniform estimate of the solution in Section 6.  

\bigskip
\noindent
{\bf Notation}. \ 
We put $\Omega=\mathbb{G}\times(0,1)$ and $\Gamma=\mathbb{G}\times\{y=1\}$, 
where $\mathbb{G}$ is the flat torus $\mathbb{T}=\mathbb{R}/\mathbb{Z}$ or $\mathbb{R}$. 
For a Banach space $X$, we denote by $\|\cdot\|_X$ the norms in $X$. 
For $1\le p\le \infty$, we put $\|u\|_{L^p}=\|u\|_{L^p(\Omega)}$, $\|u\|=\|u\|_{L^2}$, 
$|u|_{L^p}=\|u(\cdot, 1)\|_{L^p(\mathbb{G})}$, and $|u|_0=|u|_{L^2}$. 
We denote by $(\cdot, \cdot)_\Omega$ and $(\cdot, \cdot)_\Gamma$ the inner products of 
$L^2(\Omega)$ and $L^2(\Gamma)$, respectively. 
For $s\ge 0$, we denote by $H^s(\Omega)$ and $H^s(\Gamma)$ the $L^2$ Sobolev spaces of order $s$ 
on $\Omega$ and $\Gamma$, respectively. 
The norms of these spaces are denoted by $\|\cdot\|_s$ and $|\cdot|_s$. 
For a function $u=u(x,y)$ on $\Omega$, 
a Fourier multiplier $P(D_x)$ ($D_x=-{\rm{i}}\partial_x$) is defined by 
$$
(P(D_x)u)(x,y)=
\begin{cases}
\displaystyle\sum_{n\in\mathbb{Z}}P(n)\hat{u}_n(y){\rm{e}}^{2\pi {\rm{i}}nx}\quad&{\rm in\  the\  case}\ \mathbb{G}=\mathbb{T}, \\
\displaystyle\int_{\mathbb{R}}P(\xi) \hat{u}(\xi, y){\rm{e}}^{2\pi {\rm{i}}\xi x}{\rm d}\xi\quad&{\rm in\  the\  case}\ \mathbb{G}=\mathbb{R}, \\
\end{cases}
$$
where $\hat{u}_n(y)=\int_0^1 u(x,y){\rm{e}}^{-2\pi {\rm i}nx}\,{\rm{d}}x$ is the Fourier coefficient and $\hat{u}(\xi, y)=\int_{\mathbb{R}} u(x,y) {\rm e}^{-2\pi{\rm i}\xi x}{\rm d}x$ is the Fourier transform in $x$. 
We put $\nabla_\delta=(\delta\partial_x, \partial_y)^{\rm T}$, $\Delta_\delta=\nabla_\delta\cdot\nabla_\delta$, and $D_\delta^kf=\{(\delta\partial_x)^i\partial_y^j f \,|\, i+j=k\}$. 
For operators $A$ and $B$, we denote by $[A, B]=AB-BA$ the commutator. 
We put 
\[
\partial_y^{-1} f(x,y)=-\int_y^1 f(x,z){\rm{d}} z.
\]
$f\lesssim g$ means that there exists a non-essential positive constant $C$ such that $f\le Cg$ holds.

\section{Reformulation of the problem and main result}
We first rewrite \eqref{ns0}--\eqref{0bc0} in a non-dimensional form. 
We will consider fluctuations on the stationary laminar flow given by \eqref{laminar flow}, 
so that we rescale the independent and dependent variables by 
$$
\begin{cases}
 x=l_0 x',\quad y=h_0y',\quad t=t_0 t',\quad \\
 \eta=a_0\eta',\quad u=U_0(\bar{u}'+\varepsilon u'),\quad v=\varepsilon V_0v',\quad p=p_0+\varepsilon P_0p',
\end{cases}
$$
where $U_0=\rho gh_0^2\sin\alpha/2\mu$, $V_0=(h_0/l_0)U_0$, $t_0=l_0/U_0$, 
$\bar{u}'=2y'-y'^2$, and $P_0=\rho gh_0\sin\alpha$. 
Putting these into \eqref{ns0}--\eqref{0bc0} and dropping the prime sign in the notation, we obtain 
\begin{equation}\label{ns}
\left\{
 \begin{array}{lcl}
  \delta\bm{u}^\delta_t+\bigl((\bm{U}+\varepsilon\bm{u}^\delta)\cdot\nabla_\delta\bigr)\bm{u}^\delta
   +(\bm{u}^\delta\cdot\nabla_\delta)\bm{U}+\dfrac{2}{{\rm R}}\nabla_\delta p
   -\dfrac{1}{{\rm R}}\Delta_\delta\bm{u}^\delta=\bm{0}
  & \mbox{in} & \Omega_{\varepsilon}(t), \; t>0, \\ 
  \nabla_\delta\cdot\bm{u}^\delta = 0 
  & \mbox{in} & \Omega_{\varepsilon}(t), \; t>0,
 \end{array}
\right.
\end{equation}
\begin{equation}\label{bc1}
\left\{
 \begin{array}{lcl}
  \bigl(\bm{D}_\delta(\varepsilon\bm{u}^\delta+\bm{U})-\varepsilon p\bm{I}\bigr)\bm{n}^\delta & & \\
  \quad
   = \bigg(-\dfrac{1}{\tan\alpha}\varepsilon\eta
    +\dfrac{\delta^2{\rm W}}{\sin\alpha}\dfrac{\varepsilon\eta_{xx}}{(1+(\varepsilon\delta\eta_x)^2)^{\frac32}}\bigg)
    \bm{n}^\delta
    & \mbox{on} & \Gamma_{\varepsilon}(t), \; t>0, \\[3mm]
  \eta_t+\bigl(1-(\varepsilon\eta)^2+\varepsilon u\bigr)\eta_x-v = 0 
    & \mbox{on} & \Gamma_{\varepsilon}(t), \; t>0, 
 \end{array}
\right.
\end{equation}
\begin{equation}\label{bc0}
\bm{u}^\delta=\bm{0} \quad\mbox{on}\quad \Sigma, \; t>0, 
\end{equation}
where $\bm{u}^\delta=(u,\delta v)^{\rm T}$, $\bm{U}=(\bar{u},0)^{\rm T}$, $\bar{u}=2y-y^2$, 
$\bm{D}_\delta\bm{f}=\frac{1}{2}\bigl\{\nabla_\delta(\bm{f}^{\rm T})+\bigl(\nabla_\delta(\bm{f}^{\rm T})\bigr)^{\rm T}\bigr\}$, 
$\bm{n}^\delta=(-\varepsilon\delta\eta_x,1)^{\rm T}$, ${\rm R}=\rho U_0 h_0/\mu$ is the Reynolds number, 
and ${\rm W}=\sigma/\rho g h_0^2$ is the Weber number. 
In this scaling, the liquid domain $\Omega_{\varepsilon}(t)$ and the liquid surface $\Gamma_{\varepsilon}(t)$ are of the 
forms 
$$
\left\{
\begin{array}{l} 
\Omega_{\varepsilon}(t)=\{(x,y)\in\mathbb{R}^2\ |\ 0<y<1+\varepsilon\eta(x,t)\}, \\
\Gamma_{\varepsilon}(t)=\{(x,y)\in\mathbb{R}^2\ |\ y=1+\varepsilon\eta(x,t)\}.
\end{array}
\right.
$$

Next, we transform the problem in the moving domain $\Omega_{\varepsilon}(t)$ 
to a problem in the fixed domain $\Omega$ by using an appropriate diffeomorphism 
$\Phi:\Omega \rightarrow \Omega_{\varepsilon}(t)$ defined by 
\begin{equation}\label{diffeomorphism}
\Phi(x,y,t)=\bigl(x, y(1+\varepsilon\tilde{\eta}(x, y, t))\bigr),
\end{equation}
where $\tilde{\eta}$ is an extension of $\eta$ to $\Omega$. 
We need to choose the extension $\tilde{\eta}$ carefully and in this paper 
we adopt the following extension. 
For $\phi\in H^s(\Gamma)$, we define its extension $\tilde{\phi}$ to $\Omega$ by 
\begin{equation}
\tilde{\phi}(x,y)
 =
\begin{cases}
 \displaystyle\sum_{n\in \mathbb{Z}}\frac{\hat{\phi}_n}{1+(\delta n(1-y)y)^4}{\rm e}^{2\pi {\rm i}nx}\quad&{\rm in\ the\ case}\ \mathbb{G}=\mathbb{T},\\
\displaystyle\int_{\mathbb{R}}\frac{\hat{\phi}(\xi)}{1+(\delta \xi(1-y)y)^4}{\rm e}^{2\pi {\rm i}\xi x}{\rm d}\xi\quad&{\rm in\ the\ case}\ \mathbb{G}=\mathbb{R}.
\end{cases}
\label{def-extension}
\end{equation}
By the definition, it is easy to see that 
\begin{equation}\label{extension0}
\partial_y^j \tilde{\phi}(x,1)=\partial_y^j \tilde{\phi}(x,0)=0\qquad {\rm for}\ 1\le j\le 3.
\end{equation}
As usual, this extension operator has a regularizing effect so that $\tilde{\phi}\in H^{s+\frac12}(\Omega)$. 
However, if we use such a regularizing property, then we need to pay the cost of a power of $\delta$. 
Moreover, in this extension, $\partial_y$ corresponds to $\delta\partial_x$. 
More precisely, we have the following lemma, whose proof is quite standard, so we omit it.

\begin{lem}\label{extension}
Let $i$ and $j$ be non-negative integers such that $j\le4$. 
Then, for the extension \eqref{def-extension} we have 
\begin{align}
& \|\partial_x^i\partial_y^j \tilde{\phi}\| \lesssim \delta^j |\partial_x^{i+j}\phi|_0, \label{extension2} \\
& \|\partial_x^i\partial_y^j \tilde{\phi}\|_{L^\infty} \lesssim \delta^j |\partial_x^{i+j}\phi|_1. 
 \label{extension4}
\end{align}
If, in addition,  $i+j\ge1$, then 
\begin{align}
 \|\partial_x^i\partial_y^j \tilde{\phi}\| \lesssim \delta^{j-\frac12} ||D_x|^{i+j-\frac{1}{2}}\phi|_0. 
   \label{extension1} 
\end{align}
\end{lem}

The solenoidal condition on the velocity field is destroyed in general by the transformation. 
To keep the condition, following Beale \cite{Beale}, 
we also change the dependent variables and introduce new unknown functions $(u', v', p')$ 
defined in $\Omega$ by 
$$
u'={J}(u\circ\Phi), \quad
v'=v\circ\Phi-y\varepsilon\tilde{\eta}_x(u\circ\Phi), \quad
p'=p\circ\Phi,
$$
where ${J}=1+\varepsilon(y\tilde{\eta})_y$ is the Jacobian of the diffeomorphism $\Phi$. 
Putting $a_1=-yJ^{-1}\varepsilon\delta\tilde{\eta}_x$, $b_1={J^{-1}}-1$, and 
$$
A_1
=
\begin{pmatrix}
 1+b_1 & 0 \\
 -a_1 & 1
\end{pmatrix}
=N_1+I, \quad 
\bm{u}'^{\delta}=\binom{u'}{\delta v'},
$$
we have
\begin{equation}\label{u-trans}
 \bm{u}^\delta\circ\Phi=A_1\bm{u}'^{\delta}.
\end{equation}
Here $N_1$ is the nonlinear part of $A_1$. 
We note that $b_1$ is the term which is hard to handle 
because it contains the term without $\delta$ in the coefficient. 
Then, the second equation in \eqref{bc1} is transformed to 
\begin{equation}\label{kinematic}
\eta_t+\eta_x-v'=h_3,
\end{equation}
where $h_3=\varepsilon^2\eta^2\eta_x$. 
We easily obtain that
\begin{numcases}
{}
(\nabla_\delta\phi)\circ\Phi
 = A_2\nabla_\delta(\phi\circ\Phi), \label{nd} \\
(\Delta_\delta \phi)\circ\Phi
 = \delta^2(\phi\circ\Phi)_{xx}+(1+b_2)(\phi\circ\Phi)_{yy}+P_\delta(\tilde{\eta},D)(\phi\circ\Phi), \label{ld} \\
\delta(\phi_t\circ\Phi)
 = \delta(\phi\circ\Phi)_t- y {J^{-1}}\varepsilon\delta\tilde{\eta}_t(\phi\circ\Phi)_y, \label{dt}
\end{numcases}
where
$$
A_2=
\begin{pmatrix}
 1 & a_1 \\
 0 & 1+b_1
\end{pmatrix}
= N_2+I, 
$$
$N_2$ is the nonlinear part of $A_2$, $b_2=a_1^2+2b_1+b_1^2$, 
and $P_\delta(\tilde{\eta},D)$ is a second order differential operator 
defined by 
$P_\delta(\tilde{\eta}, D)f
 = 2\delta a_1f_{xy}+\bigl\{\delta a_{1x}+a_1 a_{1y}+(1+b_1)b_{1y}\bigr\}f_y$.

We proceed to transform the equations. 
We begin to transform the equations in \eqref{ns}. 
By \eqref{u-trans} and \eqref{dt}, we obtain 
\begin{equation}\label{stationary}
\delta\bm{u}_t^\delta\circ\Phi=\delta A_1\bm{u}'^{\delta}_t+\bm{f}_1,
\end{equation}
where 
$$
\bm{f}_1
= 
\delta A_{1t}\bm{u}'^{\delta}- y {J^{-1}}\varepsilon\delta\tilde{\eta}_t(A_1\bm{u}'^{\delta})_y.
$$
By \eqref{u-trans} and \eqref{nd}, we obtain 
\begin{equation}\label{advection}
\bigl\{\bigl((\bm{U}+\varepsilon\bm{u}^\delta)\cdot\nabla_\delta\bigr)\bm{u}^\delta
 +(\bm{u}^\delta\cdot\nabla_\delta)\bm{U}\bigr\}\circ\Phi
=(\bm{U}\cdot\nabla_\delta)\bm{u}'^{\delta}+(\bm{u}'^{\delta}\cdot\nabla_\delta)\bm{U}+\bm{f}_2,
\end{equation}
where
\begin{align*}
\bm{f}_2
&= (\bm{U}\cdot\nabla_\delta)N_1\bm{u}'^{\delta}+(\bm{U}\cdot N_2\nabla_\delta)A_1\bm{u}'^{\delta}
  +\bigl((\bm{V}+\varepsilon A_1\bm{u}'^{\delta})\cdot A_2\nabla_\delta\bigr)A_1\bm{u}'^{\delta} \\
&\quad
  +(\bm{u}'^{\delta}\cdot N_2\nabla_\delta)\bm{U}+\bigl((N_1\bm{u}'^{\delta})\cdot (A_2\nabla_\delta)\bigr)\bm{U}
  +\bigl((A_1\bm{u}'^{\delta})\cdot (A_2\nabla_\delta)\bigr)\bm{V}, \\
\bm{V}
&= \binom{2\varepsilon y\tilde{\eta}-2\varepsilon y^2\tilde{\eta}-(\varepsilon y\tilde{\eta})^2}{0}. 
\end{align*}
By \eqref{nd}, we have 
\begin{equation}\label{press}
(\nabla_\delta p)\circ\Phi=A_2\nabla_\delta p'.
\end{equation}
By \eqref{u-trans} and \eqref{ld}, we obtain 
\begin{equation}\label{viscosity}
(\Delta_\delta\bm{u}^\delta)\circ\Phi
= A_1\bigl(\delta^2\bm{u}'^{\delta}_{xx}+(I+A_3)\bm{u}'^{\delta}_{yy}\bigr)+\bm{f}_3,
\end{equation}
where
\begin{align*}
A_3
&= 
\begin{pmatrix}
 b_2  &0 \\
 0 & 0
\end{pmatrix}, \\
\bm{f}_3
&= [\delta^2\partial_x^2,A_1]\bm{u}'^{\delta}+(1+b_2)\bigl[\partial_y^2,A_1\bigr]\bm{u}'^{\delta}
  +P_\delta(\tilde{\eta},D)\bigl(A_1\bm{u}'^{\delta}\bigr)+A_1\binom{0}{\delta b_2 v'_{yy}}.
\end{align*}
Thus combining \eqref{stationary}--\eqref{viscosity}, we transform the first equation in \eqref{ns} to 
\begin{equation}\label{Ns}
\delta \bm{u}'^{\delta}_t+(\bm{U}\cdot\nabla_\delta)\bm{u}'^{\delta}
 +(\bm{u}'^{\delta}\cdot\nabla_\delta)\bm{U}+\dfrac{2}{{\rm R}}(I+A_4)\nabla_\delta p'
 -\frac{1}{{\rm R}}\left(\delta^2\bm{u}'^{\delta}_{xx}+(I+A_3)\bm{u}'^{\delta}_{yy}\right)
=\bm{f},
\end{equation}
where 
\begin{align}
\bm{f}
&= -N\bigl\{(\bm{U}\cdot\nabla_\delta)\bm{u}'^{\delta}+(\bm{u}'^{\delta}\cdot\nabla_\delta)\bm{U}
 \bigr\}+A_1^{-1}\biggl(-\bm{f}_1-\bm{f}_2+\frac{1}{{\rm R}}\bm{f}_3\biggr),\label{NLdomain}\\
A_4
&= A_1^{-1}A_2-I=
\begin{pmatrix}\label{NLp}
 (y\varepsilon\tilde{\eta})_y & -y\varepsilon\delta\tilde{\eta}_x\\
 -y\varepsilon\delta\tilde{\eta}_x 
  & J^{-1}\bigl((y\varepsilon\delta\tilde{\eta}_x)^2-(y\varepsilon\tilde{\eta})_y\bigr)
\end{pmatrix},
\end{align}
and $N$ is the nonlinear part of $A_1^{-1}$. 
We remark that $\bm{f}$ is a collection of nonlinear terms, which does not contain 
$\bm{u}'^{\delta}_t$, $u'_{yy}$, $\nabla_\delta p'$, nor any function of $\tilde{\eta}$ only.

Next, we transform the boundary conditions. 
By \eqref{u-trans} and \eqref{nd}, we see that 
$$
\bigr\{\bigl(\bm{D}_\delta(\varepsilon\bm{u}^\delta+\bm{U})-\varepsilon pI\bigr)\bm{n}^\delta\bigr\}\circ\Phi
= \binom{\frac{1}{2}\varepsilon(\delta^2v'_x+u'_y-2\eta)}{\varepsilon\delta v'_y}-\varepsilon p'\bm{n}^\delta+\bm{h}
 \quad\mbox{on}\quad \Gamma, 
$$
where 
\begin{align*}
\bm{h}
&= -\binom{\varepsilon^2\delta^2\eta_x u'_x}{\frac{1}{2}\varepsilon^2\delta\eta_x(\delta^2v'_x+u'_y-2\eta)} \\
&\quad
  +\frac{\varepsilon}{2}\bigl\{\nabla_\delta(N_1\bm{u}'^{\delta})^{\rm T}
  +\bigl(\nabla_\delta(N_1\bm{u}'^{\delta})^{\rm T}\bigr)^{\rm T}
  +N_2\nabla_\delta(A_1\bm{u}'^{\delta})^{\rm T}
  +\bigl(N_2\nabla_\delta(A_1\bm{u}'^{\delta})^{\rm T}\bigr)^{\rm T}\bigr\}\bm{n}^\delta. 
\end{align*}
Taking the inner product of a tangential vector $\bm{t}^\delta=(1, \varepsilon\delta\eta_x)^{\rm T}$ 
with the first equation in \eqref{bc1}, we obtain 
\begin{equation}\label{0tan-bc1}
\delta^2v'_x+u'_y-2\eta=h_4 \quad\mbox{on}\quad \Gamma,
\end{equation}
where
$$
h_4=-\frac{2}{\varepsilon}(\varepsilon^2\delta^2\eta_x v'_y+\bm{h}\cdot\bm{t}^\delta).
$$
On the other hand, taking the inner product of a normal vector $\bm{n}^\delta$ 
with the first equation in \eqref{bc1}, we obtain 
\begin{equation}\label{nor-bc1}
p'-\delta v'_y-\dfrac{1}{\tan\alpha}\eta+\dfrac{\delta^2{\rm W}}{\sin\alpha}\eta_{xx}=h_2
 \quad\mbox{on}\quad \Gamma,
\end{equation}
where 
\begin{align}\label{NLnor}
h_2 &= \frac{1}{\varepsilon}
 \bigg\{ -\frac{(\varepsilon\delta\eta_x)^2}{1+(\varepsilon\delta\eta_x)^2}\varepsilon\delta v'_y
  +\frac{1}{1+(\varepsilon\delta\eta_x)^2}
  \biggl(-\frac{1}{2}\varepsilon^2\delta\eta_x(\delta^2v'_x+u'_y-2\eta)+\bm{h}\cdot\bm{n}^\delta\biggr) \bigg\}\\
 &\phantom{=\ } +\dfrac{\delta^2{\rm W}}{\sin\alpha}\biggl(1-\frac{1}{(1+(\varepsilon\delta\eta_x)^2)^\frac{3}{2}}\biggr)
  \eta_{xx}\notag\\
 & =:h_{2,1}+\delta^2{\rm W}h_{2,2}\notag
\end{align}
and $h_2$ does not contain $p'$ nor any function of $\eta$ only. Note that the term $\delta^2{\rm W}h_{2,2}$ is the only nonlinear term which contains ${\rm W}$. 
Here, by a straightforward calculation we see that 
$$
\bm{h}\cdot\bm{t}^\delta=\varepsilon(b_4u'_y+h_5),
$$
where 
\begin{align*}
b_4
&= -\frac{1}{2}(\varepsilon\delta\eta_x)^2
 +\biggl\{
  \begin{pmatrix}
   a_1(1+b_1) & \frac{1}{2}(-a_1^2+b_1(2+b_1)) \\
   \frac{1}{2}(-a_1^2+b_1(2+b_1)) & -a_1(1+b_1)
  \end{pmatrix}
  \bm{n}^\delta\biggr\}\cdot\bm{t}^\delta, \\
h_5
&= -\varepsilon\delta^2\eta_xu_x'-\frac12(\varepsilon\delta\eta_x)^2(\delta^2v_x'-2\eta) \\
&\quad
 +\biggl\{
 \begin{pmatrix}
  \delta(b_1u')_x & \frac{1}{2}\{\delta(-a_1u')_x-a_1a_{1y} u'+\delta a_1 v_y'\} \\
  \frac{1}{2}\{\delta(-a_1u')_x-a_1a_{1y} u'+\delta a_1 v_y'\} & -a_{1y} (1+b_1)u'+\delta b_1 v_y'
 \end{pmatrix}
 \bm{n}^{\delta}\biggr\}\cdot\bm{t}^{\delta}, 
\end{align*}
and $h_5$ does not contain $u'_y$. 
Thus we can rewrite \eqref{0tan-bc1} as 
\begin{equation}\label{tan-bc1}
\delta^2 v'_x+u'_y-(2+b_3)\eta=h_1 \quad\mbox{on}\quad \Gamma,
\end{equation}
where 
\begin{align}
b_3 &= -\frac{4b_4}{1+2b_4}, \label{def-b3} \\
h_1 &= \frac{2b_4}{1+2b_4}\delta^2v'_x-\frac{2}{1+2b_4}(\varepsilon\delta^2\eta_x v'_y+h_5). \label{NLtan}
\end{align}
Note that $h_1$ does not contain $u'_y$, $p'$, nor any function of $\eta$ only.

Summarizing \eqref{kinematic}, \eqref{Ns}, \eqref{nor-bc1}, and \eqref{tan-bc1} 
and dropping the prime sign in the notation, we have 
\begin{equation}\label{NS}
\left\{
 \begin{array}{lcl}
  \delta \bm{u}^\delta_t+(\bm{U}\cdot\nabla_\delta)\bm{u}^\delta+(\bm{u}^\delta\cdot\nabla_\delta)\bm{U} \\
  \qquad +\dfrac{2}{{\rm R}}(I+A_4)\nabla_\delta p
   -\dfrac{1}{{\rm R}}\big\{\delta^2\bm{u}^\delta_{xx}+(I+A_3)\bm{u}^\delta_{yy}\big\}
  =\bm{f} & \mbox{in} & \Omega, \; t>0, \\[1mm]
  u_x+v_y=0 & \mbox{in} & \Omega, \; t>0,
 \end{array}
\right.
\end{equation}
\begin{equation}\label{BC1}
\left\{
 \begin{array}{lcl}
  \delta^2v_x+u_y-(2+b_3)\eta=h_1 & \mbox{on} & \Gamma, \; t>0, \\[1mm]
  p-\delta v_y-\dfrac{1}{\tan\alpha}\eta+\dfrac{\delta^2{\rm W}}{\sin\alpha}\eta_{xx}=h_2
   & \mbox{on} & \Gamma, \; t>0, \\[1mm]
  \eta_t+\eta_x-v=h_3 & \mbox{on} & \Gamma, \; t>0,
 \end{array}
\right.
\end{equation}
\begin{equation}\label{BC0}
u=v=0 \quad\mbox{on}\quad \Sigma, \; t>0.
\end{equation}
In the following, we will consider the initial value problem to \eqref{NS}--\eqref{BC0} 
under the initial conditions 
\begin{equation}\label{initial}
\eta|_{t=0}=\eta_0\quad {\rm on}\ \ \Gamma, \quad (u, v)^{\rm T}|_{t=0}=(u_0, v_0)^{\rm T}\quad{\rm in}\ \ \Omega. 
\end{equation}
Here we denote $b_3$ and $h_1$ determined from the initial data by $b_3^{(0)}$ and $h_1^{(0)}$, respectively.

Now, we are ready to state our main result in this paper.

\begin{thm}$($Uniform estimate$)$ \label{main thm}
There exist small positive constants ${\rm R}_0$ and $\alpha_0$ such that the following statement holds: 
Let $m$ be an integer satisfying $m\ge2$, $0<{\rm R}_1\le {\rm R}_0$, $0<{\rm W}_1\le{\rm W}_2$, and $0<\alpha\le \alpha_0$. There exist positive constants $c_0$ and $T$ such that if the initial data $(\eta_0, u_0, v_0)$ 
and the parameters $\delta$, $\varepsilon$, ${\rm R}$, and ${\rm W}$ satisfy the compatibility conditions 
$$
\left\{
\begin{array}{lcl}
 u_{0x}+v_{0y}=0 & \mbox{in} & \Omega, \\
 u_{0y}+\delta^2 v_{0x}-(2+b_3^{(0)})\eta_0=h_1^{(0)} & \mbox{on} & \Gamma, \\
 u_0=v_0=0 & \mbox{on} & \Sigma,
\end{array}
\right.
$$
and
$$
\begin{cases}
|(1+\delta|D_x|)^2\eta_0|_2+\|(1+|D_x|)^2(u_0, \delta v_0)^{\rm T}\| +\|(1+|D_x|)^2D_{\delta}(u_0, \delta v_0)^{\rm T}\|\\
\quad+\|(1+|D_x|)^2D_{\delta}^2(u_0, \delta v_0)^{\rm T}\|+\delta^2{\rm W}|(1+\delta|D_x|)\eta_{0x}|_3+\sqrt{\delta^2{\rm W}}\|(1+|D_x|)^2\delta v_{0xy}\|\le c_0, \\
|(1+\delta|D_x|)^2\eta_0|_m+\|(1+|D_x|)^m(u_0, \delta v_0)^{\rm T}\| +\|(1+|D_x|)^mD_{\delta}(u_0, \delta v_0)^{\rm T}\|\\
\quad+\|(1+|D_x|)^mD_{\delta}^2(u_0, \delta v_0)^{\rm T}\|+\delta^2{\rm W}|(1+\delta|D_x|)\eta_{0x}|_{m+1}+\sqrt{\delta^2\rm{W}}\|(1+|D_x|)^m\delta v_{0xy}\|\le M, \\
0<\delta, \varepsilon\le1, \quad {\rm R_1}\le{\rm R}\le{\rm R_0}, \quad {\rm W}_1\le{\rm W}\le\delta^{-2}{\rm W}_2,
\end{cases}
$$
then the initial value problem \eqref{NS}--\eqref{initial} has a unique solution $(\eta, u, v ,p)$ 
on the time interval $[0, T/\varepsilon]$ and the solution satisfies the estimate 
\begin{align*}
& |(1+\delta|D_x|)^2\eta(t)|_m^2+\delta^2|\eta_t(t)|_m^2+\delta^2{\rm W}\bigl\{|(1+\delta|D_x|)^2\eta_{x}(t)|_m^2+\delta^2|\eta_{tx}(t)|_m^2\bigr\} \\
& \quad
  +\|(1+|D_x|)^m(1+\delta|D_x|)^2\bm{u}^\delta(t)\|^2+\|(1+|D_x|)^m\bm{u}^\delta_y(t)\|^2
  +\delta^2\|(1+|D_x|)^m\bm{u}^\delta_t(t)\|^2 \\
& +\int_0^t\big\{\delta|\eta_x(\tau)|_m^2+\delta|(1+\delta|D_x|)^{\frac52}\eta_t(\tau)|_m^2 \\
& \makebox[3em]{} +(\delta^2{\rm W})\delta|\eta_{xx}(\tau)|_m^2+(\delta^2{\rm W})^2\bigl\{\delta|\eta_{xxx}(\tau)|_m^2+\delta^2||D_x|^{\frac72}\eta(\tau)|_m^2\bigr\} \\
& \makebox[3em]{}
  +\delta\|(1+|D_x|)^m\bm{u}_x^{\delta}(\tau)\|^2
  +\delta\|(1+|D_x|)^m(1+\delta|D_x|)\nabla_\delta\bm{u}_x^\delta(\tau)\|^2+\delta\|(1+|D_x|)^m\nabla_\delta\bm{u}^\delta_t(\tau)\|^2 \\
&\makebox[3em]{}
  +\delta^{-1}\|(1+|D_x|)^m(1+\delta|D_x|)\bm{u}^\delta_{yy}(\tau)\|^2
  +\delta\|(1+|D_x|)^m\partial_y^{-1}p_x(\tau)\|^2 \\
&\makebox[3em]{}
  +\delta^{-1}\|(1+|D_x|)^m(1+\delta|D_x|)\nabla_\delta p(\tau)\|^2
  +\delta\|(1+|D_x|)^{m-1}\nabla_\delta p_t(\tau)\|^2
  \big\}{\rm d}\tau\le C
\end{align*}
for $0\le t\le T/\varepsilon$ with a constant $C=C({\rm R}_1, {\rm W}_1, {\rm W}_2, \alpha, M)$ independent of 
$\delta$, $\varepsilon$, ${\rm R}$, and ${\rm W}$. 
Moreover, the following uniform estimate holds. 
\begin{align}\label{uniform estimate}
|\eta(t)|_m
&+\|(1+|D_x|)^{m-1}u(t)\|_1+\|\partial_x^m u_y(t)\| \\
&+\|(1+|D_x|)^{m-2}v(t)\|_1+\|\partial_x^{m-1}v_{yy}(t)\| \le C \notag
\end{align}
for $0\le t\le T/\varepsilon$. 
If, in addition, $0\le\varepsilon\lesssim\delta$, 
then the solution can be extended for all $t\ge0$ and the above estimates hold for $t\ge0$. 
\end{thm}

\noindent
{\bf Remark 2.1.}
In the case $\varepsilon\simeq1$, this theorem gives a uniform boundedness of the solution only for a 
short time interval $[0,T]$. 
However, this is essential and we cannot extend this uniform estimate for all $t\geq0$ in general, 
because by \eqref{Benney} we see that the limiting equation for $\eta$ as $\delta\to0$ becomes 
a nonlinear hyperbolic conservation law of the form 
$$
\eta_t+2(1+\varepsilon\eta)^2\eta_x=0,
$$
whose solution will have a singularity in finite time in general. 

\medskip
\noindent
{\bf Remark 2.2.}
In the case where $\mathbb{G}=\mathbb{T}, \varepsilon\lesssim\delta$, and $\int_0^1\eta_0(x){\rm d}x=0$, we also obtain the following exponential decay in time  property of the solution.
\begin{align}\label{decay}
& |(1+\delta|D_x|)^2\eta(t)|_m^2+\delta^2|\eta_t(t)|_m^2+\delta^2{\rm W}\bigl\{|(1+\delta|D_x|)^2\eta_{x}(t)|_m^2+\delta^2|\eta_{tx}(t)|_m^2\bigr\} \\
& 
  +\|(1+|D_x|)^m(1+\delta|D_x|)^2\bm{u}^\delta(t)\|^2+\|(1+|D_x|)^m\bm{u}^\delta_y(t)\|^2
  +\delta^2\|(1+|D_x|)^m\bm{u}^\delta_t(t)\|^2 \le C{\rm e}^{-c\delta t}.\notag
\end{align}

\medskip
\noindent
{\bf Remark 2.3.}
In order to derive a uniform estimate in $\rm R$, the constant $C$ in the above estimate   is required to depend on a lower bound  $\rm R_1$ of $\rm R$ for a technical reason. However, for a justification of the thin film approximation this restriction matters little because we are interested in the case where $\rm R$ is close enough to $\rm R_c$.
\section{Energy estimates}
We recall two fundamental inequalities which have a key role in this paper.

\begin{lem}$($Korn's inequality$)$\label{Korn}
There exists a constant $K$ independent of $\delta$ such that for any $0<\delta\leq1$ and 
$\bm{u}=(u,v)^{\rm T}$ satisfying 
$$
\left\{
 \begin{array}{lll}
  u_x+v_y=0 & \hbox{in} & \Omega, \\
  u=v=0 & \hbox{on} & \Sigma,
 \end{array} 
\right.
$$
we have
$$
\int\!\!\!\int_{\Omega}(\delta^2u_x^2+u_y^2+\delta^4v_x^2+\delta ^2v_y^2){\rm d}x{\rm d}y 
\leq K \int\!\!\!\int_{\Omega}\bigl(2\delta^2u_x^2+(u_y+\delta ^2v_x)^2+2\delta^2v_y^2\bigr){\rm d}x{\rm d}y. 
$$
\end{lem}

\noindent
{\bf Remark 3.1.}\ 
Teramoto and Tomoeda \cite{Tomoeda} proved that the best constant of $K$ is $3$. 
Note that in the case of $\delta=1$, this inequality is well-known.

\begin{lem}$($Trace theorem$)$\label{trace}
For $0<\delta\le1$, we have
\[
| f|_0 ^2 + \delta | |D_x|^{\frac12}f|_0 ^2 \lesssim \| f\| ^2 + \delta ^2 \| f_x\| ^2 + \| f_y\| ^2.
\]
\end{lem}

\medskip
\noindent
{\bf Remark 3.2.}\ 
This trace theorem is also well-known in the case of $\delta=1$.

\medskip
\noindent
We omit the proofs of the above lemmas because we only have to modify slightly the proofs in the case of $\delta=1$. 

The following proposition is a slight modification of the energy estimate obtained in \cite{Nishida}.

\begin{prop}\label{highEest}
There exists a positive constant ${\rm R}_0$ such that if $0<{\rm R}\le{\rm R}_0$, 
then the solution $(\eta, u, v, p)$ of \eqref{NS}--\eqref{BC0} satisfies 
\begin{align}\label{lowest}
& \dfrac{\delta}{2}\dfrac{{\rm d}}{{\rm d}t}\bigg\{\|\bm{u}^\delta\|^2
  +\dfrac{2}{{\rm R}}\bigg(\dfrac{1}{\tan\alpha}|\eta|_0^2
  +\dfrac{\delta^2{\rm W}}{\sin\alpha}|\eta_{x}|_0^2\bigg)\bigg\}
 +\frac{1}{4KR}\|\nabla_\delta\bm{u}^\delta\|^2 \\
&\le \frac{4K}{\rm R}(|\eta|_0^2+|b_3\eta|_0^2)+\dfrac{1}{{\rm R}}(h_1,u)_\Gamma
  -\dfrac{2}{{\rm R}}( h_2, \delta v)_\Gamma\notag \\
&\quad 
 +\dfrac{2}{{\rm R}}(\dfrac{1}{\tan\alpha}\eta-\dfrac{\delta^2{\rm W}}{\sin\alpha}\eta_{xx},\delta h_3)_\Gamma
 +(\bm{F}_1, \bm{u}^\delta)_\Omega, \notag
\end{align}
where $K$ is the constant in Korn's inequality and 
\begin{equation}\label{NLdomainF}
\bm{F}_1=\bm{f}-\dfrac{2}{{\rm R}} A_4\nabla_\delta p+\dfrac{1}{{\rm R}} \binom{b_2 u_{yy}}{0}.
\end{equation}
\end{prop}

\noindent
{\it Proof}. \ 
Note that Lemma \ref{Korn} implies 
\begin{equation}\label{Korn2}
\|\nabla_\delta\bm{u}^\delta\|^2\le K
|\mkern-1.5mu|\mkern-1.5mu| \bm{u}^\delta |\mkern-1.5mu|\mkern-1.5mu|^2,
\end{equation}
where $|\mkern-1.5mu|\mkern-1.5mu| \bm{u}^\delta |\mkern-1.5mu|\mkern-1.5mu|^2 
= 2\|\delta u_x\|^2+\|u_y+\delta^2 v_x\|^2+2\|\delta v_y\|^2$. 
Taking the inner product of $\bm{u}^\delta$ with the first equation in \eqref{NS}, we have 
\begin{equation}\label{enecal}
\frac{\delta}{2}\frac{{\rm d}}{{\rm d}t}\|\bm{u}^\delta\|^2
 +(u, \bar{u}_y\delta v)_\Omega+\dfrac{1}{{\rm R}}(2\nabla_\delta p
 -\Delta_\delta\bm{u}^\delta, \bm{u}^\delta)_\Omega = (\bm{F}_1, \bm{u}^\delta)_\Omega.
\end{equation}
Using the second equation in \eqref{NS} and integration by parts in $x$ and $y$, we see that 
\begin{align*}
& (2\nabla_\delta  p-\Delta_\delta\bm{u}^\delta, \bm{u}^\delta)_\Omega \\
&= 2( p, \delta v)_\Gamma-(2\delta^2 u_{xx}+\delta^2  v_{xy}+ u_{yy}, u)_\Omega
  -(\delta^3 v_{xx}+2\delta v_{yy}+\delta u_{xy}, \delta v)_\Omega \\
&= 2( p, \delta  v)_\Gamma+2\|\delta u_x\|^2+(\delta^2 v_x+ u_y,  u_y)_\Omega-(\delta^2 v_x+ u_y,  u)_\Gamma \\
&\quad
  + 2\|\delta v_y\|^2-2(\delta v_y, \delta v)_\Gamma+(\delta^2 v_x+ u_y, \delta^2 v_x)_\Omega \\
&= |\mkern-1.5mu|\mkern-1.5mu| \bm{u}^\delta|\mkern-1.5mu|\mkern-1.5mu|^2
 +2( p-\delta v_y, \delta v)_\Gamma-(\delta^2 v_x+ u_y, u)_\Gamma.
\end{align*}
By \eqref{BC1} and integration by parts in $x$, the boundary terms in the right-hand side of the above equality 
are calculated as 
\begin{align}\label{nishida}
2( p-\delta v_y, \delta v)_\Gamma
&= 2(\dfrac{1}{\tan\alpha}\eta-\dfrac{\delta^2{\rm W}}{\sin\alpha}\eta_{xx}, \delta(\eta_t+\eta_x-h_3))_\Gamma
  +2(h_2, \delta v)_\Gamma \\
&= \delta\frac{{\rm d}}{{\rm d}t}\bigg\{\dfrac{1}{\tan\alpha}|\eta|_0^2
  +\dfrac{\delta^2{\rm W}}{\sin\alpha}|\eta_x|_0^2\bigg\}+2(h_2, \delta v)_\Gamma \notag \\
&\quad
  -2(\dfrac{1}{\tan\alpha}\eta-\dfrac{\delta^2{\rm W}}{\sin\alpha}\eta_{xx}, \delta h_3)_\Gamma \notag
\end{align}
and $-(\delta^2 v_x+u_y, u)_\Gamma=-((2+b_3)\eta, u)_\Gamma-(h_1, u)_\Gamma$. 
Moreover, by the Cauchy--Schwarz and Poincar$\acute{\rm e}$'s inequalities we see that 
$|(u, \bar{u}_y\delta v)_\Omega| \leq 2\|u\|\|\delta v\| \leq \|\bm{u}^{\delta}\|^2 \leq \|\bm{u}^{\delta}_y\|^2
\leq \|\nabla_{\delta}\bm{u}^{\delta}\|^2$ and that 
$\frac{2}{\rm R}|(\eta,u)_{\Gamma}| \leq \frac{2}{\rm R}|\eta|_0\|u_y\| 
\leq \frac{1}{4KR}\|u_y\|^2+\frac{4K}{R}|\eta|_0^2$. 
Here, we used the inequality $|u(\cdot,1)|_0=|u(\cdot,1)-u(\cdot,0)|_0\leq\|u_y\|$ 
thanks to the boundary condition \eqref{BC0}. 
In the following, we use frequently this type of inequality without any comment. 
Thus we can rewrite \eqref{enecal} as 
\begin{align*}
& \dfrac{\delta}{2}\dfrac{{\rm d}}{{\rm d}t}\bigg\{\|\bm{u}^\delta\|^2
  +\dfrac{2}{{\rm R}}\bigg(\dfrac{1}{\tan\alpha}|\eta|_0^2
  +\dfrac{\delta^2{\rm W}}{\sin\alpha}|\eta_{x}|_0^2\bigg)\bigg\}
 +\dfrac{1}{2K{\rm R}}\|\nabla_{\delta}\bm{u}^\delta\|^2 \\
&\le \|\nabla_{\delta}\bm{u}^\delta\|^2+\frac{4K}{\rm R}(|\eta|_0^2+|b_3\eta|_0^2)
 +\dfrac{1}{{\rm R}}(h_1,u)_\Gamma-\dfrac{2}{{\rm R}}( h_2, \delta v)_\Gamma \\
&\quad
 +\dfrac{2}{{\rm R}}(\dfrac{1}{\tan\alpha}\eta-\dfrac{\delta^2{\rm W}}{\sin\alpha}\eta_{xx}, \delta h_3)_\Gamma
 +(\bm{F}_1, \bm{u}^\delta)_\Omega,
\end{align*}
where we used Korn's inequality \eqref{Korn2}. 
Therefore, taking ${\rm R}_0$ sufficiently small so that $4K{\rm R}_0\leq 1$, for $0<{\rm R}\le{\rm R}_0$ 
we obtain the desired energy estimate. $\quad\square$

\bigskip
Note that we can take the tangential and time derivatives of the boundary conditions. 
Applying $\delta\partial_x$, $\delta^2\partial_x^{2}$, and $\delta\partial_t$ to \eqref{NS}--\eqref{BC0} 
and using the above proposition, we obtain 
\begin{align}\label{higher}
& \dfrac{1}{2}\dfrac{{\rm d}}{{\rm d}t}\bigg\{\delta^2\|\bm{u}^\delta_x\|^2
 +\dfrac{2}{{\rm R}}\bigg(\dfrac{1}{\tan\alpha}\delta^2|\eta_x|_0^2
  +\dfrac{\delta^2\rm W}{\sin\alpha}\delta^2|\eta_{xx}|_0^2\bigg)\bigg\}
 +\frac{1}{4K{\rm R}}\delta\|\nabla_\delta\bm{u}^\delta_x\|^2 \\
&\le \frac{4K}{\rm R}(\delta|\eta_x|_0^2+\delta|(b_3\eta)_x|_0^2)
 +\dfrac{1}{\rm R}\delta(h_{1x},u_x)_\Gamma-\dfrac{2}{\rm R}\delta(h_{2x},\delta v_x)_\Gamma \notag \\
&\quad
 +\dfrac{2}{\rm R}\delta(\dfrac{1}{\tan\alpha}\eta_x-\dfrac{\delta^2\rm W}{\sin\alpha}\eta_{xxx},
  \delta h_{3x})_\Gamma+\delta(\bm{F}_{1x}, \bm{u}^\delta_x)_\Omega, \notag
\end{align}
\begin{align}\label{highest}
& \dfrac{1}{2}\dfrac{{\rm d}}{{\rm d}t}\bigg\{\delta^4\|\bm{u}^\delta_{xx}\|^2
 +\dfrac{2}{\rm R}\bigg(\dfrac{1}{\tan\alpha}\delta^4|\eta_{xx}|_0^2
  +\dfrac{\delta^2\rm W}{\sin\alpha}\delta^4|\eta_{xxx}|_0^2\bigg)\bigg\}
 +\frac{1}{4K{\rm R}}\delta^3\|\nabla_\delta\bm{u}^\delta_{xx}\|^2 \\
&\le \frac{4K}{\rm R}(\delta^3|\eta_{xx}|_0^2+\delta^3|(b_3\eta)_{xx}|_0^2)
 +\dfrac{1}{\rm R}\delta^3(h_{1xx},u_{xx})_\Gamma-\dfrac{2}{\rm R}\delta^3(h_{2xx},\delta v_{xx})_\Gamma \notag \\
&\quad
 +\dfrac{2}{\rm R}\delta^3(\dfrac{1}{\tan\alpha}\eta_{xx}-\dfrac{\delta^2{\rm W}}{\sin\alpha}\eta_{xxxx},
  \delta h_{3xx})_\Gamma+\delta^3(\bm{F}_{1xx},\bm{u}^\delta_{xx})_\Omega, \notag
\end{align}
\begin{align}\label{time}
& \dfrac{1}{2}\dfrac{{\rm d}}{{\rm d}t}\bigg\{\delta^2\|\bm{u}^\delta_t\|^2
 +\dfrac{2}{\rm R}\bigg(\dfrac{1}{\tan\alpha}\delta^2|\eta_t|_0^2
  +\dfrac{\delta^2{\rm W}}{\sin\alpha}\delta^2|\eta_{tx}|_0^2\bigg)\bigg\}
 +\frac{1}{4K{\rm R}}\delta\|\nabla_\delta\bm{u}^\delta_t\|^2 \\
&\le \frac{4K}{\rm R}(\delta|\eta_t|_0^2+\delta|(b_3\eta)_t|_0^2)
 +\dfrac{1}{\rm R}\delta(h_{1t},u_t)_\Gamma-\dfrac{2}{\rm R}\delta(h_{2t},\delta v_t)_\Gamma \notag \\
&\quad
 +\dfrac{2}{\rm R}\delta(\dfrac{1}{\tan\alpha}\eta_t-\dfrac{\delta^2\rm W}{\sin\alpha}\eta_{txx},
 \delta h_{3t})_\Gamma \notag \\
&\quad
 +\delta(\bm{f}_t, \bm{u}^\delta_t)_\Omega-\dfrac{2}{\rm R}\delta((A_4\nabla_\delta p)_t,\bm{u}^\delta_t)_\Omega
 +\dfrac{1}{\rm R}\delta((b_2u_{yy})_t,  u_t)_\Omega. \notag
\end{align}
For later use, we will compute 
$-\dfrac{2}{{\rm R}}\delta(\partial_x^k\big(A_4\nabla_\delta p\big)_t, \partial_x^k \bm{u}^\delta_t)_\Omega$ 
for a nonnegative integer $k$. 
Applying $\delta\partial_t$ to the first equation in \eqref{NS}, we have 
\begin{equation}\label{timetime}
\delta^2\bm{u}^\delta_{tt}
= -\dfrac{2}{{\rm R}}\delta(I+A_4)\nabla_\delta p_t
 -\dfrac{2}{{\rm R}}\delta A_{4t}\nabla_\delta p+\delta\bm{F}_{3t},
\end{equation}
where 
\begin{equation}\label{def-F3}
\bm{F}_3
= -(\bm{U}\cdot\nabla_\delta)\bm{u}^\delta-(\bm{u}^\delta\cdot\nabla_\delta)\bm{U}
 +\dfrac{1}{{\rm R}}\bigl(\delta^2\bm{u}^\delta_{xx}+(I+A_3)\bm{u}^\delta_{yy}\bigr)+\bm{f}.
\end{equation}
Moreover, we can rewrite \eqref{NS} as 
\begin{equation}\label{A4p}
\dfrac{2}{{\rm R}} A_4\nabla_\delta p = -\delta A_5\bm{u}^\delta_t+A_5\bm{F}_3,
\end{equation}
where $A_5=A_4(I+A_4)^{-1}$. 
Note that $A_5$ is a symmetric matrix due to the symmetry of $A_4$ (see \eqref{NLp}). 
Applying $\delta\partial_x^k\partial_t$ to the above equation, we have 
\begin{align*}
\dfrac{2}{{\rm R}}\delta\partial_x^k(A_4\nabla_\delta p)_t
&= -\delta^2A_5\partial_x^k\bm{u}^\delta_{tt}-\delta^2\partial_x^k(A_{5t}\bm{u}^\delta_t)
 -\delta^2[\partial_x^k, A_5]\bm{u}^\delta_{tt}+\delta\partial_x^k(A_5\bm{F}_3)_t.
\end{align*}
This together with \eqref{timetime} yields 
\begin{align}\label{commutator0}
-\dfrac{2}{{\rm R}}\delta(\partial_x^k(A_4\nabla_\delta p)_t, \partial_x^k \bm{u}^\delta_t)_\Omega
&= \frac{1}{2}\frac{{\rm d}}{{\rm d}t}\delta^2(A_5\partial_x^k\bm{u}^\delta_t, \partial_x^k\bm{u}^\delta_t)_\Omega \\
&\quad
 +\delta(\partial_x^k\{\frac12\delta A_{5t}\bm{u}^\delta_t-(A_5\bm{F}_3)_t\}, \partial_x^k\bm{u}^\delta_t)_\Omega
 +\delta(\bm{G}_k, \partial_x^k\bm{u}^\delta_t)_\Omega, \notag
\end{align}
where 
\begin{equation}\label{commutator}
\bm{G}_k = [\partial_x^k, A_5]
 \bigg\{-\dfrac{2}{{\rm R}}(I+A_4)\nabla_\delta p_t-\dfrac{2}{{\rm R}} A_{4t}\nabla_\delta p+\bm{F}_{3t}\bigg\}
 +\frac12\delta[\partial_x^k,A_{5t}]\bm{u}_t^{\delta}.
\end{equation}
In Particular, in the case of $k=0$, we have 
\begin{align*}
-\dfrac{2}{{\rm R}}\delta((A_4\nabla_\delta p)_t,  \bm{u}^\delta_t)_\Omega
&= \frac{1}{2}\frac{{\rm d}}{{\rm d}t}\delta^2(A_5\bm{u}^\delta_t, \bm{u}^\delta_t)_\Omega
 +\delta(\frac{1}{2}\delta A_{5t}\bm{u}^\delta_t-(A_5\bm{F}_{3})_t, \bm{u}^\delta_t)_\Omega.
\end{align*}
By substituting this into \eqref{time}, we get 
\begin{align}\label{time2}
& \dfrac{1}{2}\dfrac{{\rm d}}{{\rm d}t}\bigg\{\delta^2((I-A_5)\bm{u}^\delta_t, \bm{u}^\delta_t)_\Omega
 +\dfrac{2}{{\rm R}}\bigg(\dfrac{1}{\tan\alpha}\delta^2|\eta_t|_0^2
  +\dfrac{\delta^2{\rm W}}{\sin\alpha}\delta^2|\eta_{tx}|_0^2\bigg)\bigg\}
 +\frac{1}{4K{\rm R}}\delta\|\nabla_\delta\bm{u}^\delta_t\|^2 \\
&\le \frac{4K}{\rm R}(\delta|\eta_t|_0^2+\delta|(b_3\eta)_t|_0^2)+\dfrac{1}{{\rm R}}\delta(h_{1t},u_t)_\Gamma
 -\dfrac{2}{{\rm R}}\delta(h_{2t},\delta v_t)_\Gamma\notag \\
&\quad
 +\dfrac{2}{{\rm R}}\delta(\dfrac{1}{\tan\alpha}\eta_t-\dfrac{\delta^2{\rm W}}{\sin\alpha}\eta_{txx},
  \delta h_{3t})_\Gamma+\delta(\bm{F}_2, \bm{u}^\delta_t)_\Omega, \notag
\end{align}
where 
\begin{equation}\label{F2}
\bm{F}_2 = \bm{f}_t+\dfrac{1}{{\rm R}}\binom{(b_2 u_{yy})_t}{0}
 +\frac{1}{2}\delta A_{5t}\bm{u}^\delta_t-(A_5\bm{F}_3)_t. 
\end{equation}
Note that $I-A_5$ is positive definite for small solutions.

The lowest order energy obtained in \eqref{lowest} is not appropriate in order to get 
the uniform estimate in $\delta$, which is our goal in this paper. 
We thereby need to modify the lowest energy estimate. 
Now it follows from the first and second equations in \eqref{NS} that 
$$
\delta^2 v_t+\bar{u}\delta^2 v_x+\dfrac{2}{{\rm R}} p_y-\dfrac{1}{{\rm R}}\delta(\delta^2 v_x+ u_y)_x
 -\dfrac{2}{{\rm R}}\delta v_{yy} = f_1,
$$
where 
\begin{equation}\label{NLdomainf2}
f_1 = \bigg(\bm{f}-\frac{2}{{\rm R}} A_4\nabla_\delta p\bigg)\cdot\bm{e}_2. 
\end{equation}
Taking the inner product of $\delta v$ with the above equation, we obtain 
\begin{align*}
\frac{\delta}{2}\frac{{\rm{d}}}{{\rm{d}} t}\delta^2\| v\|^2-\dfrac{2}{{\rm R}}( p,\delta v_y)_\Omega
 +\dfrac{1}{{\rm R}}(\delta^2 v_x+ u_y, \delta^2 v_x)_\Omega
 +\dfrac{2}{{\rm R}}\delta^2\| v_y\|^2+\dfrac{2}{{\rm R}}( p-\delta v_y, \delta v)_\Gamma
= (f_1,\delta v)_\Omega.
\end{align*}
Thus using the second equation in \eqref{NS} and integration by parts in $x$, we have 
\begin{align}\label{enelow}
& \frac{\delta}{2}\frac{{\rm{d}}}{{\rm{d}} t}\delta^2\| v\|^2
 +\dfrac{2}{{\rm R}}( p-\delta v_y, \delta v)_\Gamma+\dfrac{1}{{\rm R}}\delta^4\| v_x\|^2
 +\dfrac{2}{{\rm R}}\delta^2\| v_y\|^2 \\
&=\dfrac{2}{{\rm R}}(\delta p_x,  u)_\Omega+\dfrac{1}{{\rm R}}(\delta u_{xy}, \delta  v)_\Omega
 +(f_1, \delta v)_\Omega. \notag
\end{align}

\begin{lem}\label{enelem1}
The following inequality holds. 
\begin{align*}
& \dfrac{2}{{\rm R}}(\delta p_x, u)_\Omega
 +\frac{1}{3{\rm R}}\bigg(\frac{1}{\tan^2\alpha}\delta^2|\eta_x|_0^2
  +\frac{2\delta^2{\rm W}}{\tan\alpha\sin\alpha}\delta^2|\eta_{xx}|_0^2
  +\frac{(\delta^2{\rm W})^2}{\sin^2\alpha}\delta^2|\eta_{xxx}|_0^2\bigg)
 +\frac{1}{{\rm R}}\delta^2\|\partial_y^{-1} p_x\|^2 \\
&\qquad
 \le I_1+I_2+I_3, 
\end{align*}
where 
$$
\begin{cases}
I_1=-\dfrac{2}{{\rm R}}(\delta\partial_y^{-1} p_x, (2+b_3)\eta)_\Omega, \\[2mm]
I_2=-\dfrac{2}{{\rm R}}(\delta\partial_y^{-1} p_x, 
 -\delta^2 v_x(\cdot, 1)+ h_1+\partial_y^{-1}( u_{yy}-2\delta  p_x))_\Omega, \\[2mm]
I_3=\dfrac{1}{{\rm R}}(2\delta^4| u_{xx}|_0^2+2\delta^2| h_{2x}|_0^2+3\delta^2\|\partial_y^{-2} p_{xy}\|^2).
\end{cases}
$$
\end{lem}

\noindent
{\it Proof}. \ 
By the first equation in \eqref{BC1} and \eqref{BC0}, we see that 
\begin{align}\label{enelem11}
\dfrac{2}{{\rm R}}(\delta p_x,  u)_\Omega
&= -\dfrac{2}{{\rm R}}(\partial_y^{-1}\delta p_x,  u_y)_\Omega 
 = -\dfrac{2}{{\rm R}}(\partial_y^{-1} \delta  p_x,  u_y(\cdot,1)+\partial_y^{-1}  u_{yy})_\Omega \\
&= -\dfrac{2}{{\rm R}}(\partial_y^{-1}\delta  p_x, (2+b_3)\eta-\delta^2 v_x(\cdot,1)+ h_1
  +2\partial_y^{-1} \delta p_x+\partial_y^{-1}( u_{yy}-2\delta p_x))_\Omega \notag \\
&= -\frac{4}{{\rm R}}\delta^2\|\partial_y^{-1}  p_x\|^2+I_1+I_2. \notag
\end{align}
On the other hand, it follows from the second equations in \eqref{NS} and \eqref{BC1} that 
\begin{align}\label{pm}
p(x,y)
&=  p(x, 1)+(\partial_y^{-1} p_y)(x,y) \\
&= -\delta u_x(x,1)+\dfrac{1}{\tan\alpha}\eta-\dfrac{\delta^2{\rm W}}{\sin\alpha}\eta_{xx}
 +h_2+(\partial_y^{-1} p_y)(x,y). \notag
\end{align}
Thus applying $\delta{\rm R}^{-\frac12}\partial_y^{-1}\partial_x$ to the above equation, we obtain 
\begin{align*}
& \frac{y-1}{{\rm R}^{\frac12}}\bigg(\dfrac{1}{\tan\alpha}\delta\eta_x
 -\dfrac{\delta^2{\rm W}}{\sin\alpha}\delta\eta_{xxx}\bigg) \\
&= \frac{\delta}{{\rm R}^{\frac12}}(\partial_y^{-1} p_x)(x,y)
 +\frac{y-1}{{\rm R}^{\frac12}}(\delta^2 u_{xx}(x,1)-\delta h_{2x})
 -\frac{\delta}{{\rm R}^{\frac12}}(\partial_y^{-2} p_{xy})(x,y).
\end{align*}
Squaring both sides of the above equation and integrating the resulting equality on $\Omega$, we have 
$$
\frac{1}{3{\rm R}}\bigg(\frac{1}{\tan^2\alpha}\delta^2|\eta_x|_0^2
 +\frac{2\delta^2{\rm W}}{\tan\alpha\sin\alpha}\delta^2|\eta_{xx}|_0^2
 +\frac{(\delta^2{\rm W})^2}{\sin^2\alpha}\delta^2|\eta_{xxx}|_0^2\bigg)
\le \frac{3}{{\rm R}}\delta^2\|\partial_y^{-1} p_x\|^2+I_3,
$$
where we used integration by parts in $x$.
This and \eqref{enelem11} lead to the desired inequality. 
$\quad\square$

\bigskip
This lemma together with \eqref{nishida} and \eqref{enelow} implies that 
\begin{align}\label{new energy estimate}
& \dfrac{1}{2}\dfrac{{\rm d}}{{\rm d}t}\bigg\{\delta^2\|v\|^2
  +\dfrac{2}{{\rm R}}\bigg(\dfrac{1}{\tan\alpha}|\eta|_0^2
  +\dfrac{\delta^2{\rm W}}{\sin\alpha}|\eta_{x}|_0^2\bigg)\bigg\}  
 +\frac{1}{\rm R}(\delta^3\|v_x\|^2+2\delta\|v_y\|^2+\delta\|\partial_y^{-1}p_x\|^2) \\
&\quad
 +\frac{1}{3{\rm R}}\bigg(\frac{1}{\tan^2\alpha}\delta|\eta_x|_0^2
  +\frac{2\delta^2{\rm W}}{\tan\alpha\sin\alpha}\delta|\eta_{xx}|_0^2
  +\frac{(\delta^2{\rm W})^2}{\sin^2\alpha}\delta|\eta_{xxx}|_0^2\bigg) \notag \\
&\le -\frac{2}{\rm R}(h_2,v)_\Gamma+\frac{1}{\rm R}\delta(u_{xy},v)_\Omega+(f_1,v)_\Omega
 +\frac{2}{\rm R}(\frac{1}{\tan\alpha}\eta-\frac{\delta^2\rm W}{\sin\alpha}\eta_{xx},h_3)_\Gamma \notag \\
&\quad
 +\delta^{-1}(I_1+I_2+I_3). \notag
\end{align}
The first three terms in the right-hand side are estimated as 
$$
-\frac{2}{\rm R}(h_2,v)_\Gamma+\frac{1}{\rm R}\delta(u_{xy},v)_\Omega+(f_1,v)_\Omega
\le \frac{1}{\rm R}\delta\|v_y\|^2+\frac{1}{\rm R}(2\delta^{-1}|h_2|_0^2+\delta\|u_{xy}\|^2)+{\rm R}\|f_1\|^2
$$
and the first term in the right-hand side can be absorbed in the left-hand side of \eqref{new energy estimate}. 
We proceed to estimate $I_1$, $I_2$, and $I_3$. 
By \eqref{pm} and integration by parts in $x$, $I_1$ is rewritten as 
\begin{align}\label{enelow2}
I_1
&= -\dfrac{2}{{\rm R}}(\delta\partial_y^{-1}\bigg(-\delta u_x(\cdot, 1)+\dfrac{1}{\tan\alpha}\eta
  -\dfrac{\delta^2{\rm W}}{\sin\alpha}\eta_{xx}+ h_2+\partial_y^{-1} p_y\bigg)_x, (2+b_3)\eta)_\Omega\\
&= I_{4}+I_5, \notag
\end{align}
where 
\begin{align}
I_4 
&= \dfrac{2}{{\rm R}}((y-1)(-\delta u_x(\cdot, 1)+ h_2)+\partial_y^{-2} p_y, 
  \delta\bigl((2+b_3)\eta\bigr)_x)_\Omega, 
  \notag \\
I_5
&= -\frac{1}{\rm R}(\dfrac{1}{\tan\alpha}\eta-\dfrac{\delta^2{\rm W}}{\sin\alpha}\eta_{xx}, \delta(b_3\eta)_x)_\Gamma. 
  \label{I5}
\end{align}
Here we used identities $(\eta, \eta_x)_\Gamma=(\eta_{xx}, \eta_x)_\Gamma=0$. 
We estimate $I_2$, $I_3$, and $I_4$ as follows.

\begin{lem}\label{enelem2}
There exists a positive constant $C$ independent of $\delta$, ${\rm R}$, ${\rm W}$, and $\alpha$ such that 
the following estimates hold. 
\begin{align*}
|I_2|
&\le \dfrac{1}{2{\rm R}}\delta^2\|\partial_y^{-1}p_x\|^2
 + C\Big\{\dfrac{1}{{\rm R}}(\delta^4\| v_{xy}\|^2+| h_1|_0^2+\delta^4\| u_{xx}\|^2) \\
&\qquad
 +{\rm R}(\delta^2\|u_{ty}\|^2+\delta^2\|u_x\|^2+\delta^2\| v_y\|^2+\|f_2\|^2)\Big\}, \\
|I_3|
&\le C\Bigl\{\dfrac{1}{{\rm R}}(\delta^4\|u_{xxy}\|^2+\delta^2|h_{2x}|_0^2
  +\delta^8\|v_{xxx}\|^2+\delta^4\| v_{xyy}\|^2) \\
&\qquad
 +{\rm R}(\delta^6\| v_{tx}\|^2+\delta^6\| v_{xx}\|^2 +\delta^2\| f_{1x}\|^2)\Bigr\}, \\
|I_{4}
&|\le \frac{1}{6{\rm R}\tan^2\alpha}(\delta^2|\eta_x|_0^2+\delta^2|(b_3\eta)_x|_0^2) \\
&\quad
 +C\Bigl\{\frac{\tan^2\alpha}{\rm R}(\delta^2\|u_{xy}\|^2+\delta^6\|v_{xx}\|^2
  +\delta^2\|v_{yy}\|^2+|h_2|_0^2) \\
&\qquad
 +R\tan^2\alpha(\delta^4\| v_{ty}\|^2+\delta^4\| v_x\|^2+\| f_1\|^2)\Bigr\},
\end{align*}
where 
\begin{equation}\label{NLdomainf3}
f_2 = -\frac{b_2}{1+b_2}\bigg(\delta u_t+\bar{u}\delta u_x+\bar{u}_y\delta v-\dfrac{1}{{\rm R}}\delta^2u_{xx}\bigg)
 -\frac{2b_2}{{\rm R}(1+b_2)}\delta p_x-\dfrac{1}{1+b_2}f_3
\end{equation}
and $f_3=(\bm{f}-\frac{2}{{\rm R}} A_4\nabla_\delta p)\cdot\bm{e}_1$. 
\end{lem}

\noindent
{\it Proof}. \ 
We can easily estimate $I_3$ and $I_4$ by using the second component of the first equation in \eqref{NS} 
so as to eliminate $p_y$. 
As for $I_2$, by the first component of the first equation in \eqref{NS}, we have
$$
\dfrac{1}{{\rm R}}\bigg(u_{yy}-\frac{2}{1+b_2}\delta p_x\bigg)
= \frac{1}{1+b_2}\bigg(\delta u_t+\bar{u}\delta u_x+\bar{u}_y\delta v
 -\dfrac{1}{{\rm R}}\delta^2u_{xx}\bigg)-\frac{1}{1+b_2}f_3.
$$
Substituting the above equation into $I_2$, we easily obtain the desired estimate.
$\quad\square$

\bigskip
Combining \eqref{new energy estimate}, \eqref{enelow2}, and Lemma \ref{enelem2}, we obtain 
\begin{align}\label{original}
& \frac{1}{2}\frac{{\rm{d}}}{{\rm{d}} t}\bigg\{\delta^2\| v\|^2
 +\dfrac{2}{{\rm R}}\bigg(\dfrac{1}{\tan\alpha}|\eta|_0^2
  +\dfrac{\delta^2{\rm W}}{\sin\alpha}|\eta_{x}|_0^2\bigg)\bigg\}
 +\dfrac{1}{{\rm R}}\biggl(\delta\|\bm{u}^\delta_x\|^2+\frac12\delta\|\partial_y^{-1} p_x\|^2\biggr) \\
&\quad
 +\frac{1}{3{\rm R}}\bigg(\frac{1}{2\tan^2\alpha}\delta|\eta_x|_0^2
 +\frac{2\delta^2{\rm W}}{\tan\alpha\sin\alpha}\delta|\eta_{xx}|_0^2
 +\frac{(\delta^2{\rm W})^2}{\sin^2\alpha}\delta|\eta_{xxx}|_0^2\bigg) \notag \\
&\le C_1\Bigl\{\frac{1}{\rm R}\bigl((1+\tan^2\alpha)\delta\|\nabla_{\delta}\bm{u}_x^{\delta}\|^2
  +\delta^3\|\nabla_{\delta}\bm{u}_{xx}^{\delta}\|^2 \notag \\
&\qquad\qquad
  +\delta^{-1}|h_1|_0^2+(1+\tan^2\alpha)\delta^{-1}|h_2|_0^2+\delta|h_{2x}|_0^2\bigr) \notag \\
&\qquad 
 +{\rm R}\bigl(\delta\|\nabla_{\delta}\bm{u}_x^{\delta}\|^2
  +(1+\tan^2\alpha)\delta\|\nabla_{\delta}\bm{u}_t^{\delta}\|^2 \notag \\
&\qquad\qquad
 +(1+\tan^2\alpha)\delta^{-1}\|f_1\|^2+\delta^{-1}\|f_2\|^2+\delta\|f_{1x}\|^2\bigr)\Bigr\} \notag \\
&\quad
 +\frac{2\delta^2{\rm W}}{{\rm R}\sin\alpha}\delta^{-1}|(\eta_{xx},\delta h_3)_{\Gamma}|
 +\frac{1}{6{\rm R}\tan^2\alpha}\delta|(b_3\eta)_x|_0^2+\delta^{-1}I_5, \notag
\end{align}
where we used the second equation in \eqref{NS} and 
$(\eta, h_3)_\Gamma=(\eta, \varepsilon^2\eta^2\eta_x)_\Gamma=0$. 
Here the constant $C_1$ does not depend on $\delta$, ${\rm R}$, ${\rm W}$, nor $\alpha$. 
This is the modified energy estimate. 
In the left-hand side, we have a new term $\delta\|\partial_y^{-1} p_x\|^2$, 
which plays an important role in this paper.  

In view of the energy estimates obtained in this section, we define an energy function $E_0$, 
a dissipation function $F_0$, and a collection of the nonlinear terms $N_0$ by 
\begin{align}
E_0(\eta, \bm{u}^\delta)
&= \delta^2\| v\|^2+\dfrac{2}{\rm R}\biggl(\frac{1}{\tan\alpha}|\eta|_0^2
  +\frac{\delta^2{\rm W}}{\sin\alpha}|\eta_x|_0^2\biggr) \\
&\quad
 +\beta_1\bigg\{\delta^2\|\bm{u}^\delta_x\|^2
  +\frac{2}{\rm R}\biggl(\dfrac{1}{\tan\alpha}\delta^2|\eta_x|_0^2
  +\dfrac{\delta^2{\rm W}}{\sin\alpha}\delta^2|\eta_{xx}|_0^2\biggr)\bigg\} \notag \\
&\quad
 +\beta_2\bigg\{\delta^4\|\bm{u}^\delta_{xx}\|^2
  +\frac{2}{\rm R}\biggl(\dfrac{1}{\tan\alpha}\delta^4|\eta_{xx}|_0^2
  +\dfrac{\delta^2{\rm W}}{\sin\alpha}\delta^4|\eta_{xxx}|_0^2\biggr)\bigg\} \notag \\
&\quad
 +\beta_3\bigg\{\delta^2((I-A_5)\bm{u}^\delta_t, \bm{u}^\delta_t)_\Omega
  +\frac{2}{\rm R}\biggl(\dfrac{1}{\tan\alpha}\delta^2|\eta_t|_0^2
  +\dfrac{\delta^2{\rm W}}{\sin\alpha}\delta^2|\eta_{tx}|_0^2\biggr)\bigg\},\notag
\end{align}
\begin{align}
F_0(\eta, \bm{u}^\delta, p)
&= \frac{1}{2{\rm R}}\biggl(\delta\|\bm{u}^\delta_x\|^2+\frac12\delta\|\partial_y^{-1} p_x\|^2\biggr) \\
&\quad
 +\frac{1}{6{\rm R}}\biggl(\frac{1}{2\tan^2\alpha}\delta|\eta_x|_0^2
  +\frac{2\delta^2{\rm W}}{\tan\alpha\sin\alpha}\delta|\eta_{xx}|_0^2
  +\frac{(\delta^2{\rm W})^2}{\sin^2\alpha}\delta|\eta_{xxx}|_0^2\biggr) \notag \\
&\quad
 +\frac{1}{8K{\rm R}}(\beta_1\delta\|\nabla_\delta\bm{u}^\delta_x\|^2
  +\beta_2\delta^3\|\nabla_\delta\bm{u}^\delta_{xx}\|^2
  +\beta_3\delta\|\nabla_\delta\bm{u}^\delta_t\|^2), \notag
\end{align}
\begin{align}
N_0(Z)
&= \delta^{-1}|h_1|_0^2+\delta^{-1}|h_2|_0^2+\delta|h_{1x}|_0^2+\delta| h_{2x}|_0^2 \\
&\quad
 +\delta|h_3|_0^2+\delta^3| h_{3t}|_0^2+\delta^3| h_{3x}|_0^2+\delta^5| h_{3xx}|_0^2 \notag \\
&\quad
 +\delta^2||D_x|^{\frac12}h_{1x}|_0^2+\delta^2||D_x|^{\frac12}h_{2x}|_0^2
 +\delta|( h_{1t},  u_t)_\Gamma|+\delta|( h_{2t}, \delta v_t)_\Gamma| \notag \\
&\quad
 +\delta|(b_3\eta)_x|_0^2+\delta^3|(b_3\eta)_{xx}|_0^2+\delta|(b_3\eta)_t|_0^2
 +|(\eta,(b_3\eta)_x)_\Gamma|\notag\\
&\quad
+\delta^2{\rm W}\bigl\{\delta^{-1}|(\eta_{xx},\delta h_{3}+\delta(b_3\eta)_x)_\Gamma|+\delta^3|(\eta_{xxxx},\delta h_{3xx})_\Gamma|+\delta|(\eta_{xxt},\delta h_{3t})_\Gamma|\bigr\} \notag \\
&\quad
 +\delta^{-1}\|f_1\|^2+\delta^{-1}\|f_2\|^2+\delta\|f_{1x}\|^2 \notag \\
&\quad
 +\delta|(\bm{F}_{1x}, \bm{u}^\delta_x)_\Omega|+\delta^3|(\bm{F}_{1xx},\bm{u}^\delta_{xx})_\Omega|
  +\delta|(\bm{F}_2, \bm{u}^\delta_t)_\Omega|, \notag
\end{align}
where $Z=(\eta, \bm{u}^\delta, h_1, h_2, h_3, b_3\eta,  f_1, f_2, \bm{F}_1, \bm{F}_2)$ and 
we will determine the constants $\beta_1$, $\beta_2$, and $\beta_3$ later. 
Note that the terms $|(\eta,(b_3\eta)_x)_\Gamma|$ and $(\delta^2{\rm W})\delta^{-1}|(\eta_{xx}, \delta(b_3\eta)_x)_\Gamma|$ come from $I_5$. 
Summarizing our energy estimates, we obtain the following proposition.

\begin{prop}\label{energy0}
Let ${\rm W}_1$ be a positive constant. There exists a positive constant $\alpha_0$ such that if $0<{\rm R_1}\le{\rm R}\le{\rm R}_0$, ${\rm W}_1\le{\rm W}$, and $0<\alpha\le\alpha_0$, 
then the solution $(\eta, u, v, p)$ of \eqref{NS}--\eqref{BC0} satisfies 
$$
\frac{{\rm d}}{{\rm d}t}E_0+F_0\le C_2N_0, 
$$
where ${\rm R}_0$ is the constant in Proposition \ref{highEest} and the constant $C_2({\rm R}_1, {\rm W}_1, \alpha)$ is independent of  $\delta$, ${\rm R}$,  and ${\rm W}$. 
\end{prop}

\noindent
{\it Proof}. \ 
Multiplying \eqref{higher}, \eqref{highest}, and \eqref{time2} by $\beta_1$, $\beta_2$, and $\beta_3$, 
respectively, and adding these and \eqref{original}, we see that 
$$
\frac{{\rm d}}{{\rm d}t}E_0+2F_0\le L+C(N+N_0),
$$
where 
\begin{align*}
L&= \frac{4K}{\rm R}\bigl((\beta_1+3\beta_3)\delta|\eta_x|_0^2+\beta_2\delta^3|\eta_{xx}|_0^2\bigr) 
 +\biggl\{C_1\biggl(\frac{1+\tan^2\alpha}{\rm R}+{\rm R}\biggr)+\frac{12K}{\rm R}\beta_3\biggr\}
  \delta\|\nabla_{\delta}\bm{u}_x^{\delta}\|^2 \\
&\quad
 +\frac{C_1}{\rm R}\delta^3\|\nabla_{\delta}\bm{u}_{xx}^{\delta}\|^2
 +C_1{\rm R}(1+\tan^2\alpha)\delta\|\nabla_{\delta}\bm{u}_t^{\delta}\|^2, \\
N&= \delta|(h_{1x},u_x)_\Gamma|+\delta|(h_{2x},\delta v_x)_\Gamma|
 +\delta|(\eta_x,\delta h_{3x})_\Gamma|+|((\delta^2{\rm W})\delta^{1/2}\eta_{xxx},\delta^{3/2} h_{3x})_\Gamma| \\
&\quad
 +\delta^3|(h_{1xx},u_{xx})_\Gamma|+\delta^3|(h_{2xx},\delta v_{xx})_\Gamma|
 +\delta^3|(\eta_{xx},\delta h_{3xx})_\Gamma| 
 +\delta|(\eta_t,\delta h_{3t})_\Gamma|+\delta^{-1}|I_5|. 
\end{align*}
Here we used $|\eta_t|_0\le |\eta_x|_0+\|u_{xy}\|+|h_3|_0$, 
which comes from the second equation in \eqref{NS}, the third equation in \eqref{BC1}, 
and Poincar\'e's inequality. 
Moreover, it is easy to see that for any $\epsilon>0$ there exists a constant $C_\epsilon>0$ such that 
$N\le \epsilon F_0+C_\epsilon N_0$. 
Therefore, if we take $(\beta_1, \beta_2, \beta_3)$ so that
\begin{equation}\label{coefficients}
\begin{cases}
\dfrac{4K}{\rm R}(\beta_1+3\beta_3) < \dfrac{1}{12{\rm R}\tan^2\alpha}, 
 \quad \dfrac{4K}{\rm R}\beta_2 < \dfrac{{\rm W}}{3{\rm R}\tan\alpha\sin\alpha}, \\[2mm]
C_1\biggl(\dfrac{1+\tan^2\alpha}{\rm R}+{\rm R}\biggr)+\dfrac{12K}{\rm R}\beta_3 < \dfrac{\beta_1}{8K{\rm R}}, 
 \quad \dfrac{C_1}{\rm R} < \dfrac{\beta_2}{8K{\rm R}},
 \quad C_1{\rm R}(1+\tan^2\alpha) < \dfrac{\beta_3}{8K{\rm R}},
\end{cases}
\end{equation}
and if we choose $\epsilon>0$ sufficiently small, then we obtain $L+CN \leq F_0+C_{\epsilon}N_0$. 
Here taking $(\beta_1, \beta_2, \beta_3)$ as
\[
\beta_2:=16KC_1,\quad \beta_3:=16KC_1{\rm R}_0^2(1+\tan^2\alpha),\quad \beta_1:=16K\bigl\{C_1(1+\tan^2\alpha+{\rm R}_0^2)+12K\beta_3\bigr\},
\]
we see that \eqref{coefficients} is equivalent to 
\[
48K(\beta_1+3\beta_3)\tan^2\alpha<1,\quad 12K\beta_2\tan\alpha\sin\alpha<{\rm W}_1.
\] 
Thus there exists a small constant $\alpha_0$ which depends on ${\rm W}_1$ such that 
\eqref{coefficients} is fulfilled and we obtain the desired energy inequality. 
$\quad\square$

\bigskip
Hereafter, $m$ is an integer satisfying $m\ge2$. 
We define a higher order energy and a dissipation functions $E_m$ and $F_m$ 
and a collection of the nonlinear terms $N_m$ by 
\begin{equation}\label{E and F}
E_m = \sum_{k=0}^m E_0(\partial_x^k\eta, \partial_x^k\bm{u}^\delta), \quad 
F_m = \sum_{k=0}^m F_0(\partial_x^k\eta, \partial_x^k\bm{u}^\delta, \partial_x^k p), 
\end{equation}
\begin{equation}\label{def-N}
N_m = \sum_{k=0}^mN_0(\partial_x^k Z)
  +\sum_{k=1}^m\bigl(\delta|(\bm{G}_k, \partial_x^k\bm{u}^\delta_t)_\Omega|
   +|(\partial_x^k\eta,\partial_x^kh_3)_\Gamma|\bigr).
\end{equation}
Here, we note that $\delta|(\bm{G}_k, \partial_x^k\bm{u}^\delta_t)_\Omega|$ is the term appearing in 
\eqref{commutator0} and that $(\eta,h_3)_\Gamma=0$. 
Under an appropriate assumption of the solution, we have the following equivalence uniformly in $\delta$. 
\begin{align*}
E_m
&\simeq |(1+\delta|D_x|)^2\eta|_m^2+\delta^2|\eta_t|_m^2++\delta^2{\rm W}\bigl\{|(1+\delta|D_x|)^2\eta_x|_m^2+\delta^2|\eta_{tx}|_m^2\bigr\} \\
&\quad
 +\delta^2\|(1+|D_x|)^mv\|^2+\delta^2\|(1+|D_x|)^m(1+\delta|D_x|)\bm{u}_x^{\delta}\|^2
 +\delta^2\|(1+|D_x|)^m\bm{u}_t^{\delta}\|^2 \\
&\simeq |\eta|_m^2
 +\delta^2\bigl\{|(\eta_x,\eta_t)|_m^2+\|(1+|D_x|)^m(v,u_x,u_t)\|^2\bigr\} \\
&\quad
 +\delta^4\bigl\{|(\eta_{xx},\eta_{tx})|_m^2+\|(1+|D_x|)^m(v_x,u_{xx},v_t)\|^2\bigr\}+\delta^6\|(1+|D_x|)^mv_{xx}\|^2\\
&\quad
 +\delta^2{\rm W}\bigl\{|\eta_x|_m^2+\delta^2|(\eta_{xx},\eta_{tx})|_m^2+\delta^4|\eta_{xxx}|_m^2\bigr\},
\end{align*}
\begin{align*}
F_m
&\simeq \delta|\eta_x|_m^2+(\delta^2{\rm W})\delta|\eta_{xx}|_m^2+(\delta^2{\rm W})^2\delta|\eta_{xxx}|_m^3+\delta\|(1+|D_x|)^m\partial_y^{-1}p_x\|^2 \\
&\quad
 +\delta\|(1+|D_x|)^m\bm{u}_x^{\delta}\|^2+\delta\|(1+|D_x|)^m(1+\delta|D_x|)\nabla_{\delta}\bm{u}_x^{\delta}\|^2
 +\delta\|(1+|D_x|)^m\nabla_{\delta}\bm{u}_t^{\delta}\|^2 \\
&\simeq \delta\bigl\{|\eta_x|_m^2+\|(1+|D_x|)^m(v_y,u_x,u_{xy},u_{ty},\partial_y^{-1}p_x)\|^2\bigr\} \\
&\quad
 +\delta^3\|(1+|D_x|)^m(v_x,v_{xy},v_{ty},u_{xx},u_{xxy},u_{tx})\|^2 \\
&\quad
 +\delta^5\|(1+|D_x|)^m(v_{xx},v_{xxy},v_{tx},u_{xxx})\|^2+\delta^7\|(1+|D_x|)^mv_{xxx}\|^2\\
&\quad
+(\delta^2{\rm W})\delta|\eta_{xx}|_m^2+(\delta^2{\rm W})^2\delta|\eta_{xxx}|_m^2.
\end{align*}
Applying $\partial_x^k$ to \eqref{NS}--\eqref{BC0}, using Proposition \ref{energy0}, 
and adding the resulting inequalities for $0\le k\le m$, we obtain a higher order energy estimate 
\begin{equation}\label{est-energy}
\frac{{\rm d}}{{\rm d}t}E_m+F_m\le C_2 N_m.
\end{equation}

\section{Estimate for the pressure}
We will use an elliptic estimate for the pressure $p$. 
First, we derive an equation for $p$. 
Applying $\nabla_\delta\cdot$ to the first equation in \eqref{ns} 
and using the second equation in \eqref{ns}, we have 
\begin{align*}
\dfrac{2}{{\rm R}}\Delta_\delta p
&= -\{\varepsilon(\delta u_x)^2+2\delta^2 v_x(\varepsilon u_y+\bar{u}_y)+\varepsilon(\delta v_y)^2\} \\
&= -\varepsilon^{-1}{\rm{tr}}\bigl(\nabla_\delta(\varepsilon\bm{u}^\delta+\bm{U})^{\rm T}\bigr)^2 =: f.
\end{align*}
We transform this by the diffeomorphism $\Phi$ introduced by \eqref{diffeomorphism} and obtain 
\begin{equation}\label{press2}
\nabla_\delta\cdot A_6\nabla_\delta p' = \frac12{\rm R}J(f\circ\Phi)=:g, 
\end{equation}
where $p'=p\circ\Phi$ and $A_6={J} A_2^{\rm T} A_2$. 
On the other hand, by the definition of $f$ and \eqref{u-trans}, we have 
$f\circ\Phi=-\varepsilon^{-1}{\rm tr}\bigl((A_2\nabla_\delta)(\varepsilon A_1\bm{u}'^{\delta}+\bm{U}')^{\rm T}\bigr)^2$, 
where $\bm{u}'^{\delta}$ is defined by \eqref{u-trans} and $\bm{U}'=(U',0)^{\rm T}:=\bm{U}\circ\Phi$. 
Here we see that 
\begin{align*}
(A_2\nabla_\delta)(\varepsilon A_1\bm{u}'^{\delta}+\bm{U}')^{\rm T}
= \binom{\delta\partial_x+a_1\partial_y}{{J^{-1}}\partial_y}
  (\varepsilon{J^{-1}} u'+U', -\varepsilon a_1u'+\varepsilon\delta v')
= \varepsilon F_1 u'_y+F_2,
\end{align*}
where 
\begin{align*}
F_1&:= 
 \begin{pmatrix}
  a_1J^{-1} & -a_1^2 \\
  J^{-2} & -a_1J^{-1}
 \end{pmatrix}, \\
F_2&:= 
 \begin{pmatrix}
  \delta(\varepsilon J^{-1}u')_x+\varepsilon a_1(J^{-1})_y u'
   & \varepsilon\delta(-a_1u'+\delta v')_x-\varepsilon a_1a_{1y} u'+\varepsilon\delta a_1 v'_y \\
  \varepsilon J^{-1}(J^{-1})_y u'+J^{-1} U'_y 
   & -\varepsilon{J^{-1}} a_{1y} u'+\varepsilon\delta J^{-1} v'_y
 \end{pmatrix}. 
\end{align*}
Here, in the above calculation, we used the identity $\delta U_x'+a_1U_y'=0$. 
It follows from $F_1^2=O$ that 
\begin{equation}\label{def-g}
g=-\frac12{\rm R}J\{{\rm tr}(F_1F_2+F_2F_1)u'_y+\varepsilon^{-1}{\rm tr}(F_2^2)\}, 
\end{equation}
where $F_1$ and $F_2$ do not contain $u'_y$.

Next, as for the boundary condition on $\Gamma$, by the second equation in \eqref{BC1}, we obtain 
\begin{equation}\label{pressD}
p'=-\delta u'_x+\dfrac{1}{\tan\alpha}\eta-\dfrac{\delta^2{\rm W}}{\sin\alpha}\eta_{xx}+h_2=:\phi'
 \qquad {\rm on}\quad \Gamma. 
\end{equation}
Moreover, as for the boundary condition on $\Sigma$, taking the trace of the second component of 
the first equation in \eqref{ns} on $\Sigma$, we obtain $(p+\frac12u_x)_y=0$ on $\Sigma$. 
In view of \eqref{u-trans} and \eqref{nd}, this is transformed into 
$$
{J^{-1}}\bigg\{p'+\frac{\delta}{2}({J^{-1}} u')_x+\frac{1}{2} a_1({J^{-1}} u')_y\bigg\}_y=0
 \qquad \mbox{on}\quad \Sigma. 
$$
Recalling $a_1=-yJ^{-1}\varepsilon\delta\tilde{\eta}_x$, $J=1+\varepsilon(y\tilde{\eta})_y$, 
and \eqref{extension0}, we have $a_1|_{y=0}=0$, $a_{1yy}|_{y=0}=0$, and $J_y|_{y=0}=0$, 
so that we obtain $(a_1({J^{-1}} u')_y)_y|_{y=0}=({J^{-1}} a_{1y} u')_y|_{y=0}$. 
Therefore we have 
\begin{equation}\label{pressN}
(p'+g_0)_y=0 \qquad\mbox{on}\quad \Sigma,
\end{equation}
where 
\begin{equation}\label{g0}
g_0=\frac{1}{2}\{\delta({J^{-1}} u')_x+{J^{-1}} a_{1y} u'\}. 
\end{equation}
Summarizing \eqref{press2}, \eqref{pressD}, and \eqref{pressN}, we have 
\begin{equation}\label{divformp}
\left\{
 \begin{array}{lcl}
  \nabla_\delta\cdot A_6 \nabla_\delta p=g & \mbox{in} & \Omega, \\
  p=\phi & \mbox{on} & \Gamma, \\
  (p+g_0)_y=0 & \mbox{on} & \Sigma.
 \end{array}
\right.
\end{equation}
Here we dropped the prime sign in the notation.

We proceed to derive an elliptic estimate for $p$. 
To this end, we will consider the following boundary value problem 
\begin{equation}\label{p-system3}
\left\{
 \begin{array}{lcl}
  \Delta_\delta q=g+\nabla_\delta\cdot\bm{g} & \mbox{in} & \Omega, \\
  q=\psi_1 & \mbox{on} & \Gamma, \\
  q_y=0 & \mbox{on} & \Sigma,
 \end{array}
\right.
\end{equation}
and show the following lemma.

\begin{lem}\label{q-elliptic}
For any $g,\bm{g}\in L^2(\Omega)$ and $\psi_1\in H^{\frac12}(\Gamma)$, 
there exists a unique solution $q\in H^1(\Omega)$ of \eqref{p-system3} satisfying 
$$
\|\nabla_\delta q\|^2\lesssim \|g\|^2+\|\bm{g}\|^2+\delta||D_x|^\frac{1}{2}\psi_1|_0^2.
$$
\end{lem}

\noindent
{\it Proof}. \ 
First, we will construct a solution of the following equation 
\begin{equation}\label{q1}
\Delta_\delta q_1=g+\nabla_\delta\cdot\bm{g} \qquad\mbox{in}\quad \Omega.
\end{equation}
We extend $g$ and $g_1:=\bm{g}\cdot\bm{e}_1$ as even and $4$-periodic functions in $y$ satisfying 
$\int_0^4 g(x, y){\rm d}y=\int_0^4 g_1(x,y){\rm d}y=0$ 
and $g_2:=\bm{g}\cdot\bm{e}_2$ as an odd and $4$-periodic function. 
By these extension and Fourier series expansion in $x$ and $y$, 
we can construct a solution of \eqref{q1} satisfying 
\begin{align}
& q_{1y}(x,0)=0,\label{q10} \\
& \|q_1\|^2+\|\nabla_\delta q_1\|^2\lesssim\|g\|^2+\|\bm{g}\|^2.\label{estq1}
\end{align}
Next, let us seek the solution of \eqref{p-system3} in the form $q=q_1+q_2$, 
where $q_2$ should be the solution of the following boundary value problem 
$$
\left\{
 \begin{array}{lcl}
  \Delta_\delta q_2=0 & \mbox{in} & \Omega, \\
  q_2=\psi_2 & \mbox{on} & \Gamma, \\
  q_{2y}=0 & \mbox{on} & \Sigma,
 \end{array}
\right.
$$
where $\psi_2=\psi_1-q_1|_{y=1}$ and we used \eqref{q10}. 
By Fourier series expansion in $x$, we easily construct a solution of the above problem satisfying 
\begin{equation}\label{estq2}
\|\nabla_\delta q_2\|^2\lesssim\delta||D_x|^\frac{1}{2}\psi_2|_0^2.
\end{equation}
Here, Lemma \ref{trace} yields $\delta||D_x|^\frac{1}{2}q_1|_0^2\lesssim \|q_1\|^2+\|\nabla_\delta q_{1}\|^2$, 
which together with \eqref{estq1} and \eqref{estq2} implies the desired estimate. 
The uniqueness of the solution is well-known, so that the proof is complete. 
$\quad\square$

\bigskip
Now, we rewrite \eqref{divformp} as 
\begin{equation}\label{p-system4}
\left\{
 \begin{array}{lcl}
  \Delta_\delta q=g+\nabla_\delta\cdot (\nabla_\delta g_0-N_6\nabla_\delta p) & \mbox{in} & \Omega, \\
  q=\phi+g_0 & \mbox{on} & \Gamma, \\
  q_y=0 & \mbox{on} & \Sigma,
 \end{array}
\right.
\end{equation}
where $q=p+g_0$ and $N_6$ is a nonlinear part of $A_6$, that is, $A_6=I+N_6$. 
Applying Lemma \ref{q-elliptic} to the above boundary value problem, we have 
\begin{equation}\label{p-elliptic}
\|\nabla_\delta p\|^2\lesssim
 \|g\|^2+\|g_0\|^2+\|\nabla_\delta g_0\|^2+\|N_6\nabla_\delta p\|^2+\delta||D_x|^{\frac12}\phi|_0^2, 
\end{equation}
where we used $\delta||D_x|^{\frac12} g_0|_0^2\lesssim\|g_0\|^2+\|\nabla_\delta g_0\|^2$ 
which comes from Lemma \ref{trace}. 
Differentiating \eqref{p-system4} in $x$ and $t$, likewise we deduce 
\begin{equation}\label{pxpt}
\begin{cases}
\delta\|\nabla_\delta p_x\|^2\lesssim
 \delta\|g_x\|^2+\delta\|g_{0x}\|^2+\delta\|\nabla_\delta g_{0x}\|^2
  +\delta\|(N_6\nabla_\delta p)_x\|^2+\delta^2||D_x|^{\frac12}\phi_x|_0^2,\\
\delta\|\nabla_\delta p_t\|^2\lesssim
 \delta\|g_t\|^2+\delta\|g_{0t}\|^2+\delta\|\nabla_\delta g_{0t}\|^2
  +\delta\|(N_6\nabla_\delta p)_t\|^2+\delta^2||D_x|^{\frac12}\phi_t|_0^2.
\end{cases}
\end{equation}
Here, for the same reason as the modification of the lowest order energy, 
we need to modify \eqref{p-elliptic}, that is, we estimate $\delta^{-1}\|\nabla_\delta p\|^2$ in a different way. 
As for $\delta^{-1}\|p_y\|^2$, by using the second component of the first equation in \eqref{NS}, we see that 
\begin{equation}\label{est-py}
\delta^{-1}\|p_y\|^2\lesssim F_0+\delta^{-1}\|f_1\|^2,
\end{equation}
where $f_1$ is defined by \eqref{NLdomainf2}. 
To estimate $\delta\|p_x\|^2$ in terms of the dissipation function $F_0$, 
we use the term $\delta\|\partial_y^{-1} p_x\|^2$ in the following way. 
We compute 
\begin{align*}
\delta\|p_x\|^2
&= \delta\int\!\!\!\int_\Omega p_x(x,y)
  \biggl(\frac{\partial}{\partial y}\int_0^y p_x(x,z){\rm d}z\biggr){\rm d}x{\rm d}y \\
&= -\delta\int\!\!\!\int_\Omega p_{xy}(x,y)
  \biggl(\int_0^y p_x(x,z){\rm d}z\biggr){\rm d}x{\rm d}y+\delta\int_0^1 p_x(x,1)
  \bigg(\int_0^1 p_x(x, z){\rm d}z\bigg){\rm d}x \\
&\leq \delta\|p_{xy}\|(\|p_x\|+\|\partial_y^{-1}p_x\|)+\delta|p_x|_0\|p_x\|
\end{align*}
so that we have 
\begin{equation}\label{pcal1}
\delta\|p_x\|^2 \lesssim \delta\|p_{xy}\|^2+\delta\|\partial_y^{-1}p_x\|^2+|p_x|_0^2
\end{equation}
Here, it follows from the second equation in \eqref{BC1} that $\delta| p_x|_0^2 \lesssim F_0+\delta| h_{2x}|_0^2$. 
This together with \eqref{est-py} and \eqref{pcal1} yields 
$$
\delta^{-1}\|\nabla_\delta p\|^2\lesssim F_0+\delta\|p_{xy}\|^2+\delta|h_{2x}|_0^2+\delta^{-1}\|f_1\|^2.
$$
This is the modified estimate of $\delta^{-1}\|\nabla_\delta p\|^2$. 

By differentiating \eqref{p-system4} with respect to $x$ and applying the above argument and \eqref{pxpt}, 
we obtain the following lemma.

\begin{lem}\label{est-p}
For $0\le k\le m$ and $1\le l\le m $, we have 
\begin{align}
& \delta^{-1}\|\nabla_\delta\partial_x^k p\|^2
 \lesssim F_m+\delta\|\partial_x^k p_{xy}\|^2+\delta|\partial_x^k h_{2x}|_0^2
  +\delta^{-1}\|\partial_x^k f_1\|^2,\label{p-elliptic-lower} \\
& \delta\|\nabla_\delta\partial_x^k p_x\|^2
 \lesssim \delta\|\partial_x^{k}g_x\|^2+\delta\|\partial_x^{k}g_{0x}\|^2
  +\delta\|\nabla_\delta\partial_x^{k}g_{0x}\|^2\label{p-elliptic-higher} \\
& \phantom{\delta\|\nabla_\delta\partial_x^k p_x\|^2\lesssim}
 +\delta\|\partial_x^{k}(N_6\nabla_\delta p)_x\|^2+\delta^2||D_x|^{k+\frac12}\phi_x|_0^2, \notag \\
& \delta\|\nabla_\delta\partial_x^{l-1} p_t\|^2
 \lesssim \delta\|\partial_x^{l-1}g_t\|^2+\delta\|\partial_x^{l-1}g_{0t}\|^2
  +\delta\|\nabla_\delta\partial_x^{l-1}g_{0t}\|^2\label{p-elliptic-time} \\
& \phantom{\delta\|\nabla_\delta\partial_x^{l-1} p_t\|^2\lesssim}
 +\delta\|\partial_x^{l-1}(N_6\nabla_\delta p)_t\|^2+\delta^2||D_x|^{l-\frac12}\phi_t|_0^2. \notag
\end{align}
\end{lem}

\section{Estimate for nonlinear terms}

We modify the energy and the dissipation functions $E_m$ and $F_m$ defined by \eqref{E and F} as 
\begin{align}
\tilde{E}_m &= E_m+\|(1+|D_x|)^m u\|^2+\|(1+|D_x|)^m u_y\|^2,\label{tildeE} \\
\tilde{F}_m &= F_m+\label{tildeF} 
 \delta|(1+\delta|D_x|)^{\frac52}\eta_t|_m^2+(\delta^2{\rm W})^2\delta^2||D_x|^{\frac72}\eta|_m^2 \\
&\quad
 +\delta^{-1}\|(1+|D_x|)^m(1+\delta|D_x|)\nabla_\delta p\|^2
 +\delta\|(1+|D_x|)^{m-1}\nabla_\delta p_t\|^2. \notag 
\end{align}
We also introduce another energy function $D_m$ by 
\begin{align}
D_m =& |(1+\delta|D_x|)^2\eta|_m^2+\|(1+|D_x|)^m\bm{u}^\delta\|^2\label{Dk}+\|(1+|D_x|)^mD_{\delta}\bm{u}^\delta\|^2\\
&+\|(1+|D_x|)^mD_{\delta}^2\bm{u}^\delta\|^2+(\delta^2{\rm W})^2|(1+\delta|D_x|)\eta_x|_{m+1}^2+(\delta^2{\rm W})\delta^2\|(1+|D_x|)^mv_{xy}\|^2,\notag
\end{align}
which does not include any time derivatives. 
Moreover, we have the following equivalence uniformly in $\delta$. 
\begin{align*}
D_m
&\simeq |\eta|_m^2+\|(1+|D_x|)^m(u,u_y,u_{yy})\|^2 \\
&\quad
 +\delta^2\bigl\{|\eta_x|_m^2+\|(1+|D_x|)^m(v,v_y,u_x,u_{xy},v_{yy})\|^2\bigr\} \\
&\quad
 +\delta^4\bigl\{|\eta_{xx}|_m^2+\|(1+|D_x|)^m(v_x,v_{xy},u_{xx})\|^2\bigr\} \\
&\quad
 +\delta^6\bigl\{|\eta_{xxx}|_m^2+\|(1+|D_x|)^mv_{xx}\|^2\bigr\}\\
&\quad
+\delta^2{\rm W}\bigl\{|\eta_x|_m^2+\delta^2|\eta_{xx}|_m^2+\delta^4|\eta_{xxx}|_m^2+\delta^2\|(1+|D_x|)^m v_{xy}\|^2\bigr\}\\
&\quad
+(\delta^2{\rm W})^2\bigl\{|\eta_{xx}|_m^2+\delta^2|\eta_{xxx}|_m^2\bigr\}.
\end{align*}

Since the proof of nonlinear estimates derived in this section is particular long, we give a guiding principle of the proof.
A goal of Section 5 is to estimate the nonlinear terms in terms of $\tilde{E}_2\tilde{F}_m, \tilde{F}_2\tilde{E}_m$, and $D_2D_m$. 
As for $\tilde{E}_2\tilde{F}_m$, by using a smallness of the energy this term can be absorbed in the right-hand side of the energy inequality \eqref{est-energy}.
As for $\tilde{F}_2\tilde{E}_m$,  using a boundedness of $\int_0^t \tilde{F}_2(\tau){\rm d}\tau$ and a standard Gronwall's inequality we can estimate this term.
As for $D_2D_m$, we use this estimate in order to estimate an initial energy $E(0)$.
Here, what we should be careful is that if we  use the Sobolev embedding theorem in $\Omega$, that is, $\|u\|_{L^\infty}\lesssim \|u\|_{H^2}$ and   Poincar\'e's inequality for $\eta$, that is, $|\eta|_{L^\infty}\lesssim |\eta_x|_0$, we cannot obtain uniform estimates in $\delta$. Therefore,  we have to estimate nonlinear terms carefully with  fundamental inequalities  which are described below. 

%

\begin{lem}\label{LI0}
If $f(x,0)=0$, then we have 
$$
\|f\|_{L^{\infty}} \lesssim \|f_y\|+\|f_{xy}\|. 
$$
\end{lem}

\noindent
{\it Proof}. \ 
By the Sobolev embedding theorem, we see that 
\begin{align*}
|f(x,y)|^2
&= |f(x,y)-f(x,0)|^2 = \Bigl|\int_0^y f_y(x,y){\rm d}y \Bigr|^2 \\
&\le \int_0^1 |f_y(x,y)|^2{\rm d}y
 \lesssim \int_0^1 \bigl(\|f_y(\cdot,y)\|_{L^2(\mathbb{G})}^2+\|f_{xy}(\cdot,y)\|_{L^2(\mathbb{G})}^2\bigr){\rm d}y,
\end{align*}
which is the desired inequality. 
$\quad\square$

\begin{lem}\label{NLlem}
For any integer $k\ge0$, we have 
\begin{align*}
& \|\partial_x^k(af)\|
 \lesssim \|a\|_{L^{\infty}}\|\partial_x^k f\|+(\|\partial_x^k a\|+\|\partial_x^k a_y\|)(\|f\|+\|f_x\|), \\
& \|\partial_x^k(abf)\|
 \lesssim \|a\|_{L^{\infty}}\|b\|_{L^{\infty}}\|\partial_x^k f\|+\|b\|_{L^{\infty}}(\|\partial_x^k a\|
  +\|\partial_x^k a_y\|)(\|f\|+\|f_x\|) \\
&\qquad\qquad\qquad
 +\|a\|_{L^{\infty}}(\|\partial_x^k b\|+\|\partial_x^k b_y\|)(\|f\|+\|f_x\|). \notag
\end{align*} 
\end{lem}

\noindent
{\it Proof}. \ 
By the well-known inequality 
$$
\|\partial_x^k(af)(\cdot, y)\|_{L^2(\mathbb{G})}
 \lesssim \|a(\cdot, y)\|_{{L^{\infty}}(\mathbb{G})}\|\partial_x^k f(\cdot, y)\|_{L^2(\mathbb{G})}
  +\|f(\cdot, y)\|_{{L^{\infty}}(\mathbb{G})}\|\partial_x^k a(\cdot, y)\|_{L^2(\mathbb{G})},
$$
and the Sobolev embedding theorem, we see that 
\begin{align*}
\|\partial_x^k(af)\|^2
&\lesssim \int_0^1\bigl(\|a(\cdot, y)\|_{{L^{\infty}}(\mathbb{G})}^2\|\partial_x^k f(\cdot, y)\|_{L^2(\mathbb{G})}^2
 +\|f(\cdot, y)\|_{{L^{\infty}}(\mathbb{G})}^2\|\partial_x^k a(\cdot, y)\|_{L^2(\mathbb{G})}^2\bigr){\rm d}y \\
&\lesssim \|a\|_{L^{\infty}}^2\|\partial_x^k f\|^2
 +\sup_{y\in(0,1)}\|\partial_x^k a(\cdot, y)\|_{L^2(\mathbb{G})}^2
  \int_0^1\|f(\cdot, y)\|_{{L^{\infty}}(\mathbb{G})}^2{\rm d}y \\
&\lesssim \|a\|_{L^{\infty}}^2\|\partial_x^k f\|^2+(\|\partial_x^k a\|^2+\|\partial_x^k a_y\|^2)(\|f\|^2+\|f_x\|^2).
\end{align*}
We can prove the second inequality in a similar way. 
$\quad\square$

\begin{lem}\label{commu}
For any integer $k\ge1$, we have 
$$
\|[\partial_x^k, a]f\| 
 \lesssim \|a_x\|_{L^{\infty}}\|\partial_x^{k-1}f\|+(\|\partial_x^k a\|+\|\partial_x^k a_y\|)(\|f\|+\|f_x\|). 
$$
\end{lem}

\noindent
{\it Proof}. \ 
In view of the well-known inequality 
\begin{align*}
\|[\partial_x^k,a]f(\cdot, y)\|_{L^2(\mathbb{G})}
\lesssim \|a(\cdot, y)\|_{{L^{\infty}}(\mathbb{G})}\|\partial_x^{k-1}f(\cdot,y)\|_{L^2(\mathbb{G})}
+\|f(\cdot,y)\|_{{L^{\infty}}(\mathbb{G})}\|\partial_x^k a(\cdot,y)\|_{L^2(\mathbb{G})}, 
\end{align*}
the desired inequality follows in a similar way as the proof of Lemma \ref{NLlem}. 
$\quad\square$

\bigskip

Throughout this section, we assume that 
\begin{equation}\label{smallness}
\tilde{E}_2(t)\le c_1 \quad {\rm for}\quad t\in[0, T/\varepsilon],
\end{equation}
where $T$ and $c_1$ will be determine later. 
We also assume that $(\eta, u, v, p)$ is a solution of \eqref{NS}--\eqref{BC0}, 
$0<\delta, \varepsilon\le1$, ${\rm W_1}\le{\rm W}\le\delta^{-2}{\rm W_2}$, $k$ and $l$ are integers satisfying $0\le k\le m$ and $1\le l\le m$. 
Moreover, we will generally denote smooth functions of $\bm{f}$ by the same symbol $\Phi=\Phi(\bm{f})$ 
and $\Phi_0$ is such a function satisfying $\Phi_0(\bm{0})=0$. 
We also use such a function $\Phi_0$ depending also on $y\in[0,1]$ and denote it by $\Phi_0(\bm{f};y)$, that is, 
$\Phi_0(\bm{0};y)\equiv 0$. 

We prepare several lemmas to proceed nonlinear estimates.

\begin{lem}\label{LI}
The following estimates hold. 
\begin{align}
& \|\tilde{\eta}\|_{L^{\infty}}^2 \lesssim \min\{\tilde{E}_2, D_2\}, 
 \quad \|D_\delta^i\tilde{\eta}\|_{L^{\infty}}^2 \lesssim 
  \min\{\delta^2\tilde{E}_2, \delta^2D_2, \delta \tilde{F}_2\}\quad\mbox{for}\quad 1\leq i\leq 4,\label{LI-E} \\
& \delta^2{\rm W}\|D_\delta^i\tilde{\eta}_x\|_{L^\infty}^2+\delta^4{\rm W}\|D_\delta^i\tilde{\eta}_{xx}\|_{L^\infty}^2\lesssim\min\{\tilde{E}_2, D_2, \delta\tilde{F}_2\}\quad\mbox{for}\quad i=0,1, \label{WL}\\
& 
\begin{cases}
 \|\partial_x^i\bm{u}^\delta\|_{L^\infty}^2 \lesssim \min\{\tilde{E}_2, D_2\}, 
  \quad \delta\|\partial_x^i\bm{u}_x^\delta\|_{L^{\infty}}^2+\delta\|\partial_x^i\bm{u}_t^\delta\|_{L^{\infty}}^2
  \lesssim \tilde{F}_2 \quad\mbox{for}\quad i=0,1, \\
 \delta^4\|v_{xx}\|_{L^{\infty}}^2 \lesssim \min\{\tilde{E}_2, D_2\}, 
\end{cases} \label{LI-u} \\
& 
\begin{cases}
 \|\tilde{\eta}_t\|_{L^{\infty}}^2 \lesssim \min\{\tilde{E}_2,D_2\}, \quad 
  \delta\|\tilde{\eta}_t\|_{L^{\infty}}^2 \lesssim \tilde{F}_2, \\
 \|D_\delta^i\tilde{\eta}_t\|_{L^{\infty}}^2 \lesssim \min\{\tilde{E}_2, D_2, \delta\tilde{F}_2\}
  \quad\mbox{for}\quad i=1, 2, 
  \qquad \|D_\delta^3\tilde{\eta}_t\|_{L^{\infty}}^2 \lesssim \delta\tilde{F}_2,
 \end{cases} \label{LI-v} \\
& \delta\|D_\delta^i\tilde{\eta}_{tt}\|_{L^{\infty}}^2 \lesssim \tilde{F}_2, \quad\mbox{for}\quad i=0,1.
  \label{est-vt}
\end{align}
In particular, we have 
\begin{equation}\label{LI2}
 \begin{cases}
  \|(\tilde{\eta}, \tilde{\eta}_t, D_\delta\tilde{\eta}_t, D_\delta^2\tilde{\eta}_t,
   \bm{u}^\delta, \bm{u}^\delta_x)\|_{L^{\infty}}^2 \lesssim \min\{\tilde{E}_2, D_2\}, \\
  \|(D_\delta\tilde{\eta}, D_\delta^2\tilde{\eta}, D_\delta^3\tilde{\eta})\|_{L^{\infty}}^2
   \lesssim \delta\min\{\tilde{E}_2, D_2\}.
 \end{cases}
\end{equation}
\end{lem}

\noindent
{\bf{Remark 5.1.}} \ 
Using \eqref{LI-E} and taking $c_1$ sufficiently small, we see that $J=1+\varepsilon(y\tilde{\eta})_y$ 
and $I-A_5$ are positive definite.

\bigskip
\noindent
{\it Proof}. \ 
By \eqref{extension4}  in Lemma \ref{extension}, we have 
\begin{align*}
& \|\tilde{\eta}\|_{L^{\infty}}^2\lesssim |\eta|_1^2\lesssim \min\{\tilde{E}_2, D_2\} \\
& \|D_\delta^i\tilde{\eta}\|_{L^\infty}^2 \lesssim \delta^{2i}|\partial_x^{i}\eta|_1^2
 \lesssim \min\{\delta^2\tilde{E}_2, \delta^2D_2, \delta \tilde{F}_2\}
 \quad\mbox{for}\quad 1\le i\le 4.
\end{align*}
Thus \eqref{LI-E} holds. Similarly, we obtain \eqref{WL}. 
\eqref{LI-u} is obtained from Lemma \ref{LI0} and the second equation in \eqref{NS}. 
By the second equation in \eqref{NS}, we have $|v|_1 \lesssim \|(1+|D_x|)v_y\|=\|(1+|D_x|)u_x\|$. 
In view of the assumption \eqref{smallness}, we have 
$$
|\partial_x^jh_3|_0 \lesssim \varepsilon^2|\eta|_{L^{\infty}}^2|\partial_x^{j+1}\eta|_0
 \lesssim |\partial_x^{j+1}\eta|_0 \quad\mbox{for}\quad j\ge0.
$$
Therefore, by \eqref{extension4}  in Lemma \ref{extension} 
and the third equation in \eqref{BC1}, we see that 
\begin{align*}
& \|\tilde{\eta}_t\|_{L^{\infty}} \lesssim |\eta_t|_1 \lesssim |v|_1+|\eta_x|_1+|h_3|_1
 \lesssim \|(1+|D_x|)u_x\|+|\eta_x|_1, \\
& \|D_\delta^i\tilde{\eta}_t\|_{L^{\infty}} \lesssim \delta^i|\partial_x^{i}\eta_t|_1
 \lesssim \delta^i(|\partial_x^{i}v|_1+|\partial_x^{i+1}\eta|_1+|\partial_x^{i}h_3|_1)
 \lesssim \delta^i(\|(1+|D_x|)\partial_x^{i+1}u\|+|\partial_x^{i+1}\eta|_1).
\end{align*}
These estimates give \eqref{LI-v}. 
Similarly, we see that 
\begin{align*}
& \|\tilde{\eta}_{tt}\|_{L^{\infty}} \lesssim |\eta_{tt}|_1
 \lesssim |v_t|_1+|\eta_{tx}|_1+|h_{3t}|_1 \lesssim \|(1+|D_x|)u_{tx}\|+|\eta_t|_2, \\
& \|D_{\delta}\tilde{\eta}_{tt}\|_{L^{\infty}} \lesssim \delta|\partial_x\eta_{tt}|_1
 \lesssim \delta(|\partial_xv_t|_1+|\partial_x^2\eta_t|_1+|\partial_xh_{3t}|_1)
 \lesssim \delta(\|(1+|D_x|)\partial_xu_{tx}\|+|\eta_{tx}|_2+|\eta_x|_2).
\end{align*}
Here, we used $|\eta_t|_{L^{\infty}}\leq\|\tilde{\eta}_t\|_{L^{\infty}}\lesssim \sqrt{\tilde{E}_2}$ and 
the assumption \eqref{smallness}. 
Thus \eqref{est-vt} holds. 
The proof is complete. 
$\quad\square$

\begin{lem}\label{L2estimates}
The following estimates hold. 
\begin{align}
& \|\partial_x^k\tilde{\eta}\|^2\lesssim\min\{\tilde{E}_m, D_m\}, 
 \qquad \|\partial_x^kD_\delta^i\tilde{\eta}\|^2
  \lesssim \min\{\tilde{E}_m, D_m, \delta\tilde{F}_m\} \quad\mbox{for}\quad i=1,2,3, \label{dxkE} \\
&\delta^2{\rm W}\|\partial_x^kD_\delta^i\tilde{\eta}_x\|^2
  \lesssim \min\{\tilde{E}_m, D_m\} \quad\mbox{for}\quad i=1,2\label{Wk}\\
& \|\partial_x^kD_\delta^4\tilde{\eta}\|^2 \lesssim \delta\tilde{F}_m, \label{est-highE} \\
& \delta^2\|\partial_x^kD_\delta^i \tilde{\eta}_t\|^2
  \lesssim \min\{\tilde{E}_m, D_m, \delta \tilde{F}_m\} \quad\mbox{for}\quad i=0,1,2,
 \qquad \delta\|\partial_x^kD_\delta^3 \tilde{\eta}_t\|^2 \lesssim \tilde{F}_m, \label{vtilde} \\
& \delta^3\|\partial_x^kD_\delta^i\tilde{\eta}_{tt}\|^2 \lesssim \tilde{F}_m \quad\mbox{for}\quad i=0,1.
 \label{second time derivative}
\end{align}
\end{lem}

\noindent
{\it Proof}. \ 
By \eqref{extension2} and \eqref{extension1} in Lemma \ref{extension}, we have 
\begin{align*}
& \|\partial_x^kD_\delta^i\tilde{\eta}\|\lesssim\delta^i|\partial_x^{k+i}\eta|_0
 \quad\mbox{for}\quad i\ge0, \\
& \|\partial_x^kD_{\delta}^4\tilde{\eta}\| \lesssim \delta^{\frac72}||D_x|^{k+\frac72}\eta|_0,
\end{align*}
which give \eqref{dxkE} and \eqref{est-highE}, respectively. Similarly, we obtain \eqref{Wk}.
By \eqref{extension2} in Lemma \ref{extension} and a similar argument in the proof of Lemma \ref{LI}, 
we see that 
\begin{align*}
\|\partial_x^kD_\delta^i \tilde{\eta}_t\|
&\lesssim \delta^i|\partial_x^{k+i}\eta_t|_0
 \lesssim \delta^i(\|\partial_x^{k+i+1}u\|+|\partial_x^{k+i+1}\eta|_0) \\
&\lesssim \delta^i\bigl(\|(1+|D_x|)^m\partial_x^{i+1}u\|+|\partial_x^{i+1}\eta|_m\bigr).
\end{align*}
By \eqref{extension1} in Lemma \ref{extension}, Lemma \ref{trace}, Poincar\'e's inequality, and the estimate 
\begin{equation}\label{estimate of h3}
||D_x|^{k+\frac52}h_3|_0 \lesssim |(\eta,\eta_x)|_{L^{\infty}}^2||D_x|^{k+\frac72}\eta|_0
 \lesssim ||D_x|^{k+\frac72}\eta|_0,
\end{equation}
we see that 
\begin{align*}
\|\partial_x^kD_\delta^3 \tilde{\eta}_t\|
&\lesssim \delta^{\frac52}||D_x|^{\frac12}\partial_x^{k+2}\eta_t|_0
 \le \delta^{\frac52}(||D_x|^{\frac12}\partial_x^{k+2}v|_0+||D_x|^{\frac12}\partial_x^{k+2}\eta_x|_0
  +||D_x|^{\frac12}\partial_x^{k+2}h_3|_0) \\
&\lesssim \delta^3\|\partial_x^kv_{xxx}\|+\delta^2\|\partial_x^kv_{xxy}\|
 +\delta^{\frac52}||D_x|^{k+\frac72}\eta|_0 \\
&\lesssim \delta^3\|(1+|D_x|)^mv_{xxx}\|+\delta^2\|(1+|D_x|)^mv_{xxy}\|
 +\delta^{\frac52}||D_x|^{k+\frac72}\eta|_0. 
\end{align*}
These estimates give \eqref{vtilde}. 
It is easy to see that 
$$
|\partial_x^jh_{3t}|_0 \lesssim \varepsilon^2(
 |\eta|_{L^{\infty}}^2|\partial_x^{j+1}\eta_t|_0+|\eta|_{L^{\infty}}|\eta_t|_{L^{\infty}}|\partial_x^{j+1}\eta|_0)
 \lesssim |\partial_x^{j+1}\eta_t|_0+|\partial_x^{j+1}\eta|_0
  \quad\mbox{for}\quad j\ge0.
$$
Therefore, by \eqref{extension2} in Lemma \ref{extension} and the third equation in \eqref{BC1}, we see that 
\begin{align*}
\|\partial_x^k\tilde{\eta}_{tt}\| 
&\lesssim |\partial_x^k\eta_{tt}|_0 \lesssim \|(1+|D_x|)^mu_{tx}\|_0+|\eta_{tx}|_m+|h_{3t}|_m \\
&\lesssim \|(1+|D_x|)^mu_{tx}\|_0+|\eta_{tx}|_m+|\eta_x|_m. 
\end{align*}
Similarly, by \eqref{extension1} in Lemma \ref{extension} we obtain 
\begin{align*}
\|\partial_x^kD_\delta\tilde{\eta}_{tt}\| 
&\lesssim \delta^{\frac12}||D_x|^{\frac12}\partial_x^k\eta_{tt}|_0
 \lesssim \delta^{\frac12}(||D_x|^{\frac12}\partial_x^kv_t|_0+||D_x|^{\frac12}\partial_x^{k+1}\eta_t|_0
  +||D_x|^{\frac12}\partial_x^kh_{3t}|_0) \\
&\lesssim \delta\|\partial_x^kv_{tx}\|+\|\partial_x^kv_{ty}\|
 +\delta|\partial_x^k\eta_{txx}|_0+|\partial_x^k\eta_{tx}|_0
 +\delta|\partial_x^{k+1}h_{3t}|_0+|\partial_x^kh_{3t}|_0 \\
&\lesssim \delta\|(1+|D_x|)^mv_{tx}\|+\|(1+|D_x|)^mv_{ty}\|
 +\delta|\eta_{txx}|_m+|(\eta_{xx},\eta_{tx},\eta_x)|_m. 
\end{align*}
These estimates give \eqref{second time derivative}. 
The proof is complete. 
$\quad\square$

\begin{lem}\label{additonal lemma}
The following estimates hold. 
\begin{align}
& \delta^{2i+1}||D_x|^{i+\frac12}\eta|_m^2 \lesssim \min\{\tilde{E}_m, D_m\},
 \quad \delta^{2i+2}||D_x|^{i+\frac12}\eta_x|_m^2 \lesssim \tilde{F}_m
  \quad\mbox{for}\quad i=0, 1, 2, \label{interpo} \\
& \delta^{3}{\rm W}||D_x|^{\frac12}\eta_x|_m^2 \lesssim \min\{\tilde{E}_m, D_m\},
 \quad \delta^{2i+4}{\rm W}||D_x|^{i+\frac12}\eta_{xx}|_m^2 \lesssim \tilde{F}_m
  \quad\mbox{for}\quad i=0, 1,\label{interpoW}\\
& \delta|\bm{u}^\delta|_{m+\frac12}^2 \lesssim \min\{\tilde{E}_m, D_m\},
 \quad \delta^{3}|\bm{u}^\delta_x|_{m+\frac12}^2 \lesssim \min\{D_m, \delta \tilde{F}_m\}, 
 \quad \delta^{5}|\bm{u}^\delta_{xx}|_{m+\frac12}^2 \lesssim \delta \tilde{F}_m, \label{fracu} \\
& |\bm{u}^\delta|_m^2 \lesssim \min\{\tilde{E}_m, D_m\},
 \quad \delta^{2i}|\partial_x^i\bm{u}^\delta|_m^2
  \lesssim \delta \tilde{F}_m \quad\mbox{for}\quad i=1,2, \label{gammaku} \\
& \delta^2|\bm{u}^\delta_t|^2_{m+\frac12} \lesssim \tilde{F}_m. \label{gammaut}
\end{align}
\end{lem}

\noindent
{\it Proof}. \ 
By an interpolation inequality, we have 
$\delta^{2i+1}||D_x|^{i+\frac12}\eta|_m^2 \lesssim 
 \delta^{2i}|\partial_x^{i}\eta|_m^2+\delta^{2i+2}|\partial_x^i\eta_x|_m^2$, 
which gives the first estimate in \eqref{interpo}. 
Similarly, we can show the second estimate in \eqref{interpo} for $i=0,1$, 
and the case $i=2$ follows directly from the definition of $\tilde{F}_m$. Likewise, we obtain \eqref{interpoW}.
By Lemma \ref{trace} and Poincar$\acute{\rm e}$'s inequality, we see that 
\begin{align*}
\delta^{2i+1}|\partial_x^i\bm{u}^\delta|_{m+\frac12}^2
&\lesssim \delta^{2i+1}|\partial_x^i\bm{u}^\delta|_0^2
 +\delta^{2i+1}||D_x|^{\frac12}\partial_x^{m+i}\bm{u}^\delta|_0^2 \\
&\lesssim \delta^{2i+1}\|\partial_x^i\bm{u}_y^{\delta}\|^2
 +\delta^{2(i+1)}\|\partial_x^{m+i}\bm{u}^\delta_x\|^2+\delta^{2i}\|\partial_x^{m+i}\bm{u}^\delta_y\|^2
\end{align*}
for $i\ge0$, which leads to \eqref{fracu}. 
Similarly, we can show \eqref{gammaut}. 
Poincar$\acute{\rm e}$'s inequality and the second equation in \eqref{NS} yield \eqref{gammaku}. 
The proof is complete. 
$\quad\square$

\bigskip
In view of Lemmas \ref{LI} and \ref{L2estimates} and the inequality 
$\|\partial_x^k\Phi_0(\bm{f};y)\|\le C(\|\bm{f}\|_{L^{\infty}})\|\partial_x^k\bm{f}\|$, 
we obtain the following lemma.

\begin{lem}\label{NLlem2}
For $j=0,1$, the following estimates hold. 
\begin{align}
& \|\Phi_0(\tilde{\eta}, D_\delta\tilde{\eta}, D_\delta^2\tilde{\eta}, D_\delta^3\tilde{\eta},
 \delta\tilde{\eta}_t, \delta D_\delta\tilde{\eta}_t, \delta D_\delta^2\tilde{\eta}_t,
 \bm{u}^\delta, \delta\bm{u}^\delta_x, \delta^3v_{xx};y)\|^2_{L^{\infty}}
  \lesssim \min\{\tilde{E}_2, D_2\}, \label{NLlem2-1} \\ 
& \|\partial_x^k\Phi_0(\tilde{\eta}, D_\delta\tilde{\eta}, D_\delta^2\tilde{\eta}, D_\delta^3\tilde{\eta}, 
 \delta\tilde{\eta}_t, \delta D_\delta\tilde{\eta}_t, \delta D_\delta^2\tilde{\eta}_t, 
 \bm{u}^\delta, \delta\bm{u}^\delta_x,\delta^3v_{xx};y)\|^2 \lesssim \min\{\tilde{E}_m, D_m\}, \label{NLlem2-3} \\ 
& \|\partial_x^k\partial_y^j\Phi_0(\tilde{\eta}, D_\delta\tilde{\eta}, D_\delta^2\tilde{\eta},
 \delta\tilde{\eta}_t, \delta D_\delta\tilde{\eta}_t, 
 \bm{u}^\delta, \delta^2v_x;y)\|^2 \lesssim \min\{\tilde{E}_m, D_m\} \label{NLlem2-5} \\
& \delta\|\partial_x^l\partial_y^j\Phi_0(\tilde{\eta}, D_\delta\tilde{\eta}, D_\delta^2\tilde{\eta}, 
 D_\delta^3\tilde{\eta}, \delta\tilde{\eta}_t, \delta D_\delta\tilde{\eta}_t, \delta D_\delta^2\tilde{\eta}_t, 
 \bm{u}^\delta, \delta\bm{u}^\delta_x;y)\|^2 \lesssim \tilde{F}_m. \label{NLlem2-4}
\end{align}
\end{lem}

\noindent
{\bf Remark 5.2.} \ 
As for \eqref{NLlem2-4}, if $\Phi_0$ does not contain $\tilde{\eta}$ and $u$, 
then $\delta$ appearing in the coefficient of the term $\|\partial_x^l\partial_y^j\Phi_0\|^2$ 
is unnecessary and we can replace $l$ with $k$.

\bigskip
This lemma together with Lemma \ref{trace} gives the following lemma.

\begin{lem}\label{NLlem3}
The following estimates hold. 
\begin{align}
& |\Phi_0(\eta, \delta\eta_x, \delta^2\eta_{xx}, \bm{u}^\delta|_\Gamma, \delta^2 v_x|_\Gamma)|^2_{L^{\infty}}
 \lesssim \min\{\tilde{E}_2, D_2\}, \label{NLlem3-1} \\
& \delta|\Phi_0(\eta, \delta\eta_x, \delta^2\eta_{xx}, \bm{u}^\delta|_\Gamma, \delta^2 v_x|_\Gamma)|^2_{m+\frac12}
 \lesssim \min\{\tilde{E}_m, D_m\}, \label{NLlem3-3} \\
& |\Phi_0(\eta, \delta\eta_x, \delta^2\eta_{xx}, \bm{u}^\delta|_\Gamma, \delta^2 v_x|_\Gamma)|^2_m
 \lesssim \min\{\tilde{E}_m, D_m\}. \label{NLlem3-5}
\end{align}
\end{lem}

By \eqref{WL} in Lemma \ref{LI}, \eqref{Wk} in Lemma \ref{L2estimates} and Lemma \ref{trace}, we obtain the following lemma.

\begin{lem}\label{NLlemW}
The following estimates hold. 
\begin{align}
& {\rm W}|\Phi_0(\delta\eta_x, \delta^2\eta_{xx})|^2_{L^{\infty}}
 \lesssim \min\{\tilde{E}_2, D_2\}, \label{NLlemW-1} \\
& \delta{\rm W}|\Phi_0(\delta\eta_x, \delta^2\eta_{xx})|^2_{m+\frac12}
 \lesssim \min\{\tilde{E}_m, D_m\}, \label{NLlemW-3} \\
& {\rm W}|\Phi_0(\delta\eta_x, \delta^2\eta_{xx})|^2_m
 \lesssim \min\{\tilde{E}_m, D_m\}. \label{NLlemW-5}
\end{align}
\end{lem}

We set
\begin{equation}\label{def-w}
W = (w_1,\ldots,w_7)
 := (D_\delta\tilde{\eta}, D_\delta^2\tilde{\eta}, \delta\tilde{\eta}_t, 
  \delta D_\delta\tilde{\eta}_t, D_\delta^3\tilde{\eta}, \delta D_\delta^2\tilde{\eta}_t, \delta\bm{u}^\delta_x). 
\end{equation}

\begin{lem}\label{NLlemw}
For $j=0,1$, the following estimates hold. 
\begin{align}
& \delta^{-1}\|w_\lambda\|_{L^{\infty}}^2 \lesssim \min\{\delta\tilde{E}_2, \tilde{F}_2\}
 \makebox[-8em]{}& \mbox{for}\quad 1\le \lambda\le7, \label{w-1} \\
& \delta^{-1}\|\partial_x^k\partial_y^j w_\lambda\|^2 \lesssim \tilde{F}_m
 \makebox[-8em]{}&\mbox{for}\quad 1\le \lambda\le7, \label{w-2} \\
& \delta^{-2}\|\partial_x^{l-1}\partial_y^j w_\lambda\|^2 \lesssim \tilde{E}_m
 \makebox[-8em]{}&\mbox{for}\quad  1\le \lambda\le4. \label{w-3}
\end{align}
\end{lem}

\noindent
{\it Proof}. \ 
\eqref{w-1} and \eqref{w-2} follow from Lemmas \ref{LI} and \ref{LI-E}, respectively. 
In the same way as the proof of Lemma \ref{LI-E}, we can show \eqref{w-3}. 
$\quad\square$

\bigskip
We begin to estimate the nonlinear terms. 
First, we will estimate $h_1$, $h_2$, $h_3$, and $b_3\eta$. 
By the explicit form of $h_1$ defined by \eqref{NLtan}, $h_1$ is consist of terms in the form 
\begin{equation}\label{expl-h}
\begin{cases}
\Phi_0(\varepsilon\eta, \varepsilon\delta\eta_x, \varepsilon\bm{u}^\delta|_\Gamma)
 \delta^{i}\partial_x^i\eta \qquad\mbox{for}\quad i=1,2, \\
\Phi_0(\varepsilon\eta,\varepsilon\delta\eta_x)\delta\bm{u}^\delta_x|_\Gamma. 
\end{cases}
\end{equation}
Although $h_{2,1}$ contains $\Phi(\varepsilon\eta, \varepsilon\delta\eta_x)\varepsilon\delta\eta_x u_y$ 
in addition to the above terms (see \eqref{NLnor}), 
by using the boundary condition $u_y=-\delta^2 v_x+(2+b_3)\eta+h_1$ on $\Gamma$, 
we can reduce the estimate of $h_{2,1}$ to that of $h_1$. Moreover, we note that $\delta^2{\rm W}h_{2,2}$ is of the form $\delta^2{\rm W}\Phi_0(\varepsilon^2\delta^2\eta_x^2)\eta_{xx}$.

\begin{lem}\label{est-h}
For any $\epsilon>0$ there exists a positive constant $C_\epsilon$ such that we have 
\begin{align}
& \delta^{-1}|(h_1,h_2)|_m^2+\delta|(h_{1x},h_{2x})|_m^2
 \lesssim \tilde{E}_2\tilde{F}_m+\tilde{F}_2\tilde{E}_m, \label{est-hx1} \\
& \delta^2|(h_{1x},h_{2x})|_{m+\frac12}^2
 \lesssim \tilde{E}_2\tilde{F}_m+\tilde{F}_2\tilde{E}_m, \label{est-hx2} \\
& \delta^{-1}|(h_1,h_2)|_{m-\frac12}^2 \lesssim \tilde{E}_2\tilde{E}_m, \label{est-hx3} \\
& \delta|h_2|_{m+\frac12}^2 \lesssim D_2D_m,\label{est-hx4} \\
& \delta|(\partial_x^k h_{1t}, \partial_x^k u_{t})_\Gamma|
 + \delta|(\partial_x^k h_{2t}, \delta\partial_x^k v_{t})_\Gamma|
 \le \epsilon\tilde{F}_m+C_{\epsilon}(\tilde{E}_2\tilde{F}_m+\tilde{F}_2\tilde{E}_m), \label{est-ht} \\
&
\begin{cases}\label{est-h3}
 \delta|h_{3}|_m^2+\delta^3|(h_{3x},h_{3t})|_m^2+\delta^5|h_{3xx}|_m^2
  \lesssim\tilde{E}_2\tilde{F}_m+\tilde{F}_2\tilde{E}_m, \\
 \delta^6{\rm W}|(\partial_x^k\eta_{xxxx}, \partial_x^k h_{3xx})_\Gamma|
  \le \epsilon\tilde{F}_m+C_\epsilon\tilde{E}_2\tilde{F}_m,\\
\end{cases} \\
&
\begin{cases}\label{est-b3I5}
 \delta|(b_3\eta)_x|_m^2+\delta^3|(b_3\eta)_{xx}|_m^2+\delta|(b_3\eta)_t|_m^2
  \lesssim \tilde{E}_2\tilde{F}_m+\tilde{F}_2\tilde{E}_m, \\
 |(\partial_x^k\eta, \partial_x^k(b_3\eta)_x)_\Gamma|
 +|(\partial_x^k\eta,\partial_x^kh_3)_\Gamma|+\delta^2{\rm W}|(\partial_x^k\eta_{xx},\partial_x^kh_3+\partial_x^k(b_3\eta)_x)_\Gamma|\\
\qquad  +\delta^4{\rm W}|(\partial_x^k\eta_{txx}, \partial_x^k h_{3t})_\Gamma|\le \epsilon\tilde{F}_m+C_\epsilon(\tilde{E}_2\tilde{F}_m+\tilde{F}_2\tilde{E}_m+\varepsilon\sqrt{\tilde{E}_2}\tilde{E}_m).
\end{cases}
\end{align}
Moreover, if $\varepsilon\lesssim\delta$, then we have 
\begin{align}\label{global}
&|(\partial_x^k\eta, \partial_x^k(b_3\eta)_x)_\Gamma|
 +|(\partial_x^k\eta,\partial_x^kh_3)_\Gamma|+\delta^2{\rm W}|(\partial_x^k\eta_{xx},\partial_x^kh_3+\partial_x^k(b_3\eta)_x)_\Gamma|\\
&\qquad+\delta^4{\rm W}|(\partial_x^k\eta_{txx}, \partial_x^k h_{3t})_\Gamma| \le \epsilon\tilde{F}_m + C_\epsilon(\tilde{E}_2\tilde{F}_m+\tilde{F}_2\tilde{E}_m).\notag
\end{align}
\end{lem}

\noindent
{\bf Remark 6.3.} \ 
Concerning the terms in the left-hand side of \eqref{est-b3I5}, in the case where $\varepsilon$ is not dominated by $\delta$, 
we cannot estimate these terms by using $\tilde{F}_m$ because the power of $\delta$ of these terms is not enough. 
These are the only terms which prevent from deriving a uniform estimate of the solution for all time.

\bigskip
\noindent
{\it Proof}. \ 
Since $\varepsilon$ is the nonlinear parameter, that is, $\varepsilon$ measures the nonlinearity, 
it is sufficient to show the estimates in the case $\varepsilon=1$ except the last estimate \eqref{global}. 
Therefore, we will assume that $\varepsilon=1$ in the following. 

As for \eqref{est-hx1}, it suffices to estimate 
$$
\begin{cases}
J_1:=\delta^{2i-1}|\Phi_0^1\partial_x^i\eta|_m^2 & \mbox{for}\quad i=1,2,3, \\
J_2:=\delta^{2i-1}|\Phi_0^1\partial_x^i\bm{u}^\delta|_m^2 & \mbox{for}\quad i=1,2,\\
J_3:=\delta^{2i+3}{\rm W}^2|\Phi_0^2\partial_x^i\eta_{xx}|_m^2& \mbox{for}\quad i=0,1,
\end{cases}
$$
where $\Phi_0^1=\Phi_0(\eta, \delta\eta_x, \delta^2\eta_{xx}, \bm{u}^\delta|_\Gamma, \delta^2 v_x|_\Gamma)$ and $\Phi_0^2=\Phi_0(\delta\eta_x, \delta^2\eta_{xx})$. 
Note that we included the term $\delta^2 v_x|_\Gamma$ in $\Phi_0^1$ for later use, although we can drop it. 
In the following we use the inequality 
\begin{equation}\label{product}
|fg|_s \lesssim |f|_{L^{\infty}}|g|_s+|g|_{L^{\infty}}|f|_s. 
\end{equation}
By \eqref{product}, \eqref{LI-E} in Lemma \ref{LI}, and \eqref{NLlem3-1} and \eqref{NLlem3-5} in Lemma \ref{NLlem3}, 
we obtain $J_1\lesssim \tilde{E}_2\tilde{F}_m+\tilde{F}_2\tilde{E}_m$. 
By \eqref{product}, the second inequality in \eqref{LI-u} in Lemma \ref{LI}, \eqref{NLlem3-1} and \eqref{NLlem3-5} 
in Lemma \ref{NLlem3}, and the second inequality in \eqref{gammaku} in Lemma \ref{additonal lemma}, 
we obtain $J_2\lesssim \tilde{E}_2\tilde{F}_m+\tilde{F}_2\tilde{E}_m$.
By \eqref{product}, \eqref{WL} in Lemma \ref{LI}, and \eqref{NLlemW-1} and \eqref{NLlemW-5} in Lemma \ref{NLlemW}, 
we obtain $J_3\lesssim \tilde{E}_2\tilde{F}_m+\tilde{F}_2\tilde{E}_m$. 
Thus \eqref{est-hx1} holds. 

As for \eqref{est-hx2}, it suffices to estimate 
$$
\begin{cases}
J_4:=\delta^{2i}|\Phi^1_0\partial_x^i\eta|^2_{m+\frac12} & \mbox{for}\quad i=1,2,3, \\
J_5:=\delta^{2i}|\Phi^1_0\partial_x^i\bm{u}^\delta|^2_{m+\frac12} & \mbox{for}\quad i=1,2,\\
J_6:=\delta^{2i+4}{\rm W}^2|\Phi_0^2\partial_x^i\eta_{xx}|_{m+\frac12}^2& \mbox{for}\quad i=0,1.
\end{cases}
$$
By \eqref{product}, \eqref{LI-E} in Lemma \ref{LI}, the second inequality in \eqref{interpo} 
in Lemma \ref{additonal lemma}, and \eqref{NLlem3-1} and \eqref{NLlem3-3} in Lemma \ref{NLlem3}, 
we obtain $J_4 \lesssim \tilde{E}_2\tilde{F}_m+\tilde{F}_2\tilde{E}_m$. 
By \eqref{product}, the second inequality in \eqref{LI-u} in Lemma \ref{LI}, 
the second and third inequalities in \eqref{fracu} in Lemma \ref{additonal lemma}, 
and \eqref{NLlem3-1} and \eqref{NLlem3-3} in Lemma \ref{NLlem3}
we obtain $J_5 \lesssim \tilde{E}_2\tilde{F}_m+\tilde{F}_2\tilde{E}_m$. 
By \eqref{product}, \eqref{WL} in Lemma \ref{LI}, the second inequality in \eqref{interpoW} 
in Lemma \ref{additonal lemma}, and \eqref{NLlemW-1} and \eqref{NLlemW-3} in Lemma \ref{NLlemW}, we obtain $J_6\lesssim \tilde{E}_2\tilde{F}_m+\tilde{F}_2\tilde{E}_m$. 
Thus \eqref{est-hx2} holds.

As for \eqref{est-hx3}, it suffices to estimate 
$$
\begin{cases}
J_7:=\delta^{2i+1}|\Phi_0^1\partial_x^i\eta_x|_{m-\frac12}^2 & \mbox{for}\quad i=0,1, \\
J_8:=\delta|\Phi_0^1\bm{u}^\delta_x|_{m-\frac12}^2,\\
J_9:=\delta^{3}{\rm W}^2|\Phi_0^2\eta_{xx}|_{m-\frac12}^2.
\end{cases}
$$
By \eqref{product}, \eqref{LI-E} in Lemma \ref{LI}, the first inequality in \eqref{interpo} 
in Lemma \ref{additonal lemma}, and \eqref{NLlem3-1} and \eqref{NLlem3-5} in Lemma \ref{NLlem3}, 
we obtain $J_7 \lesssim \tilde{E}_2\tilde{E}_m$. 
By \eqref{product}, \eqref{LI-u} in Lemma \ref{LI}, the first inequality in \eqref{fracu} in 
Lemma \ref{additonal lemma}, and \eqref{NLlem3-1} and \eqref{NLlem3-5} in Lemma \ref{NLlem3}, 
we obtain $J_8 \lesssim\tilde{E}_2\tilde{E}_m$. 
By \eqref{product}, \eqref{WL} in Lemma \ref{LI}, the first inequality in \eqref{interpoW} 
in Lemma \ref{additonal lemma}, and \eqref{NLlemW-1} and \eqref{NLlemW-5} in Lemma \ref{NLlemW}, 
we obtain $J_9 \lesssim \tilde{E}_2\tilde{E}_m$. 
Thus \eqref{est-hx3} holds. 

As for \eqref{est-hx4}, it suffices to estimate 
$$
\begin{cases}
J_{10}:=\delta^{2i+1}|\Phi_0^1\partial_x^i\eta|_{m+\frac12}^2 & \mbox{for}\quad i=1,2, \\
J_{11}:=\delta^{3}|\Phi_0^1\bm{u}^\delta_x|_{m+\frac12}^2,\\
J_{12}:=\delta^{5}{\rm W}^2|\Phi_0^2\eta_{xx}|_{m+\frac12}^2.
\end{cases}
$$
By \eqref{product}, \eqref{LI-E} in Lemma \ref{LI}, the first inequality in \eqref{interpo} 
in Lemma \ref{additonal lemma}, and \eqref{NLlem3-1} and \eqref{NLlem3-3} in Lemma \ref{NLlem3}, 
we obtain $J_{10} \lesssim D_2D_m$. 
By \eqref{product}, \eqref{LI-u} in Lemma \ref{LI}, the first inequality in \eqref{fracu} in 
Lemma \ref{additonal lemma}, and \eqref{NLlem3-1} and \eqref{NLlem3-3} in Lemma \ref{NLlem3}, 
we obtain $J_{11} \lesssim D_2D_m$. 
By \eqref{product}, \eqref{WL} in Lemma \ref{LI}, the first inequality in \eqref{interpoW} 
in Lemma \ref{additonal lemma}, and \eqref{NLlemW-1} and \eqref{NLlemW-3} in Lemma \ref{NLlemW}, 
we obtain $J_{12} \lesssim D_2D_m$. 
Thus \eqref{est-hx4} holds. 

We proceed to estimate \eqref{est-ht}. 
By the third equation in \eqref{BC1}, we can reduce the estimates of the terms which contain $\eta_t$ except the terms which accompany ${\rm W}$
to those of $J_1$, $J_2$, and $J_4$. 
Thus it suffices to estimate 
$$
\begin{cases}
J_{13}:=\delta^9{\rm W}^2|\Phi^3\eta_x\eta_{xx}\eta_{tx}|_m^2,\\
J_{14}:=\delta^9{\rm W}^2|\Phi^3\eta_x^2\eta_{txx}|_m^2,\\
J_{15}:=\delta|\Phi^1_0\bm{u}^\delta_t|_m^2,\\
J_{16}:=\delta^2|(\partial_x^k(\Phi^4_0\bm{u}^\delta_{tx}), \partial_x^k\bm{u}^\delta_t)_\Gamma|,
\end{cases}
$$
where $\Phi^3=\Phi(\delta\eta_x)$, $\Phi_0^4=\Phi_0(\eta, \delta\eta_x)$ and we used $h_{2,2}=\Phi_0(\delta^2\eta_x^2)\eta_{xx}=\Phi(\delta\eta_x)\delta^2\eta_x^2\eta_{xx}$.
Taking into account that $\delta^9{\rm W}^2\lesssim\delta^5$, by the third equation in \eqref{BC1}, we can reduce the estimates of $J_{13}$ and $J_{14}$ 
to those of $J_1$, $J_2$, and $J_4$. 
By \eqref{product}, the second inequality in \eqref{LI-u} in Lemma \ref{LI}, 
\eqref{NLlem3-1} and \eqref{NLlem3-5} in Lemma \ref{NLlem3}, and 
$\delta|\bm{u}^{\delta}_t|_m^2 \lesssim \delta\|(1+|D_x|)^m\bm{u}^{\delta}_{ty}\|^2 \lesssim \tilde{F}_m$, we obtain 
$J_{15} \lesssim \tilde{E}_2\tilde{F}_m+\tilde{F}_2\tilde{E}_m$. 
By Lemma \ref{trace}, we see that 
\begin{align*}
J_{16}
&= \delta^2|(\partial_x^k\{(\Phi^4_0\bm{u}^\delta_t)_x-\Phi^4_{0x}\bm{u}^\delta_t\}, 
 \partial_x^k\bm{u}^\delta_t)_\Gamma| \\
&\le \delta^2||D_x|^{\frac12}\partial_x^k(\Phi^4_0\bm{u}^\delta_t)|_0||D_x|^{\frac12}\partial_x^k\bm{u}^\delta_t|_0
 +\delta^2|(\partial_x^k(\Phi^4_{0x}\bm{u}^\delta_t), \partial_x^k\bm{u}^\delta_t)_\Gamma| \\
&\le \epsilon(\delta\|\partial_x^k\bm{u}^\delta_t\|^2+\delta^3\|\partial_x^k\bm{u}^\delta_{tx}\|^2
  +\delta\|\partial_x^k\bm{u}^\delta_{ty}\|^2)
 +C_\epsilon\bigl(\delta^2|\Phi^4_0\bm{u}^\delta_t|^2_{m+\frac12}
  +\delta^{3}|\Phi^4_{0x}\bm{u}^\delta_t|_m^2\bigr).
\end{align*}
Here, we can reduce the estimate of $\delta^{3}|\Phi^4_{0x}\bm{u}^\delta_t|_m^2$ to that of $J_8$. 
By \eqref{product}, the second inequality in \eqref{LI-u} in Lemma \ref{LI}, 
\eqref{gammaut} in Lemma \ref{additonal lemma}, and \eqref{NLlem3-1} and \eqref{NLlem3-3} in Lemma \ref{NLlem3},  
we obtain $\delta^2|\Phi^4_0\bm{u}^\delta_t|^2_{m+\frac12} \lesssim \tilde{E}_2\tilde{F}_m+\tilde{F}_2\tilde{E}_m$. 
We thereby deduce $J_{16} \le \epsilon \tilde{F}_m+C_{\epsilon}(\tilde{E}_2\tilde{F}_m+\tilde{F}_2\tilde{E}_m)$. 
Thus \eqref{est-ht} holds.

As for \eqref{est-h3}, since $h_3=\eta^2\eta_x$ is contained in the first term in \eqref{expl-h}, 
we have already checked that the first inequality holds. 
As for the second inequality, we have 
$$
\delta^6{\rm W}|(\partial_x^k\eta_{xxxx}, \partial_x^k h_{3xx})_\Gamma|
 \le \epsilon\delta^6{\rm W}||D_x|^{k+\frac72}\eta|^2_0+C_{\epsilon}\delta^6{\rm W}||D_x|^{k+\frac52}h_3|_0^2.
$$
Here, \eqref{estimate of h3} leads to $\delta^6{\rm W}||D_x|^{k+\frac52}h_3|_0^2 \lesssim \tilde{E}_2\tilde{F}_m$. 
Therefore, we get the second inequality. 

As for \eqref{est-b3I5}, taking into account that we can write $b_3$ as $\Phi_0^4$ (see \eqref{def-b3}), 
we obtain the first inequality in the same reason as the last estimate. 
Concerning the term $|(\partial_x^k\eta, \partial_x^k(b_3\eta)_x)_\Gamma|$ in the second inequality, there exist rational functions $b_{3,1}$ and $b_{3,2}$ such that 
$b_3\eta=b_{3,1}(\eta)+b_{3,2}(\eta,\delta\eta_x)\delta\eta_x$ and $b_{3,2}(\bm{0})=0$. 
Since the term $b_{3,2}(\eta,\delta\eta_x)\delta\eta_x$ can be treated in the same way as before, 
it suffices to estimate 
$$
J_{17}:=|(\partial_x^k\eta,\partial_x^{k+1}b_{3,1}(\eta))_\Gamma|. 
$$
Here we can assume that $k\ge1$ because we have $(\eta,b_{3,1}(\eta)_x)_\Gamma=0$ in the case $k=0$. 
We see that 
$J_{17} \le |(\partial_x^k\eta,b_{3,1}'(\eta)\partial_x^{k+1}\eta)_\Gamma|
 +|(\partial_x^k\eta,[\partial_x^k,b_{3,1}'(\eta)]\eta_x)_\Gamma)|$, 
where by integration by parts we have 
$|(\partial_x^k\eta,b_{3,1}'(\eta)\partial_x^{k+1}\eta)_\Gamma|
= \frac12|(\partial_x^k\eta,b_{3,1}''(\eta)\eta_x\partial_x^k\eta)_\Gamma|
\lesssim \sqrt{\tilde{E}_2}|\eta_x|_{m-1}^2$. 
In view of 
$$
|[\partial_x^k,b_{3,1}'(\eta)]\eta_x|_0
 \le C(|\eta|_{L^{\infty}})(1+|\eta_x|_{L^{\infty}})^{k-1}|\eta_x|_{L^{\infty}}|\partial_x^{k-1}\eta_x|_0, 
$$
we also have 
$|(\partial_x^k\eta,[\partial_x^k,b_{3,1}'(\eta)]\eta_x)_\Gamma| \lesssim \sqrt{\tilde{E}_2}|\eta_x|_{m-1}^2$. 
Therefore, $J_{17} \lesssim \sqrt{\tilde{E}_2}\min\{\tilde{E}_m,\delta^{-1}\tilde{F}_m\}$, so that 
we obtain 
$$
|(\partial_x^k\eta,\partial_x^k(b_3\eta)_x)_\Gamma|
\le \epsilon\tilde{F}_m+C_{\epsilon}(\tilde{E}_2\tilde{F}_m+\tilde{F}_2\tilde{E}_m)
 +C\sqrt{\tilde{E}_2}\min\{\tilde{E}_m,\delta^{-1}\tilde{F}_m\}.
$$
As for the the term $\delta^4{\rm W}|(\partial_x^k\eta_{txx}, \partial_x^k h_{3t})_\Gamma|$, integration by parts in $x$ leads to
\begin{align*}
\delta^4{\rm W}|(\partial_x^k\eta_{txx}, \partial_x^kh_{3t})_\Gamma|
 & \le \delta^4{\rm W}|(\partial_x^k\eta_{txx}, \partial_x^k(\eta^2\eta_{tx}))_\Gamma|+\delta^4{\rm W}|(\partial_x^k\eta_{txx}, \partial_x^k(2\eta\eta_t\eta_x))_\Gamma|\\
 & \le \delta^4{\rm W}|(\partial_x^k\eta_{tx}, \eta\eta_x\partial_x^k\eta_{tx})_\Gamma|+\delta^4{\rm W}|(\partial_x^k\eta_{tx}, ([\partial_x^k, \eta^2]\eta_{tx})_x)_\Gamma|\\
   &\quad+\delta^4{\rm W}|(\partial_x^k\eta_{tx}, \partial_x^k(2\eta\eta_t\eta_x)_x)_\Gamma|\\
 &=:J_{18}+J_{19}+J_{20}.
\end{align*}
Here, it follows from the third equation in \eqref{BC1} that $\delta^4{\rm W}|\partial_x^k\eta_{tx}|_0^2\lesssim\delta^2\bigl(|\partial_x^k\eta_{xx}|_0^2+|\partial_x^kv_{x}|_0^2+|\partial_x^kh_{3x}|_0^2\bigr)\lesssim\delta^{-1}\tilde{F}_m$ and $\delta^4{\rm W}|\partial_x^k\eta_{tx}|_0^2\lesssim  \tilde{E}_m$ so that we have
\begin{equation}\label{We-eta-time}
\delta^4{\rm W}|\partial_x^k\eta_{tx}|_0^2\lesssim\min\{\tilde{E}_m, \delta^{-1}\tilde{F}_m\}.
\end{equation}
By \eqref{LI-E} in Lemma \ref{LI} and \eqref{We-eta-time}, we have $J_{18}\lesssim\tilde{E}_2\min\{\tilde{E}_m, \delta^{-1}\tilde{F}_m\}$. By the estimate
\[
|([\partial_x^k, \eta^2]\eta_{tx})_x|_0\lesssim|\eta|_{L^\infty}|\eta_x|_{L^\infty}|\partial_x^k\eta_{tx}|+|\eta|_{L^\infty}|\eta_{tx}|_{L^\infty}|\partial_x^k\eta_{x}|+|\eta_x|_{L^\infty}|\eta_{tx}|_{L^\infty}|\partial_x^k\eta|,
\]
\eqref{LI-E}, \eqref{WL}, and \eqref{LI-v} in Lemma \ref{LI}, and \eqref{We-eta-time}, we easily obtain $J_{19}\le\epsilon\tilde{F}_m+C_{\epsilon}\tilde{F}_2\tilde{E}_m+C\tilde{E}_2\min\{\tilde{E}_m, \delta^{-1}\tilde{F}_m\}$. By \eqref{product} and \eqref{LI-E}, \eqref{WL}, and \eqref{LI-v} in Lemma \ref{LI}, and \eqref{We-eta-time}, we have $J_{20}\le\epsilon\tilde{F}_m+C_{\epsilon}\tilde{E}_2\tilde{F}_m+C\tilde{E}_2\min\{\tilde{E}_m, \delta^{-1}\tilde{F}_m\}$. Therefore, we get the third inequality.
Thus far, we have assumed that $\varepsilon=1$. 
Now, for general $\varepsilon\in(0,1]$ it follows easily from the above estimate that 
$$
|(\partial_x^k\eta,\partial_x^k(b_3\eta)_x)_\Gamma|
\le \epsilon\tilde{F}_m+\varepsilon^2C_{\epsilon}(\tilde{E}_2\tilde{F}_m+\tilde{F}_2\tilde{E}_m)
 +C\sqrt{\tilde{E}_2}\min\{\varepsilon\tilde{E}_m,\varepsilon\delta^{-1}\tilde{F}_m\}.
$$
The term $|(\partial_x^k\eta,\partial_x^kh_3)_\Gamma|$ is of the form $J_{17}$, so that it 
also satisfies the above estimate. Moreover, by taking into account that $\delta\sqrt{{\rm W}}|\eta_x|_{L^\infty}\lesssim\sqrt{\tilde{E}_2}$ and $\delta^2{\rm W}|\eta_x|_{m-1}^2\lesssim|\eta_x|_{m-1}^2\lesssim\min\{\tilde{E}_m,\delta^{-1}\tilde{F}_m\}$, the term $\delta^2{\rm W}|(\partial_x^k\eta_{xx},\partial_x^kh_3+\partial_x^k(b_3\eta)_x)_\Gamma|$ also satisfies the above estimate. Similarly, we obtain
\[
\delta^4{\rm W}|(\partial_x^k\eta_{txx}, \partial_x^kh_{3t})_\Gamma|\le \epsilon\tilde{F}_m+\varepsilon^2C_{\epsilon}(\tilde{E}_2\tilde{F}_m+\tilde{F}_2\tilde{E}_m)
 +C\tilde{E}_2\min\{\varepsilon\tilde{E}_m,\varepsilon\delta^{-1}\tilde{F}_m\}.
\]
Therefore, the second inequalities in \eqref{est-b3I5} and \eqref{global} hold. 
The proof is complete. 
$\quad\square$

\bigskip
Next, we will estimate $f_1$, $f_2$, $\bm{F}_1$, $\bm{F}_2$, and $\bm{G}_k$. 
By the explicit form of $\bm{f}$ (see \eqref{NLdomain}), we see that this is consist of terms in the form 
$$
\begin{cases}
\Phi_0(\tilde{\eta}, D_\delta\tilde{\eta}, \bm{u}^\delta;y)D_\delta^3\tilde{\eta}, \\
\Phi_0(\tilde{\eta}, D_\delta\tilde{\eta};y)\delta^i \partial_x^{i}\partial_y^{j}u
  & \mbox{for}\quad (i, j)=(2, 0), (1,1), \\
\Phi(\tilde{\eta}, D_\delta\tilde{\eta}, \bm{u}^\delta,y)w_\lambda u_y  & \mbox{for}\quad 1\le \lambda\le3, \\
\Phi(\tilde{\eta}, D_\delta\tilde{\eta}, D_\delta^2\tilde{\eta}, \bm{u}^\delta,y)w_\lambda\bm{u}^\delta
  & \mbox{for}\quad 1\le \lambda\le4, \\
\Phi_0(\tilde{\eta}, D_\delta\tilde{\eta}, D_\delta^2\tilde{\eta}, \delta\tilde{\eta}_t, \bm{u}^\delta;y)
 \delta\bm{u}^\delta_x,
\end{cases}
$$
where $w_\lambda$ is defined by \eqref{def-w}. 
Thus by the explicit forms of $f_1$ and $f_2$ (see \eqref{NLdomainf2} and \eqref{NLdomainf3}), 
we see that these contain the above terms, $\Phi_0(\tilde{\eta}, D_\delta\tilde{\eta};y)\nabla_\delta p$, 
and $\Phi_0(\tilde{\eta},D_\delta\tilde{\eta};y)\delta u_t$ (see also \eqref{NLp}). 
In addition to these terms, $\bm{F}_1$ contains also 
$\Phi_0(\tilde{\eta}, D_\delta\tilde{\eta};y)u_{yy}$ (see \eqref{NLdomainF}).

\begin{lem}\label{est-f}
For any $\epsilon>0$ there exists a positive constant $C_\epsilon$ such that the following estimates hold. 
\begin{align}
& \delta^{-1}\|\partial_x^k f_1\|^2+\delta^{-1}\|\partial_x^k f_2\|^2
   +\delta\|\partial_x^k f_{1x}\|^2 \label{est-fx1} \\
&\quad
  +\delta|(\partial_x^k\bm{F}_{1x}, \partial_x^k\bm{u}^\delta_x)_\Omega|
   +\delta^3|(\partial_x^k\bm{F}_{1xx}, \partial_x^k\bm{u}^\delta_{xx})_\Omega|
 \le \epsilon\tilde{F}_m+C_\epsilon(\tilde{E}_2\tilde{F}_m+\tilde{F}_2\tilde{E}_m), \notag \\
& \delta^{-2}\|\partial_x^{l-1}\bm{f}\|^2 \lesssim \tilde{E}_2\tilde{E}_m, \label{est-fx2} \\
& |(\partial_x^l\{A_4\nabla_\delta p+(b_2u_{yy}, 0)^{\rm T}\}, \partial_x^l\bm{u}^\delta)_\Omega|
  \le (\epsilon+C_\epsilon\tilde{E}_2)\tilde{E}_m, \label{est-fx4} \\
& \|\partial_x^k\bm{f}\|^2 \lesssim D_2D_m, \label{est-fx3} \\
& \delta|(\partial_x^k\bm{F}_2, \partial_x^k\bm{u}^\delta_t)_\Omega|
  \le \epsilon\tilde{F}_m+C_\epsilon(\tilde{E}_2\tilde{F}_m+\tilde{F}_2\tilde{E}_m), \label{est-ft} \\
& \delta|(\bm{G}_k, \partial_x^k\bm{u}^\delta_t)_\Omega|
  \le \epsilon\tilde{F}_m+C_\epsilon(\tilde{E}_2\tilde{F}_m+\tilde{F}_2\tilde{E}_m). \label{est-commutator}
\end{align}
\end{lem}

\noindent
{\it Proof}. \ 
As for \eqref{est-fx1}, the definition of $\bm{F}_1$ and integration by parts in $x$ imply 
\begin{align*}
\delta^3|(\partial_x^k\bm{F}_{1xx}, \partial_x^k\bm{u}^\delta_{xx})_\Omega|
&\le \epsilon\delta^5\|\partial_x^k\bm{u}^\delta_{xxx}\|^2
 +C_\epsilon\delta\|\partial_x^k(\bm{f}-\Phi_0(\tilde{\eta}, D_\delta\tilde{\eta})\nabla_\delta p)_x\|^2 \\
&\quad
 +\delta^3|(\partial_x^k(\Phi_0(\tilde{\eta}, D_\delta\tilde{\eta})u_{yy})_{xx},
  \partial_x^k\bm{u}^\delta_{xx})_\Omega|.
\end{align*}
Taking this into account, it suffices to estimate 
$$
\begin{cases}
K_1:=\delta^{-1}\|\partial_x^k(\Phi_0^5D_\delta^i\tilde{\eta})\|^2 & \mbox{for}\quad 1\le i\le4, \\
K_2:=\delta^{-1}\|\partial_x^k(\Phi_0^5\delta\tilde{\eta}_t)\|^2, \\
K_3:=\delta^{2i-1}\|\partial_x^k(\Phi_0^5\partial_x^iu_y)\|^2 & \mbox{for}\quad i=1, 2, \\
K_4:=\delta^{-1}\|\partial_x^k(\Phi^5w_\lambda\partial_y^j u)\|^2 & \mbox{for}\quad 1\le \lambda\le7, \ j=0, 1, \\
K_5:=\delta^{2i-1}\|\partial_x^k(\Phi_0^6 \partial_x^i\bm{u}^\delta)\|^2 & \mbox{for}\quad  i=1,2,3 \\
K_6:=\delta^{2i-1}|(\partial_x^{k+i}(\Phi_0^7 u_{yy}), \partial_x^{k+i}\bm{u}^\delta)_\Omega|
 & \mbox{for}\quad i=1, 2, \\
K_7:=\delta^{2i-1}\|\partial_x^{k+i}(\Phi_0^7\nabla_\delta p)\|^2 & \mbox{for}\quad i=0, 1,
\end{cases}
$$
where 
\begin{align*}
& \Phi^5=\Phi(\tilde{\eta}, D_\delta\tilde{\eta}, D_\delta^2\tilde{\eta}, 
 \delta\tilde{\eta}_t, \delta D_\delta\tilde{\eta}_t, \bm{u}^\delta;y), \\
& \Phi^6=\Phi(\tilde{\eta}, D_\delta\tilde{\eta}, D_\delta^2\tilde{\eta}, D_\delta^3\tilde{\eta}, 
 \delta\tilde{\eta}_t, \delta D_\delta\tilde{\eta}_t, \delta D_\delta^2\tilde{\eta}_t, 
 \bm{u}^\delta, \delta\bm{u}^\delta_x;y), \\
& \Phi^7=\Phi(\tilde{\eta}, D_\delta\tilde{\eta};y). 
\end{align*}
In the following we will use the well-known inequality 
\begin{equation}\label{interpolation}
\|\partial_x^k(fg)\| \lesssim \|f\|_{L^{\infty}}\|\partial_x^kg\|+\|g\|_{L^{\infty}}\|\partial_x^kf\|.
\end{equation}
By this, \eqref{LI-E} in Lemma \ref{LI}, \eqref{dxkE} and \eqref{est-highE} in Lemma \ref{L2estimates}, 
and \eqref{NLlem2-1} and \eqref{NLlem2-3} in Lemma \ref{NLlem2}, we obtain 
$K_1 \lesssim \tilde{E}_2\tilde{F}_m+\tilde{F}_2\tilde{E}_m$. 
Similarly, we get $K_2 \lesssim \tilde{E}_2\tilde{F}_m+\tilde{F}_2\tilde{E}_m$.
By Lemma \ref{NLlem}, we have 
$$
K_3 \lesssim \|\Phi_0^5\|^2_{L^{\infty}}\delta^{2i-1}\|\partial_x^{k+i}u_y\|^2
 +(\|\partial_x^k \Phi_0^5\|^2+\|\partial_x^k \Phi_{0y}^5\|^2)
  \delta^{2i-1}(\|\partial_x^iu_y\|^2+\|\partial_x^iu_{xy}\|^2),
$$
which together with \eqref{NLlem2-1} and \eqref{NLlem2-5} in Lemma \ref{NLlem2} 
gives $K_3\lesssim \tilde{E}_2\tilde{F}_m+\tilde{E}_m\tilde{F}_2$. 
By Lemma \ref{NLlem}, we have 
\begin{align*}
K_4
&\lesssim \|\Phi^5\|_{L^{\infty}}^2\delta^{-1}\|w_\lambda\|_{L^{\infty}}^2\|\partial_x^k \partial_y^ju\|^2 \\
&\quad
 + \delta^{-1}\|w_\lambda\|_{L^{\infty}}^2(\|\partial_x^k \Phi^5\|^2+\|\partial_x^k \Phi_y^5\|^2)
  (\|\partial_y^ju\|^2+\|\partial_y^ju_x\|^2)\\
&\quad
 +\|\Phi^5\|_{L^{\infty}}^2\delta^{-1}(\|\partial_x^k w_\lambda\|^2+\|\partial_x^k w_{\lambda y}\|^2)
  (\|\partial_y^ju\|^2+\|\partial_y^ju_x\|^2),
\end{align*}
which together with \eqref{NLlem2-1} and \eqref{NLlem2-5} in Lemma \ref{NLlem2} 
and \eqref{w-1} and \eqref{w-2} in Lemma \ref{NLlemw} 
gives $K_4\lesssim \tilde{E}_2\tilde{F}_m+\tilde{F}_2\tilde{E}_m$. 
As for $K_5$, it suffices to consider the case of $k\ge1$ since we can easily treat the case of $k=0$. 
By Lemma \ref{NLlem}, we have 
$$
K_5 \lesssim \|\Phi_0^6\|^2_{L^{\infty}}\delta^{2i-1}\|\partial_x^{k+i}\bm{u}^\delta\|^2
 +\delta(\|\partial_x^k \Phi_0^6\|^2+\|\partial_x^k \Phi_{0y}^6\|^2)\delta^{2(i-1)}
  (\|\partial_x^i\bm{u}^\delta\|^2+\|\partial_x^i\bm{u}^\delta_x\|^2),
$$
which together with \eqref{NLlem2-1} and \eqref{NLlem2-4} in Lemma \ref{NLlem2} gives 
$K_5 \lesssim \tilde{E}_2\tilde{F}_m$. 
As for $K_6$, we will consider the case $i=2$ only, 
because the case where $i=1$ can be treated in a similar but easier way. 
Using integration by parts in $x$ and $y$ and Lemma \ref{trace}, we have 
\begin{align*}
K_6
&= \delta^3|(\partial_x^k\{(\Phi_0^7u_y)_{xxy}-(\Phi_{0y}^7u_y)_{xx}\}, \partial_x^k\bm{u}^\delta_{xx})_\Omega| \\
&\le \delta^3|(\partial_x^k(\Phi_0^7u_y)_{xx}, \partial_x^k\bm{u}^\delta_{xxy})_\Omega|
  +\delta^3||D_x|^{\frac12}\partial_x^k(\Phi_0^7u_y)_{x}|_0||D_x|^{\frac12}\partial_x^k \bm{u}^\delta_{xx}|_0 \\
&\quad
 +\delta^3|(\partial_x^k(\Phi_{0y}^7u_y)_{x}, \partial_x^k\bm{u}^\delta_{xxx})_\Omega| \\
&\le \epsilon\tilde{F}_m
  +C_{\epsilon}\bigl(\delta^{3}\|\partial_x^k(\Phi_0^7u_y)_{xx}\|^2
   +\delta\|\partial_x^k(\Phi_{0y}^7u_y)_x\|^2+\delta^2||D_x|^{\frac12}\partial_x^k(\Phi_0^7u_y)_x|^2_0\bigr).
\end{align*}
Here we can reduce the estimate of 
$\delta^{3}\|\partial_x^k(\Phi_0^7u_y)_{xx}\|^2+\delta\|\partial_x^k(\Phi_{0y}^7u_y)_x\|^2$ 
to those of $K_3$ and $K_4$. 
Furthermore, using the first equation in \eqref{BC1} to eliminate $u_y|_\Gamma$, 
we can reduce the estimate of $\delta^2||D_x|^{\frac12}\partial_x^k(\Phi_0^7u_y)_x|^2_0$ to those of $J_4$ and $J_5$. 
Thus combining these estimates, we obtain 
$K_6 \le \epsilon\tilde{F}_m+C_\epsilon(\tilde{E}_2\tilde{F}_m+\tilde{F}_2\tilde{E}_m)$. 
By Lemma \ref{NLlem}, we have 
$$
K_7 \lesssim \|\Phi_0^7\|_{L^{\infty}}^2\delta^{2i-1}\|\nabla_\delta\partial_x^{k+i}p\|^2
 +\delta^{2i}(\|\partial_x^{k+i}\Phi_0^7\|^2+\|\partial_x^{k+i}\Phi_{0y}^7\|^2)
 \delta^{-1}(\|\nabla_\delta p\|^2+\|\nabla_\delta p_x\|^2),
$$
which together with \eqref{NLlem2-1} and \eqref{NLlem2-5} in Lemma \ref{NLlem2} gives 
$K_7 \lesssim \tilde{E}_2\tilde{F}_m+\tilde{E}_m\tilde{F}_2$. 
Thus \eqref{est-fx1} holds. 

As for \eqref{est-fx2}, it suffices to estimate 
$$
\begin{cases}
K_8:=\delta^{-2}\|\partial_x^{l-1}(\Phi_0^5D_\delta^3\tilde{\eta})\|^2, \\
K_9:=\|\partial_x^{l-1}(\Phi_0^5 u_{xy})\|^2, \\
K_{10}:=\delta^{-2}\|\partial_x^{l-1}(\Phi^5 w_\lambda\partial_y^j\bm{u}^\delta)\|^2
 & \mbox{for}\quad 1\le\lambda\le4,\  j=0,1, \\
K_{11}:=\delta^{2i}\|\partial_x^{l-1}(\Phi_0^5\partial_x^i\bm{u}^\delta_x)\|^2
 & \mbox{for}\quad i=0,1.
\end{cases}
$$
By \eqref{extension2} in Lemma \ref{extension}, we have 
$\delta^{-2}\|\partial_x^{l-1}D_\delta^3\tilde{\eta}\|^2 \lesssim \delta^4|\partial_x^l\eta_{xx}|_0^2$. 
Therefore, by \eqref{interpolation}, \eqref{LI2} in Lemma \ref{LI}, 
and \eqref{NLlem2-1} and \eqref{NLlem2-3} in Lemma \ref{NLlem2}, we obtain 
$K_8 \lesssim \tilde{E}_2\tilde{E}_m$. 
By Lemma \ref{NLlem}, we have 
$$
K_9 \lesssim \|\Phi_0^5\|_{L^{\infty}}^2\|\partial_x^l\bm{u}^\delta_y\|^2
 +(\|\partial_x^{l-1}\Phi_0^5\|^2+\|\partial_x^{l-1}\Phi_{0y}^5\|^2)(\|u_{xy}\|^2+\|u_{xxy}\|^2),
$$
which together with \eqref{NLlem2-1} and \eqref{NLlem2-5} in Lemma \ref{NLlem2} gives 
$K_9 \lesssim \tilde{E}_2\tilde{E}_m$. 
By Lemma \ref{NLlem}, we have 
\begin{align*}
K_{10}
&\lesssim
 \|\Phi^5\|_{L^{\infty}}^2\delta^{-2}\|w_\lambda\|_{L^{\infty}}^2\|\partial_x^{l-1}\partial_y^j\bm{u}^\delta\|^2 \\
&\quad
 +\delta^{-2}\|w_\lambda\|_{L^{\infty}}^2(\|\partial_x^{l-1}\Phi^5\|^2+\|\partial_x^{l-1}\Phi_y^5\|^2)
  (\|\partial_y^j\bm{u}^\delta\|^2+\|\partial_y^j\bm{u}^\delta_x\|^2) \\
&\quad
 +\|\Phi^5\|_{L^{\infty}}^2\delta^{-2}(\|\partial_x^{l-1} w_\lambda\|^2+\|\partial_x^{l-1} w_{\lambda y}\|^2)
  (\|\partial_y^j\bm{u}^\delta\|^2+\|\partial_y^j\bm{u}^\delta_x\|^2),
\end{align*}
which together with \eqref{NLlem2-1} and \eqref{NLlem2-5} in Lemma \ref{NLlem2} 
and \eqref{w-1} and \eqref{w-3} in Lemma \ref{NLlemw} gives 
$K_{10} \lesssim \tilde{E}_2\tilde{E}_m$. 
By Lemma \ref{NLlem}, we have 
$$
K_{11} \lesssim \|\Phi_0^5\|_{L^{\infty}}^2\delta^{2i}\|\partial_x^{l+i}\bm{u}^\delta\|^2
 +(\|\partial_x^{l-1}\Phi_0^5\|^2+\|\partial_x^{l-1}\Phi_{0y}^5\|^2)
  \delta^{2i}(\|\partial_x^i\bm{u}^\delta_x\|^2+\|\partial_x^i\bm{u}^\delta_{xx}\|),
$$
which together with \eqref{NLlem2-1} and \eqref{NLlem2-5} in Lemma \ref{NLlem2} gives 
$K_{11} \lesssim \tilde{E}_2\tilde{E}_m$. 
Thus \eqref{est-fx2} holds. 

We proceed to estimate \eqref{est-fx4}. 
With the aid of \eqref{A4p}, we can express $A_4\nabla_\delta p$ in terms of the product of $\Phi_0^7$ 
and derivatives of $\bm{u}^\delta$ in addition to $\Phi_0^7\bm{f}$. 
Taking this into account and using \eqref{est-fx2}, it suffices to estimate 
$$
\begin{cases}
K_{12}:=\delta^2\|\partial_x^l(\Phi_0^7\bm{u}^\delta_t)\|^2, \\
K_{13}:=|(\partial_x^l(\Phi_0^7 \bm{u}^\delta_{yy}), \partial_x^l\bm{u}^\delta)_\Omega|.
\end{cases}
$$
By Lemma \ref{NLlem}, we have 
$$
K_{12} \lesssim \|\Phi_0^7\|_{L^{\infty}}^2\delta^2\|\partial_x^l\bm{u}^\delta_t\|^2
 +(\|\partial_x^l\Phi_0^7\|^2+\|\partial_x^l\Phi_{0y}^7\|^2)
  \delta^2(\|\bm{u}^\delta_t\|^2+\|\bm{u}^\delta_{tx}\|^2),
$$
which together with \eqref{NLlem2-1} and \eqref{NLlem2-5} in Lemma \ref{NLlem2} gives 
$K_{12} \lesssim \tilde{E}_2\tilde{E}_m$. 
Integration by parts in $y$ implies 
\begin{align*}
K_{13}
&= |(\partial_x^l\{(\Phi_0^7\bm{u}^\delta_y)_y-\Phi_{0y}^7\bm{u}^\delta_y\} ,\partial_x^l\bm{u}^\delta)_\Omega| \\
&\le \|\partial_x^l(\Phi_0^7\bm{u}^\delta)\|\|\partial_x^l\bm{u}^\delta_y\|
 +|(\partial_x^l(\Phi_0^7\bm{u}_y^{\delta}),\partial_x^l\bm{u}^{\delta})_\Gamma|
 +\|\partial_x^l(\Phi_{0y}^7\bm{u}^\delta_y)\|\|\partial_x^k\bm{u}^\delta\| \\
&\le \epsilon\tilde{E}_m
 +C_\epsilon\bigl(\|\partial_x^l(\Phi_0^7\bm{u}^\delta_y)\|^2+\|\partial_x^l(\Phi_{0y}^7\bm{u}^\delta_y)\|^2\bigr)
  +|(\partial_x^l(\Phi_0^7\bm{u}_y^{\delta}),\partial_x^l\bm{u}^{\delta})_\Gamma|.
\end{align*}
Here the estimates of 
$\|\partial_x^l(\Phi_0^7\bm{u}^\delta_y)\|^2$ and $\|\partial_x^l(\Phi_{0y}^7\bm{u}^\delta_y)\|^2$ 
are reduced to that of $\|\partial_x^l(\Phi_0^5\bm{u}^\delta_y)\|^2$. 
By Lemma \ref{NLlem}, we have 
$$
\|\partial_x^l(\Phi_0^5\bm{u}^\delta_y)\|^2
\lesssim \|\Phi_0^5\|_{L^{\infty}}^2\|\partial_x^l\bm{u}^\delta_y\|^2
 +(\|\partial_x^l\Phi_0^5\|^2+\|\partial_x^l\Phi_{0y}^5\|^2)(\|\bm{u}^\delta_y\|^2+\|\bm{u}^\delta_{xy}\|^2),
$$
which together with \eqref{NLlem2-1} and \eqref{NLlem2-5} in Lemma \ref{NLlem2} gives 
$\|\partial_x^l(\Phi_0^5\bm{u}^\delta_y)\|^2 \lesssim \tilde{E}_2\tilde{E}_m$. 
Concerning the boundary integral, by the first equation in \eqref{BC1} and the second equation in \eqref{NS}, 
we can replace $u_y$ and $\delta v_y$ by $h_1+(2+b_3)\eta-\delta^2v_x$ and by $-\delta u_x$, respectively, 
so that we obtain 
$$
|(\partial_x^l(\Phi_0^7\bm{u}_y^{\delta}),\partial_x^l\bm{u}^{\delta})_\Gamma|
 \lesssim |\Phi_0^7(2+b_3)\eta|_m|\bm{u}^{\delta}|_m
  +(|\Phi_0^7h_1|_{m-\frac12}+\delta|\Phi_0^7\bm{u}_x^{\delta}|_{m-\frac12})|\bm{u}^{\delta}|_{m+\frac12}.
$$
These terms can be treated by the estimate of $J_8$ and \eqref{fracu} and \eqref{gammaku} 
in Lemma \ref{additonal lemma}. 
Therefore, we obtain $K_{13} \le (\epsilon+C_\epsilon\tilde{E}_2)\tilde{E}_m$. 
Thus \eqref{est-fx4} holds. 

As for \eqref{est-fx3}, it suffices to estimate 
$$
\begin{cases}
K_{14}:=\|\partial_x^k(\Phi_0^5D_\delta^3\tilde{\eta})\|^2, \\
K_{15}:=\delta^{2i}\|\partial_x^k(\Phi_0^5\partial_x^i \partial_y^j \bm{u}^\delta)\|^2
 & \mbox{for}\quad 0\le i+j\le2,\ j\neq2.
\end{cases}
$$
By \eqref{interpolation}, \eqref{LI2} in Lemma \ref{LI}, \eqref{dxkE} in Lemma \ref{L2estimates}, 
and \eqref{NLlem2-1} and \eqref{NLlem2-3} in Lemma \ref{NLlem2}, we obtain $K_{14} \lesssim D_2D_m$. 
By Lemma \ref{NLlem}, we have 
$$
K_{15} \lesssim \|\Phi_0^5\|_{L^{\infty}}^2\delta^{2i}\|\partial_x^{k+i}\partial_y^j\bm{u}^\delta\|^2
 +(\|\partial_x^k\Phi_0^5\|^2+\|\partial_x^k\Phi_{0y}^5\|^2)
  \delta^{2i}(\|\partial_x^i\partial_y^j \bm{u}^\delta\|^2+\|\partial_x^i\partial_y^j \bm{u}^\delta_{x}\|^2),
$$
which together with \eqref{NLlem2-1} and \eqref{NLlem2-5} in Lemma \ref{NLlem2} gives 
$K_{15} \lesssim D_2D_m$. 
Thus \eqref{est-fx3} holds. 

As for \eqref{est-ft}, by the definition of $\bm{F}_2$ (see \eqref{F2}) and using the third equation in \eqref{BC1}, 
it suffices to estimate
$$
\begin{cases}
K_{16}:=\delta\|\partial_x^k(\Phi_0^5D_\delta^i\tilde{\eta}_t)\|^2 & \mbox{for}\quad i=1,2,3, \\
K_{17}:=\delta\|\partial_x^k(\Phi^5_0u_{ty})\|^2, \\
K_{18}:=\delta\|\partial_x^k(\Phi^5\bm{u}^\delta_t u_y)\|^2, \\
K_{19}:=\delta^{2i+1}\|\partial_x^k(\Phi^6_0\partial_x^i\bm{u}^\delta_{t})\|^2 & \mbox{for}\quad i=0,1, \\
K_{20}:=\delta^3\|\partial_x^k(\Phi^5D_\delta^i\tilde{\eta}_{tt}\partial_y^j u)\|^2
 & \mbox{for}\quad (i, j)=(0,0), (1, 0), (0, 1), \\
K_{21}:=\delta^3\|\partial_x^k(\Phi^5D_\delta^i\tilde{\eta}_{tt}\bm{u}_x^{\delta})\|^2, \\
K_{22}:=\delta^{i+2}|(\partial_x^k(\Phi^7_0\partial_x^{i+1}\partial_y^j\bm{u}^\delta_t),
 \partial_x^k\bm{u}^\delta_t)_\Omega| & \mbox{for}\quad  (i, j)=(1, 0), (0, 1), \\
K_{23}:=\delta|(\partial_x^k(\Phi_0^7 u_{yy})_{t}, \partial_x^k\bm{u}^\delta_{t})_\Omega|.
\end{cases}
$$
Here we did not list the terms which we have already estimated as $K_1,\ldots, K_5$. 
By \eqref{interpolation}, \eqref{LI-v} in Lemma \ref{LI}, \eqref{vtilde} in Lemma \ref{L2estimates}, 
and \eqref{NLlem2-1} and \eqref{NLlem2-3} in Lemma \ref{NLlem2}, 
we obtain $K_{16} \lesssim \tilde{E}_2\tilde{F}_m+\tilde{F}_2\tilde{E}_m$. 
By Lemma \ref{NLlem}, we have 
$$
K_{17} \lesssim \|\Phi^5_0\|^2_{L^{\infty}}\delta\|\partial_x^k u_{ty}\|^2
 +(\|\partial_x^k\Phi^5_0\|^2+\|\partial_x^k\Phi^5_{0y}\|^2)\delta(\|u_{ty}\|^2+\|u_{txy}\|^2), \\
$$
which together with \eqref{NLlem2-1} and \eqref{NLlem2-5} in Lemma \ref{NLlem2} gives 
$K_{17} \lesssim \tilde{E}_2\tilde{F}_m+\tilde{E}_m\tilde{F}_2$. 
By Lemma \ref{NLlem}, we have 
\begin{align*}
K_{18}
&\lesssim \|\Phi^5\|^2_{L^{\infty}}\delta\|\bm{u}^\delta_t\|^2_{L^{\infty}}\|\partial_x^k u_y\|^2 
 +\delta\|\bm{u}^\delta_t\|^2_{L^{\infty}}(\|\partial_x^k\Phi^5\|^2+\|\partial_x^k\Phi^5_y\|^2)
  (\|u_y\|^2+\|u_{xy}\|^2) \\
&\quad
 +\|\Phi^5\|^2_{L^{\infty}}\delta(\|\partial_x^k \bm{u}^\delta_t\|^2+\|\partial_x^k \bm{u}^\delta_{ty}\|^2)
  (\|u_y\|^2+\|u_{xy}\|^2),
\end{align*}
which together with the second inequality in \eqref{LI-u} in Lemma \ref{LI} and \eqref{NLlem2-1} 
and \eqref{NLlem2-5} in Lemma \ref{NLlem2} gives $K_{18} \lesssim \tilde{E}_m\tilde{F}_2+\tilde{E}_2\tilde{F}_m$. 
By \eqref{interpolation}, the second inequality in \eqref{LI-u} in Lemma \ref{LI} 
and \eqref{NLlem2-1} and \eqref{NLlem2-3} in Lemma \ref{NLlem2}, we obtain 
$K_{19}\lesssim \tilde{E}_2\tilde{F}_m+\tilde{F}_2\tilde{E}_m$. 
As for $K_{20}$, we will consider the case $(i, j)=(0, 1)$ only, because the other cases can be treated more easily. 
By Lemma \ref{NLlem}, we have 
\begin{align*}
K_{20}
&\lesssim \|\Phi^5\|_{L^{\infty}}^2\delta^3\|\tilde{\eta}_{tt}\|_{L^{\infty}}^2\|\partial_x^k u_y\|^2 
 +\delta^3\|\tilde{\eta}_{tt}\|_{L^{\infty}}^2
  (\|\partial_x^k\Phi^5\|^2+\|\partial_x^k\Phi^5\|^2)(\|u_y\|^2+\|u_{xy}\|^2) \\ 
&\quad
 +\|\Phi^5\|_{L^{\infty}}^2\delta^3(\|\partial_x^k\tilde{\eta}_{tt}\|^2+\|\partial_x^k\tilde{\eta}_{tty}\|^2)
  (\|u_y\|^2+\|u_{xy}\|^2), 
\end{align*}
which together with \eqref{est-vt} in Lemma \ref{LI} and \eqref{NLlem2-1} and \eqref{NLlem2-5} in Lemma \ref{NLlem2} 
gives $K_{20} \lesssim \tilde{E}_2\tilde{F}_m+\tilde{F}_2\tilde{E}_m$. 
Similarly, we obtain $K_{21} \lesssim \tilde{E}_2\tilde{F}_m+\tilde{F}_2\tilde{E}_m$. 
As for $K_{22}$, integration by parts in $x$ yields 
\begin{align*}
K_{22}
&= \delta^{i+2}|(\partial_x^k\big\{(\Phi^7_0\partial_x^i\partial_y^j\bm{u}^\delta_t)_x
 -\Phi^7_{0x}\partial_x^i\partial_y^j\bm{u}^\delta_t\big\}, \partial_x^k\bm{u}^\delta_t)_\Omega| \\
&\le \delta^{i+2}|(\partial_x^k(\Phi^7_0\partial_x^i\partial_y^j\bm{u}^\delta_t),
  \partial_x^k\bm{u}^\delta_{tx})_\Omega|
 +\delta^{i+2}|(\partial_x^k(\Phi^7_{0x}\partial_x^i\partial_y^j\bm{u}^\delta_t),
  \partial_x^k\bm{u}^\delta_t)_\Omega| \\
&\le \epsilon\tilde{F}_m
 +C_\epsilon\bigl(\delta^{2i+1}\|\partial_x^k(\Phi^7_0\partial_x^i\partial_y^j\bm{u}^\delta_t)\|^2
  +\delta^{2i+3}\|\partial_x^k(\Phi^7_{0x}\partial_x^i\partial_y^j\bm{u}^\delta_t)\|^2\bigr).
\end{align*}
Since the estimate of the right-hand side of the above inequality is reduced to those of $K_{17}$ and $K_{19}$, 
we obtain $K_{22} \le \epsilon \tilde{F}_m+C_\epsilon(\tilde{E}_2\tilde{F}_m+\tilde{F}_2\tilde{E}_m)$. 
As for $K_{23}$, integration by parts in $y$ yields 
\begin{align*}
K_{23}
&= \delta|(\partial_x^k\big\{(\Phi^7_0u_y)_{y}-\Phi^7_{0y}u_y\big\}_t, \partial_x^k\bm{u}^\delta_{t})_\Omega| \\
&\le \delta|(\partial_x^k(\Phi^7_0u_y)_{t}, \partial_x^k\bm{u}^\delta_{ty})_\Omega|
 +\delta|(\partial_x^k(\Phi^7_0u_y)_t, \partial_x^k \bm{u}^\delta_{t})_\Gamma|
 +\delta|(\partial_x^k(\Phi^7_{0y}u_y)_{t}, \partial_x^k\bm{u}^\delta_{t})_\Omega| \\
&\le \epsilon\tilde{F}_m
 +C_\epsilon\bigl(\delta\|\partial_x^k(\Phi^7_0u_y)_{t}\|^2+\delta\|\partial_x^k(\Phi^7_{0y}u_y)_t\|^2\bigr)
  +\delta|(\partial_x^k(\Phi^7_{0}u_y)_{t}, \partial_x^k\bm{u}^\delta_{t})_\Gamma|.
\end{align*}
Here we can reduce the estimate of 
$\delta\|\partial_x^k(\Phi^7_0u_y)_{t}\|^2+\delta\|\partial_x^k(\Phi^7_{0y}u_y)_t\|^2$ 
to those of $K_4$ and $K_{17}$. 
Moreover, by the first equation in \eqref{BC1}, we can estimate the term 
$\delta|(\partial_x^k(\Phi^7_{0}u_y)_{t}, \partial_x^k\bm{u}^\delta_{t})_\Gamma|$ 
in the same way as the proof of \eqref{est-ht} in Lemma \ref{est-h}. 
We thereby obtain $K_{23} \le \epsilon \tilde{F}_m+C_\epsilon(\tilde{E}_2\tilde{F}_m+\tilde{F}_2\tilde{E}_m)$. 
Thus \eqref{est-ft} holds. 

As for \eqref{est-commutator}, by the definition of $\bm{G}_k$ (see \eqref{commutator}) we see that 
\begin{align*}
\delta|(\bm{G}_k, \partial_x^k\bm{u}^\delta_t)_\Omega|
&\le \epsilon\tilde{F}_m
 +C_\epsilon \delta\|[\partial_x^k, A_5]\big\{(I+A_4)\nabla_\delta p_t+  A_{4t}\nabla_\delta p\big\}\|^2 \\
&\quad
 +C_\epsilon \delta\|[\partial_x^k,A_{5t}]\bm{u}_t^{\delta}\|^2
 +\delta|([\partial_x^k, A_5]\bm{F}_{3t}, \partial_x^k\bm{u}^\delta_t)_\Omega| \\
&=: \epsilon\tilde{F}_m+K_{24}+K_{25}+K_{26}.
\end{align*}
Here we can assume $k\ge1$. 
By the fact that $A_4$ and $A_5$ are of the form $\Phi_0^7$ (see \eqref{NLp} and \eqref{A4p}), 
Lemma \ref{commu}, \eqref{LI-E} in Lemma \ref{LI}, and \eqref{NLlem2-1} and \eqref{NLlem2-5} in Lemma \ref{NLlem2}, 
we obtain 
\begin{align*}
K_{24}
&\le C_\epsilon\big\{\tilde{E}_2(\delta\|\nabla_\delta\partial_x^{k-1}p_t\|^2
 +\|\nabla_\delta\partial_x^{k-1} p\|^2) \\
&\qquad
 +\tilde{E}_m(\delta\|\nabla_\delta p_t\|^2+\delta\|\nabla_\delta p_{tx}\|^2
  +\|\nabla_\delta p\|^2+\|\nabla_\delta p_x\|^2)\big\}, 
\end{align*}
which gives $K_{24} \le C_\epsilon(\tilde{E}_2\tilde{F}_m+\tilde{E}_m\tilde{F}_2)$. 
The estimate for $K_{25}$ is reduced to that of $K_{19}$. 
Taking into account the explicit form of $\bm{F}_3$ (see \eqref{def-F3}), 
we can estimate $K_{26}$ in the same way as the proof of \eqref{est-ft}. 
Therefore the proof is complete. 
$\quad\square$

\bigskip
By Lemmas \ref{est-h} and \ref{est-f}, for the nonlinear term $N_m$ defined by \eqref{def-N}, 
we obtain the following proposition.

\begin{prop}\label{est-NL}
For any $\epsilon>0$ there exists a positive constants $C_\epsilon$ such that the following estimate holds. 
$$
N_m\le \epsilon\tilde{F}_m
 +C_\epsilon\bigl(\tilde{E}_2\tilde{F}_m+\tilde{F}_2\tilde{E}_m+\varepsilon\sqrt{\tilde{E}_2}\tilde{E}_m\bigr).
$$
Moreover, if $\varepsilon\lesssim\delta$, then we have 
$$
N_m\le \epsilon\tilde{F}_m
 +C_\epsilon(\tilde{E}_2\tilde{F}_m+\tilde{F}_2\tilde{E}_m).
$$
\end{prop}

Finally, we estimate the terms appearing in the right-hand side of \eqref{p-elliptic-higher} 
and \eqref{p-elliptic-time} in Lemma \ref{est-p}. 
By the explicit form of $g$ (see \eqref{def-g}), this consists of the terms in the form 
$$
\begin{cases}
\Phi(\tilde{\eta}, D_\delta\tilde{\eta},y)\delta\bm{u}^\delta_x u_y, \\
\Phi(\tilde{\eta}, D_\delta\tilde{\eta}, D_\delta^2\tilde{\eta}, \bm{u}^\delta,y)w_\lambda\partial_y^j\bm{u}^\delta
 & \mbox{for}\quad \lambda=1,2, \ j=0,1, \\
\Phi_0(\tilde{\eta}, D_\delta\tilde{\eta}, D_\delta^2\tilde{\eta}, \bm{u}^\delta, \delta \bm{u}^\delta_x;y)
 \delta\bm{u}^\delta_x.
\end{cases}
$$

\begin{lem}\label{est-NLpx}
The following estimates hold. 
\begin{align}
& \delta\|\partial_x^k g_x\|^2+\delta\|\partial_x^k g_{0x}\|^2
 +\delta\|\nabla_\delta\partial_x^k g_{0x}\|^2 \label{est-NLpx1} \\
&\quad
 +\delta\|\partial_x^k(N_6\nabla_\delta p)_x\|^2+\delta^2||D_x|^{k+\frac12}\phi_x|_0^2\lesssim\tilde{F}_m
 +\tilde{F}_2\tilde{E}_m, \notag \\
& \|\partial_x^k g\|^2+\|\partial_x^k g_{0}\|^2+\|\nabla_\delta\partial_x^k g_0\|^2
 +\|\partial_x^k(N_6\nabla_\delta p)\|^2+\delta||D_x|^{k+\frac12}\phi|_0^2\label{est-NLpx2} \\
&\qquad
 \lesssim (1+D_2)D_m+\tilde{E}_2\|\nabla_\delta\partial_x^k p\|^2
  +\min\{\tilde{E}_m, D_m\}(\|\nabla_\delta p\|^2+\|\nabla_\delta p_x\|^2), \notag \\
& \delta\|\partial_x^{l-1} g_t\|^2+\delta\|\partial_x^{l-1} g_{0t}\|^2
 +\delta\|\nabla_\delta\partial_x^{l-1} g_{0t}\|^2 \label{est-NLpt} \\
&\quad
  +\delta\|\partial_x^{l-1}(N_6\nabla_\delta p)_t\|^2+\delta^2||D_x|^{l-\frac12}\phi_t|_0^2
 \lesssim\tilde{F}_m+\tilde{F}_2\tilde{E}_m. \notag
\end{align}
\end{lem}

\noindent
{\it Proof}. \ 
By Lemma \ref{NLlem}, we have 
\begin{align*}
\delta^{3}\|\partial_x^k(\Phi^5\bm{u}^\delta_{xx}u_y)\|^2
&\lesssim \|\Phi^5\|^2_{L^{\infty}}\delta^{3}\|\bm{u}^\delta_{xx}\|^2_{L^{\infty}}\|\partial_x^k u_y\|^2\\
&\quad
 +\delta^{3}\|\bm{u}^\delta_{xx}\|^2_{L^{\infty}}(\|\partial_x^k \Phi^5\|^2+\|\partial_x^k \Phi_y^5\|^2)
  (\|u_y\|^2+\|u_{xy}\|^2) \\
&\quad
 +\|\Phi^5\|^2_{L^{\infty}}\delta^{3}(\|\partial_x^k\bm{u}^\delta_{xx}\|^2+\|\partial_x^k \bm{u}^\delta_{xxy}\|^2)
  (\|u_y\|^2+\|u_{xy}\|^2),
\end{align*}
which together with the second inequality in \eqref{LI-u} in Lemma \ref{LI} and 
\eqref{NLlem2-1} and \eqref{NLlem2-5} in Lemma \ref{NLlem2}, we obtain 
$$
\delta^{3}\|\partial_x^k(\Phi^5\bm{u}^\delta_{xx}u_y)\|^2
\lesssim \tilde{F}_2\tilde{E}_m+\tilde{F}_m\tilde{E}_2. 
$$
It follows from \eqref{LI2} in Lemma \ref{LI} that 
$\|\bm{u}^\delta_x\|_{L^{\infty}}^2\delta^3\|\partial_x^k u_{xy}\|^2\lesssim \tilde{E}_2 \tilde{F}_m$. 
Therefore, in the same way as the above estimate, we obtain 
$$
\delta^3\|\partial_x^k(\Phi^5\bm{u}^\delta_x u_{xy})\|^2 \lesssim \tilde{F}_2\tilde{E}_m+\tilde{F}_m\tilde{E}_2.
$$
These together with the estimates of $K_3$, $K_4$, and $K_5$ yield 
$\delta\|\partial_x^k g_x\|^2\lesssim \tilde{F}_2\tilde{E}_m+\tilde{E}_2\tilde{F}_m$. 
It follows from the explicit form of $g_0$ (see \eqref{pressN}) that 
$\delta\|\partial_x^k g_{0x}\|^2+\delta\|\nabla_\delta \partial_x^k g_{0x}\|^2
\lesssim \tilde{F}_m+\tilde{F}_2\tilde{E}_m$, where we used the estimates for $K_1,\ldots, K_5$. 
By Lemma \ref{NLlem}, we have 
\begin{align*}
\delta\|\partial_x^{k+1}(N_6\nabla_\delta p)\|^2
&\lesssim \|N_6\|_{L^{\infty}}^2\delta\|\partial_x^{k+1}\nabla_\delta p\|^2 \\
&\quad
 +\delta^2(\|\partial_x^{k+1}N_6\|^2+\|\partial_x^{k+1}N_{6y}\|)
  \delta^{-1}(\|\nabla_\delta p\|^2+\|\nabla_\delta p_x\|^2).
\end{align*}
Since $N_6$ is the nonlinear part of $A_6$, which is defined by \eqref{press2}, we see that 
$N_6$ is of the form $\Phi_0^7$. 
Thus by \eqref{NLlem2-1} and \eqref{NLlem2-5} in Lemma \ref{NLlem2}, we obtain 
\begin{equation}\label{est-n6}
\delta\|\partial_x^{k+1}(N_6\nabla_\delta p)\|^2 \lesssim \tilde{E}_2\tilde{F}_m+\tilde{E}_m\tilde{F}_2. 
\end{equation}
The definition of $\phi$ (see \eqref{pressD}), Lemma \ref{trace}, and \eqref{est-hx2} in Lemma \ref{est-h} 
imply $\delta^2||D_x|^{k+\frac12}\phi_x|_0^2\lesssim \tilde{F}_m+\tilde{F}_2\tilde{E}_m$. 
Combining the above estimates, we obtain \eqref{est-NLpx1}. 

By Lemma \ref{NLlem}, \eqref{LI2} in Lemma \ref{LI} and \eqref{NLlem2-1}, and \eqref{NLlem2-5} in Lemma \ref{NLlem2}, 
we obtain 
\begin{align*}
\delta^2\|\partial_x^k(\Phi^5\bm{u}^\delta_x u_y)\|^2
&\lesssim \|\Phi^5\|^2_{L^{\infty}}\delta^{2}\|\bm{u}^\delta_{x}\|^2_{L^{\infty}}\|\partial_x^k u_y\|^2 \\
&\quad
 +\delta^{2}\|\bm{u}^\delta_{x}\|^2_{L^{\infty}}(\|\partial_x^k \Phi^5\|^2
  +\|\partial_x^k \Phi_y^5\|^2)(\|u_y\|^2+\|u_{xy}\|^2)\\
&\quad
 +\|\Phi^5\|^2_{L^{\infty}}\delta^{2}(\|\partial_x^k\bm{u}^\delta_{x}\|^2
  +\|\partial_x^k \bm{u}^\delta_{xy}\|^2)(\|u_y\|^2+\|u_{xy}\|^2) \\
&\lesssim D_2(1+D_m).
\end{align*}
By \eqref{interpolation}, \eqref{LI2} in Lemma \ref{LI}, 
and \eqref{NLlem2-1} and \eqref{NLlem2-3} in Lemma \ref{NLlem2}, we get 
$\delta^4\|\partial_x^k\bigl(\Phi^5(\bm{u}^\delta_x)^2\bigr)\|^2 \lesssim D_2(1+D_m)$. 
These together with the estimate of $K_{15}$ yield 
$\|\partial_x^k g\|^2\lesssim D_2(1+D_m)$. 
By the estimate for $K_{15}$, we obtain 
$\|\partial_x^k g_0\|^2+\|\nabla_\delta \partial_x^k g_{0}\|^2\lesssim (1+D_2)D_m$. 
In the same way as the proof of \eqref{est-n6}, we obtain 
$$
\|\partial_x^k(N_6\nabla_\delta p)\|^2
\lesssim \tilde{E}_2\|\nabla_\delta\partial_x^k p\|^2
 +\min\{\tilde{E}_m, D_m\}(\|\nabla_\delta p\|^2+\|\nabla_\delta p_x\|^2).
$$
Lemma \ref{trace} and \eqref{est-hx4} in Lemma \ref{est-h} lead to 
$\delta||D_x|^{k+\frac12}\phi|_0^2\lesssim (1+D_2)D_m$. 
Combining the above estimates implies \eqref{est-NLpx2}. 

By Lemma \ref{NLlem}, we have 
\begin{align*}
\delta^3\|\partial_x^{l-1}(\Phi^5\bm{u}^\delta_{tx} u_y)\|^2
&\lesssim \|\Phi^5\|_{L^{\infty}}^2\delta^3\|\bm{u}^\delta_{tx}\|_{L^{\infty}}^2\|\partial_x^{l-1}u_y\|^2 \\
&\quad
 +\delta^3\|\bm{u}^\delta_{tx}\|_{L^{\infty}}^2(\|\partial_x^{l-1}\Phi^5\|^2
  +\|\partial_x^{l-1}\Phi^5_y\|^2)(\|u_y\|^2+\|u_{xy}\|^2) \\
&\quad
 +\|\Phi^5\|_{L^{\infty}}^2(\|u_y\|^2+\|u_{xy}\|^2)
  \delta^3(\|\partial_x^{l}\bm{u}^\delta_t\|^2+\|\partial_x^{l}\bm{u}^\delta_{ty}\|),
\end{align*}
which together with the second inequality in \eqref{LI-u} in Lemma \ref{LI} 
and \eqref{NLlem2-1} and \eqref{NLlem2-5} in Lemma \ref{NLlem2} gives 
$\delta^3\|\partial_x^{l-1}(\Phi^5\bm{u}^\delta_{tx} u_y)\|^2
 \lesssim \tilde{E}_2\tilde{F}_m+\tilde{F}_2\tilde{E}_m$. 
In a similar way, we get 
$\delta^3\|\partial_x^{l-1}(\Phi^5\bm{u}^\delta_{t} u_{xy})\|^2
 \lesssim \tilde{E}_2\tilde{F}_m+\tilde{F}_2\tilde{E}_m$. 
Thus, in the same way as the proof of \eqref{est-NLpx1}, we obtain \eqref{est-NLpt}. 
The proof is complete. 
$\quad\square$

\section{Proof of the main theorem}
Summarizing the estimates in the last sections, we will prove the following proposition.

\begin{prop}\label{uniform}
Let $m$ be an integer satisfying $m\ge2$, $0<{\rm R}_1\le{\rm R}_0$, $0<{\rm W}_1\le{\rm W}_2$, and $0<\alpha\le\alpha_0$, where ${\rm R}_0$ and $\alpha_0$ are constants in Propositions \ref{highEest} and \ref{energy0}. There exist positive constants $c_1$, $C_5$, $C_6$, and $C_7$ such that if the solution $(\eta, u, v, p)$ of \eqref{NS}--\eqref{BC0} and the parameters $\delta$, $\varepsilon$, ${\rm R}$, and ${\rm W}$ satisfy 
$$
\tilde{E}_2(t)\le c_1,\quad0<\delta, \varepsilon\le1, \quad {\rm R_1}\le{\rm R}\le{\rm R_0}, \quad {\rm W}_1\le{\rm W}\le\delta^{-2}{\rm W}_2,
$$
then we have 
\begin{equation}\label{energy-exp}
\tilde{E}_2(t) \le C_7 E_2(0){\rm e}^{C_6\varepsilon t}, \quad
 \tilde{E}_m(t)+\int_0^t\tilde{F}_m(\tau){\rm d}\tau
  \le C_7 E_m(0)\exp(C_5 E_2(0){\rm e}^{C_6\varepsilon t}+C_5\varepsilon t)
\end{equation}
Moreover, if $\varepsilon\lesssim\delta$, then we have 
$$
\tilde{E}_2(t) \le C_7 E_2(0), \quad 
\tilde{E}_m(t)+\int_0^t\tilde{F}_m(\tau){\rm d}\tau \le C_7 E_m(0)\exp(C_5 E_2(0)).
$$
\end{prop}

In order to prove the above proposition, we prepare the following lemma.

\begin{lem}\label{est-E and F}
Under the same assumptions of Proposition \ref{uniform}, for any integer $k$ satisfying $0\le k\le m$, 
the following estimates hold. 
\begin{align}
&\tilde{E}_m\lesssim E_m, \label{est-tildeE} \\
&\tilde{F}_m\lesssim F_m+\tilde{F}_2\tilde{E}_m, \label{est-tildeF} \\
&\|(1+|D_x|)^m\nabla_\delta p\|^2\lesssim (1+D_2)^2 D_m. \label{est-pD}
\end{align}
\end{lem}

\noindent
{\it Proof}. \ 
As for \eqref{est-tildeE}, by the definition of $\tilde{E}_m$ (see \eqref{tildeE}) and 
Poincar$\acute{\rm e}$'s inequality, it suffices to show that for any $\epsilon>0$ 
there exists a positive constant $C_\epsilon$ such that 
\begin{equation}\label{est-uy}
\|\partial_x^k u_y\|^2 \le \epsilon \tilde{E}_m+C_\epsilon(E_m+\tilde{E}_2\tilde{E}_m).
\end{equation}
Applying $\partial_x^k$ to \eqref{NS}--\eqref{BC0} and using the argument 
in the proof of Proposition \ref{highEest}, we obtain 
\begin{align*}
\frac{1}{4K}\|\nabla_\delta\partial_x^k\bm{u}^\delta\|^2
&\le -\bigg\{{\rm R}\delta(\partial_x^k\bm{u}^\delta,\partial_x^k\bm{u}_t^\delta)_\Omega
 +2\bigg(\dfrac{1}{\tan\alpha}\delta(\partial_x^k\eta,\partial_x^k\eta_t)_\Gamma
  +\dfrac{\delta^2{\rm W}}{\sin\alpha}\delta(\partial_x^k\eta_{x},\partial_x^k\eta_{tx})_\Gamma\bigg)\bigg\} \\
&\quad
 +4K\bigl(|\partial_x^k\eta|_0^2+|\partial_x^k(b_3\eta)|_0^2\bigr)
  +(\partial_x^k h_1, \partial_x^k u)_\Gamma-2(\partial_x^k h_2, \delta\partial_x^k v)_\Gamma \\
&\quad
 +2(\dfrac{1}{\tan\alpha}\partial_x^k\eta-\dfrac{\delta^2{\rm W}}{\sin\alpha}\partial_x^k\eta_{xx},
  \delta \partial_x^k h_3)_\Gamma \\
&\quad
 +{\rm R}(\partial_x^k\bm{f}, \partial_x^k\bm{u}^\delta)_\Omega
 +(\partial_x^k\{-2A_4\nabla_\delta p+(b_2u_{yy}, 0)^{\rm T}\}, \partial_x^k \bm{u}^\delta)_\Omega.
\end{align*}
Here we consider the case $k\ge1$ only, because the case $k=0$ can be treated more easily. 
Then, by Lemma \ref{trace} we obtain 
\begin{align*}
\|\nabla_\delta\partial_x^k\bm{u}^\delta\|^2
&\lesssim E_m+|b_3\eta|_m^2+\delta^{-1}|(h_1,h_2)|_{m-\frac12}^2+\delta^2|h_3|_m^2
 +\delta^{-2}\|\partial_x^{k-1}\bm{f}\|^2 \\
&\quad
 +|(\partial_x^k\{-2A_4\nabla_\delta p+(b_2u_{yy}, 0)^{\rm T}\}, \partial_x^k \bm{u}^\delta)_\Omega|.
\end{align*}
It is easy to see that $|b_3\eta|_m^2+\delta^2|h_3|_m^2 \lesssim E_m$. 
Combining these, \eqref{est-hx3} in Lemma \ref{est-h}, and \eqref{est-fx2} 
and \eqref{est-fx4} in Lemma \ref{est-f}, we obtain \eqref{est-uy}. 
Then, taking $\epsilon$ and $c_1$ sufficiently small we get \eqref{est-tildeE}. 

As for \eqref{est-tildeF}, in view of the definition of $\tilde{F}_m$ (see \eqref{tildeF}), 
it suffices to show 
\begin{align}
& \delta^{-1}\|\nabla_\delta\partial_x^k p\|^2+\delta\|\nabla_\delta\partial_x^k p_x\|^2
 +\delta\|\nabla_\delta \partial_x^{l-1} p_t\|^2 \label{est-many} \\
&\qquad
  +\delta^6||D_x|^{k+\frac72}\eta|^2_0+\delta|(1+\delta|D_x|)^{\frac52}\partial_x^k\eta_t|_0^2
 \lesssim F_m+\tilde{E}_2\tilde{F}_m+\tilde{F}_2\tilde{E}_m. \notag
\end{align}
Combining Lemma \ref{est-p}, \eqref{est-hx1} in Lemma \ref{est-h}, \eqref{est-fx1} in Lemma \ref{est-f}, 
and \eqref{est-NLpx1} and \eqref{est-NLpt} in Lemma \ref{est-NLpx}, we obtain 
\begin{equation}\label{est-pp}
\delta^{-1}\|\nabla_\delta\partial_x^k p\|^2+\delta\|\nabla_\delta\partial_x^k p_x\|^2
  +\delta\|\nabla_\delta \partial_x^{l-1} p_t\|^2
 \lesssim F_m+\tilde{E}_2\tilde{F}_m+\tilde{F}_2\tilde{E}_m.
\end{equation}
We proceed to estimate $(\delta^2{\rm W})^2\delta^2||D_x|^{k+\frac72}\eta|^2_0$.
Applying $-\delta|D_x|^{k+\frac32}$ to the second equation in $\eqref{BC1}$ and 
taking the inner product of $(\delta^2{\rm W})\delta|D_x|^{k+\frac72}\eta$ with the resulting equality, we have 
\begin{align*}
& (\dfrac{1}{\tan\alpha}\delta|D_x|^{k+\frac32}\eta
 +\dfrac{\delta^2{\rm W}}{\sin\alpha}\delta|D_x|^{k+\frac72}\eta, (\delta^2{\rm W})\delta|D_x|^{k+\frac72}\eta)_\Gamma \\
&\qquad
 =(\delta|D_x|^{k+\frac32}(p-\delta v_y-h_2), (\delta^2{\rm W})\delta|D_x|^{k+\frac72}\eta)_\Gamma,
\end{align*}
which together with Lemma \ref{trace} and the second equation in \eqref{NS} leads to 
\begin{align*}
&(\delta^2{\rm W})^2\delta^2||D_x|^{k+\frac72}\eta|^2_0\\
&\qquad\lesssim \delta^2||D_x|^{\frac12}\partial_x^k(p_x+\delta u_{xx}- h_{2x})|_0^2 \\
&\qquad\lesssim \delta\|\partial_x^k p_x\|^2+\delta\|\partial_x^k\nabla_\delta p_{x}\|^2
 +\delta^3\|\partial_x^k u_{xx}\|^2+\delta^3\|\partial_x^k\nabla_\delta u_{xx}\|+\delta^2||D_x|^{k+\frac12}h_{2x}|_0^2.
\end{align*}
Combining this,  \eqref{est-hx2} in Lemma \ref{est-h}, and \eqref{est-pp}, 
we obtain the estimate for $(\delta^2{\rm W})^2\delta^2||D_x|^{k+\frac72}\eta|^2_0$. 
Finally, the estimate for $\delta|(1+\delta|D_x|)^{\frac52}\partial_x^k\eta_t|_0^2$ follows easily from 
the third equation in \eqref{BC1} and the estimate for $\delta^6||D_x|^{k+\frac72}\eta|^2_0$. 
Thus, we obtain \eqref{est-many}. 
Then, taking $c_1$ sufficiently small we get \eqref{est-tildeF}. 

As for \eqref{est-pD}, using \eqref{p-elliptic} and \eqref{est-NLpx2} in Lemma \ref{est-NLpx} 
and taking $c_1$ sufficiently small, we have 
$$
\|\nabla_\delta\partial_x^k p\|^2\lesssim (1+D_2) D_m+\min\{\tilde{E}_m, D_m\}(\|\nabla_\delta p\|^2+\|\nabla_\delta p_x\|^2).
$$
Considering the case $m=2$ and $k=0, 1$ in the above inequality and taking $c_1$ sufficiently small yield 
$\|\nabla_\delta p\|^2+\|\nabla_\delta p_x\|^2\lesssim (1+D_2)D_2$, 
which together with the above estimates gives \eqref{est-pD}. 
The proof is complete. 
$\quad\square$

\bigskip
\noindent
{\it Proof of Proposition \ref{uniform}}. \ 
Combining \eqref{est-energy}, Proposition \ref{est-NL}, and \eqref{est-tildeE} and \eqref{est-tildeF} 
in Lemma \ref{est-E and F} and taking $\epsilon$ and $c_1$ sufficiently small, we have 
\begin{equation}\label{est-energy2}
\frac{{\rm d}}{{\rm d}t}E_m(t)+\tilde{F}_m(t)\le C_5(\tilde{F}_2(t)+\varepsilon)E_m(t)
\end{equation}
for a positive constant $C_5$ independent of $\delta$. 
Note that if $\varepsilon\lesssim\delta$, 
then we can drop the term $C_5\varepsilon E_m(t)$ from the above inequality. 
Now, let us consider the case where $m=2$. 
By taking $c_1$ sufficiently small, we have 
$$
\frac{{\rm d}}{{\rm d}t}E_2(t)+\tilde{F}_2(t)\le C_6\varepsilon E_2(t)
$$
for a positive constant $C_6$ independent of $\delta$. 
Thus, Gronwall's inequality yields 
\begin{equation}\label{gronwall1}
E_2(t)+\int_0^t\exp\big(C_6\varepsilon(t-\tau)\big)\tilde{F}_2(\tau){\rm d}\tau
\le E_2(0){\rm e}^{C_6\varepsilon t}.
\end{equation}
In particular, we have $\int_0^t \tilde{F}_2(\tau){\rm d}\tau \le E_2(0){\rm e}^{C_6\varepsilon t}$. 
By this, \eqref{est-energy2}, and Gronwall's inequality, we see that 
\begin{align*}
E_m(t)+\int_0^t\tilde{F}_m(\tau){\rm d}\tau
&\le E_m(0)\exp\bigg(C_5\int_0^t(\tilde{F}_2(\tau)+\varepsilon){\rm d}\tau\bigg) \\
&\le E_m(0)\exp\big(C_5\tilde{E}_2(0){\rm e}^{C_6\varepsilon t}+C_5\varepsilon t\big).
\end{align*}
This together with \eqref{gronwall1} and \eqref{est-tildeE} in Lemma \ref{est-E and F} 
gives the desired estimates in Proposition \ref{uniform}. 
The proof is complete. 
$\quad\square$

\bigskip
\noindent
{\it Proof of Theorem \ref{main thm}}. \ 
Since the existence theorem of the solution locally in time is now classical, 
for example see \cite{Teramoto, Nishida}, 
it is sufficient to give a priori estimate of the solution. 
The first equation in \eqref{NS} leads to 
$$
\delta^2\|\partial_x^k\bm{u}^\delta_t\|^2
\lesssim \|\partial_x^k\bm{u}^\delta\|^2+\|\nabla_\delta\partial_x^k\bm{u}^\delta\|^2
 +\|\Delta_\delta\partial_x^k\bm{u}^\delta\|^2+\|\nabla_\delta\partial_x^k p\|^2+\|\partial_x^k\bm{f}\|^2.
$$
Thus, by \eqref{est-fx3} in Lemma \ref{est-f} and \eqref{est-pD} in Lemma \ref{est-E and F}, we have 
$\delta^2\|\partial_x^k\bm{u}^\delta_t\|^2 \lesssim (1+D_2)^2D_m$. 
By this, the third equation in \eqref{BC1}, and the definitions of $E_m$ and $D_m$ 
(see \eqref{E and F} and \eqref{Dk}), we obtain 
\begin{equation}\label{est-ut}
E_m(0) \le C_8\big(1+D_2(0)\big)^2 D_m(0)
\end{equation}
for a positive constant $C_8$ independent of $\delta$. 
Thus considering the case of $m=2$ in the above inequality, 
taking $D_2(0)$ and $T$ sufficiently small so that 
$2C_7C_8\big(1+D_2(0)\big)^2 D_2(0) \le c_1$ and $\quad {\rm e}^{C_6T}\le 2$, 
and using the first inequality in \eqref{energy-exp} in Proposition \ref{uniform}, 
we see that the solution satisfies  
$$
\tilde{E}_2(t) \le c_1 \qquad\mbox{for}\quad 0\le t\le T/\varepsilon.
$$
Thus, using the second inequality in \eqref{energy-exp} in Proposition \ref{uniform} 
together with \eqref{est-ut}, we obtain 
\begin{equation}\label{uniformE}
\tilde{E}_m(t)+\int_0^t\tilde{F}_m(\tau){\rm d}\tau \le C, 
\end{equation}
where the constant $C$ depends on ${\rm R}_1$, ${\rm W_1}$, ${\rm W}_2$, $\alpha$, and $M$ but not on $\delta$, $\varepsilon$, ${\rm R}$, nor ${\rm W}$. 
By the first equation in \eqref{NS}, we easily obtain 
$\delta^{-1}\|(1+|D_x|)^m(1+\delta|D_x|)u_{yy}\|^2\lesssim\tilde{F}_m$. 
Therefore, we obtain the desired estimate in Theorem \ref{main thm}. 
In view of the explicit form of $\tilde{E}_m$, using the second equation in \eqref{NS} and 
Poincar\'e's inequality, we easily obtain \eqref{uniform estimate}. 
Moreover, in the case where $\mathbb{G}=\mathbb{T}, \varepsilon\lesssim\delta$, and $\int_0^1\eta_0(x){\rm d}x=0$, it follows from Poincar\'e's inequality that $\delta E_m(t)\lesssim F_m(t)$, which yields \eqref{decay}.
The proof is complete. 
$\quad\square$


\bigskip
Hiroki Ueno \par
Department of Mathematics, Faculty of Science and Technology, Keio University, \par
3-14-1 Hiyoshi, Kohoku-ku, Yokohama 223-8522, Japan. \par
E-mail: hueno@math.keio.ac.jp

\bigskip
Akinori Shiraishi \par
Former address: Department of Mathematics, Faculty of Science and Technology,
\par Keio University, 
3-14-1 Hiyoshi, Kohoku-ku, Yokohama 223-8522, Japan. \par

\bigskip
Tatsuo Iguchi\par
Department of Mathematics, Faculty of Science and Technology, Keio University, \par
3-14-1 Hiyoshi, Kohoku-ku, Yokohama 223-8522, Japan. \par
E-mail: iguchi@math.keio.ac.jp


\begin{thebibliography}{99}
\bibitem{Atherton}
R. W. Atherton and G. M. Homsy, 
On the derivation of evolution equations for interfacial waves, 
Chem. Eng. Comm., {\bf 2} (1976), 57--77.
%
\bibitem{Beale}
J. T. Beale, 
Large-time regularity of viscous surface waves, 
Arch. Rational Mech. Anal., {\bf 84} (1984), 307--352.
%
\bibitem{Benjamin}
T. B. Benjamin, 
Wave formation in laminar flow down an inclined plane, 
J. Fluid Mech., {\bf 2} (1957), 554--574.
%
\bibitem{Benney}
D. J. Benney, 
Long waves on liquid film, 
J. Math. Phys., {\bf 45} (1966), 150--155.
%
\bibitem{Bresch}
D. Bresch and P. Noble, Mathematical justification of a shallow water model,  Methods Appl. Anal., {\bf 14} (2007), 87--117. 
%
\bibitem{Chang}
H. -C. Chang and E. A. Demekhin, 
Complex wave dynamics on thin films, 
Studies in Interface Science, 14, Elsevier Science B.V., Amsterdam, 2002.
%
\bibitem{Craster}
R. V. Craster and O. K. Matar, 
Dynamics and stability of thin liquid films, 
Rev. Mod. Phys., {\bf 81} (2009), 1131--1198. 
%
\bibitem{Gjevik}
B. Gjevik, 
Occurrence of finite-amplitude surface waves on falling liquid films, 
Phys. Fluids, {\bf 13} (1970), 1918--1925. 
%
\bibitem{Johnson}
R. S. Johnson, 
Shallow water waves on a viscous fluid---The undular bore, 
Phys. Fluids, {\bf 15} (1972), 1693--1699. 
%
\bibitem{Kalli}
S. Kalliadasis, C. Ruyer-Quil, B. Scheid, and M. G. Velarde, 
Falling Liquid film, 
Applied Mathematical Sciences, 176, Springer, London, 2012.
%
\bibitem{Kapitza}
P. L. Kapitza, 
Wave flow in thin layers of a viscous fluid, Zh. Eksp. Teor. Fiz., {\bf 19} (1948), 105--120 
(Collected Works of P. L. Kapitza, Pergamon, Oxford, 1965).
%
\bibitem{Lin}
S. P. Lin and M. V. G. Krishna, 
Stability of a liquid film with respect to initially finite three-dimensional disturbances, 
Phys. Fluids, {\bf 20} (1977), 2005--2011.
%
\bibitem{Mei}
C. C. Mei, 
Nonlinear gravity waves in a thin sheet of viscous fluid, 
J. Math. Phys., {\bf 45} (1966), 266--288.
%
\bibitem{Nakaya}
C. Nakaya, 
Long waves on a thin fluid layer flowing down an inclined plane, 
Phys. Fluids, {\bf 18} (1975), 1407--1412.
%
\bibitem{Nishida}
T. Nishida, Y. Teramoto, and H. A. Win, 
Navier--Stokes flow down an inclined plane: downward periodic motion, 
J. Math. Kyoto Univ., {\bf 33} (1993), 787--801.
%
\bibitem{Oron}
A. Oron, S. H. Davis, and S. G. Bankoff, 
Long-scale evolution of thin liquid films, 
Rev. Mod. Phys., {\bf 69} (1997), 931--980. 
%
\bibitem{Roskes}
G. J. Roskes, 
Three-dimensional long waves on a liquid film, 
Phys. Fluids, {\bf 13} (1970), 1440--1445. 
%
\bibitem{Teramoto}
Y. Teramoto, 
On the Navier--Stokes flow down an inclined plane, 
J. Math. Kyoto Univ., {\bf 32} (1992), 593--619.
%
\bibitem{Tomoeda}
Y. Teramoto and K. Tomoeda, 
Optimal Korn's inequality for solenoidal vector fields on a periodic slab, 
Proc. Japan Acad. Ser. A Math. Sci., {\bf 88} (2012), 168--172. 
%
\bibitem{Kawahara}
J. Topper and T. Kawahara, 
Approximate equations for long nonlinear waves on a viscous fluid, 
J. Phys. Soc. Japan, {\bf 44} (1978), 663--666.
%
\bibitem{Uecker}
H. Uecker, Self-similar decay of spatially localized perturbations of the Nusselt solution for the inclined film problem, Arch. Rational Mech. Anal., {\bf 184} (2007), 401--447. 
%
\bibitem{Ueno}
H. Ueno and T. Iguchi, A mathematical justification of a thin film approximation for the flow down an inclined plane, arXiv:1506.08489.
%
\bibitem{Yih1}
C. S. Yih, 
Stability of parallel laminar flow with a free surface, 
Proceedings of the Second U. S. National Congress of Applied Mechanics, Ann Arbor, 1954, pp. 623--628. 
American Society of Mechanical Engineers, New York, 1955. 
%
\bibitem{Yih2}
C. S. Yih, 
Stability of liquid flow down an inclined plane, 
Phys. Fluids, {\bf 6} (1963), 321--334.
\end{thebibliography}
\end{document}